\documentclass[reqno]{amsart}
\usepackage{amssymb,enumerate,mathrsfs}

\numberwithin{equation}{section}

\newtheorem{theorem}[equation]{Theorem}
\newtheorem{proposition}[equation]{Proposition}
\newtheorem{lemma}[equation]{Lemma}
\newtheorem{corollary}[equation]{Corollary}

\theoremstyle{definition}
\newtheorem{definition}[equation]{Definition}
\newtheorem{remark}[equation]{Remark}

\def\Vsp{\rule[-1ex]{0ex}{3.5ex}}
\DeclareMathOperator{\sym}{ \sigma\!\!\!\sigma}

\def\A{\mathcal A}
\def\B{\mathcal B}
\def\C{\mathbb C}

\def\Dom{\mathcal D}
\def\K{\mathcal K}
\def\L{\mathscr L}
\def\N{\mathbb N}
\def\R{\mathbb R}
\def\S{\mathscr S}
\def\Sing{\mathcal E}

\def\bT{\,{}^b\hspace{-0.5pt}T}
\def\bsym{\,{}^b\!\!\sym}

\def\csym{\,{}^c\!\sym}
\def\cT{\,{}^c T}
\def\cpi{\,{}^c\hspace{-1.5pt}\pi}
\def\cn{\mathrm{cn}}

\def\e{\mathrm e}
\def\embed{\hookrightarrow}
\def\eps{\varepsilon}
\def\id{I}

\def\open#1{\smash[t]{\overset{{}_{\circ}}{#1}{}}}
\def\s{\mathrm s}
\def\set#1{\{#1\}}
\def\vp{\varphi}
\def\minus{\backslash}

\DeclareMathOperator{\bgres}{bg-res}
\DeclareMathOperator{\bgspec}{bg-spec}
\DeclareMathOperator{\Ind}{ind}

\DeclareMathOperator{\Diff}{Diff}
\DeclareMathOperator{\rg}{rg}

\DeclareMathOperator{\spec}{spec}
\DeclareMathOperator{\supp}{supp}

\begin{document}
\title{Resolvents of elliptic cone operators}
\author{Juan B. Gil}
\address{Penn State Altoona\\ 3000 Ivyside Park \\
 Altoona, PA 16601-3760}
\author{Thomas Krainer}
\address{Institut f\"ur Mathematik\\ Universit\"at Potsdam\\
 14415 Potsdam, Germany}
\author{Gerardo A. Mendoza}
\address{Department of Mathematics\\ Temple University\\
 Philadelphia, PA 19122}
\begin{abstract}
We prove the existence of sectors of minimal growth for general closed
extensions of elliptic cone operators under natural ellipticity conditions.
This is achieved by the construction of a suitable parametrix and reduction 
to the boundary. Special attention is devoted to the clarification of the
analytic structure of the resolvent.
\end{abstract}


\maketitle

\section{Introduction}

Motivated by Seeley's seminal work \cite{Seeley}, and with the same
intentions, the purpose of this paper is, first, to prove the
existence of sectors of minimal growth for general closed extensions
of elliptic cone differential operators under suitable ray
conditions on the symbols of the operator; and second, to describe
the structure of the resolvent as a pseudodifferential operator.

Previous relevant investigations in this direction assume that the
coefficients are constant near the boundary, cf. \cite{Mooers},
\cite{SchrSe03}, or the technically convenient but rather
restrictive dilation invariance of the domain, cf. \cite{BrSe91},
\cite{Le97}, \cite{GiHeat01}, \cite{LoRes01}, \cite{SchrSe03}.
Some of these works deal with special classes of operators such as
Laplacians.  In the general setting followed in this paper, the
interactions of lower order terms in the Taylor expansion of the
coefficients of the operators near the boundary lead to a domain
structure beyond the minimal domain $\Dom_{\min}$ that brings up essential 
new difficulties not present in the constant coefficients case. Thus the 
investigation of the general case entails the development of new techniques.

Let $M$ be a smooth compact $n$-manifold with boundary. Recall that
a cone differential operator is a linear differential operator with
smooth coefficients in the interior of $M$ which locally near the
boundary and in terms of coordinates $x,y_1,\dotsc,y_{n-1}$ with
$x=0$ on $\partial M$, is of the form
\begin{equation*}
x ^{-m}\sum_{k+|\alpha|\leq m} a_{k\alpha}(x,y)D_y^\alpha(xD_x)^k
\end{equation*}
with $a_{k\alpha}$ smooth up to the boundary and $m$ a positive integer. Such an operator is called $c$-elliptic if it is elliptic in the interior in the usual sense, and near the boundary, if written as above, then
\begin{equation*}
\sum_{k+|\alpha|= m} a_{k\alpha}(x,y)\eta^\alpha \xi^k
\end{equation*}
is an elliptic symbol up to $\set{x=0}$. Fix some smooth defining function $x$ for $\partial M$ with $x>0$ in the interior $\open M$ of $M$ and denote by $x^{-m}\Diff_b^m(M;E)$ the space of cone operators of order at most $m$ 
acting on sections of a Hermitian vector bundle $E\to M$. 

Cone differential operators arise when introducing polar coordinates
around a point, and for that reason they are of great interest in
the study of operators on manifolds with conical singularities  
(cf. \cite{K1}, \cite{SzNorthHolland}). In this context it is natural to base the $L^2$ theory of these operators, at least initially, on a $c$-density on $M$, which is a measure of the form $x^n\mathfrak m$ where $\mathfrak m$ is a smooth $b$-density, that is, $x \mathfrak m$ is a smooth everywhere positive density on $M$.

Let $A\in x^{-m}\Diff_b^m(M;E)$, and write $L^2_c(M;E)$ for the space
$L^2(M,x^n\mathfrak m;E)$.  There are two canonical closed extensions one can specify for the unbounded operator 
\begin{equation}\label{InitialDomain}
	A:C_0^\infty( \open M;E)\subset L^2_c(M;E)\to L^2_c(M;E),
\end{equation}
namely the closure
\begin{equation}\label{Amin}
A:\Dom_{\min}\subset L^2_c(M;E)\to L^2_c(M;E),
\end{equation}
and
\begin{equation}\label{Amax}
A:\Dom_{\max}\subset L^2_c(M;E)\to L^2_c(M;E),
\end{equation}
with
\begin{equation*}
\Dom_{\max}=\set{u\in L^2_c(M;E) : Au\in L^2_c(M;E)}.
\end{equation*}
Obviously, both $\Dom_{\min}$ and $\Dom_{\max}$ are complete in the graph norm,
\begin{equation*}
\|u\|_A=\|u\|_{L^2_c}+\|Au\|_{L^2_c},
\end{equation*}
and $\Dom_{\min}\subset\Dom_{\max}$.

Suppose that $A$ is $c$-elliptic. By a theorem of Lesch \cite{Le97},
$\Dom_{\min}$ has finite codimension in $\Dom_{\max}$, and all closed
extensions of \eqref{InitialDomain} are Fredholm and have domain $\Dom$ such
that $\Dom_{\min}\subset\Dom\subset\Dom_{\max}$. Moreover, if $A_\Dom$ denotes the closed extension with domain $\Dom$, then
\begin{equation}\label{LeschRelIndex} 
\Ind(A_\Dom)=\Ind(A_{\Dom_{\min}})+\dim	(\Dom/\Dom_{\min}), 
\end{equation}
see Lesch, op. cit. and Gil-Mendoza \cite{GiMe01}. Thus, if
$\Ind(A_{\Dom_{\min}})$ is already positive, then there is no extension of $A$
with nonempty resolvent set. In fact, a necessary and sufficient condition for the existence of a closed extension $A_\Dom$ of \eqref{InitialDomain} with nonempty resolvent set is that for some $\lambda\in\C$,
$A_{\Dom_{\min}}-\lambda$ is injective and $A_{\Dom_{\max}}-\lambda$ is
surjective, see \cite{GilKrainerMendoza1}.

Given a closed extension $A_\Dom$, we will prove in
Section~\ref{sec-Resolvents} (see Theorem~\ref{RayMinimalGrowth}) that under natural ellipticity conditions pertaining the symbol of $A$ and the model operator $A_\wedge$, cf. \eqref{Awedge}, there exists a sector
$$
\Lambda = \{z \in \C : z = re^{i\theta},\;r \geq 0,\; 
|\theta - \theta_0| \leq a\}
$$
of minimal growth, i.e.,
$$
A - \lambda : \Dom \to L_c^2(M;E)
$$
is invertible for $\lambda \in \Lambda$ with $|\lambda|$ large, and
$$
\|(A_\Dom - \lambda)^{-1}\|_{{\mathscr L}(L_c^2(M;E))} = O(|\lambda|^{-1})
\; \text{ as } |\lambda| \to \infty.
$$
More precisely, we require that $\Lambda$ is free of spectrum of the homogeneous principal $c$-symbol $\csym(A)$ of $A$ on $\cT^*M\minus\{0\}$, and that the model operator 
$$
A_{\wedge}-\lambda: \Dom_{\wedge} \to L^2_c(Y^{\wedge};\pi_Y^*E|_Y)
$$
is invertible for large $\lambda \in \Lambda$ with inverse bounded in the norm as $|\lambda| \to \infty$, where $Y^{\wedge} = \overline{\R}_+ \times Y$ is the stretched model cone with boundary $Y = \partial M$, and
$\Dom_{\wedge}$ is a domain for $A_{\wedge}$ associated with $\Dom$ in 
a natural way (see \cite{GilKrainerMendoza1}). 

The proof of this result relies on the construction of a  parameter-dependent parametrix
\begin{equation}\label{Linksparametrix}
B(\lambda): L_c^2(M;E) \to \Dom_{\min}(A),
\end{equation}
which is a left-inverse for the operator $A_{\Dom_{\min}} - \lambda$ for large $|\lambda|$. Then, in order to deal with the finite dimensional contribution of the domain $\Dom$ beyond $\Dom_{\min}$, we follow the idea of reduction to the boundary motivated by the point of view that the choice of a domain corresponds to the choice of a boundary condition for the operator $A$.

More precisely, we add a suitable operator family $K(\lambda)$ to $A_{\Dom{\min}} - \lambda$ such that
\begin{equation}\label{DirichletErsatz}
\begin{pmatrix} A_{\Dom_{\min}} - \lambda & K(\lambda) \end{pmatrix} :
\begin{array}{c} \Dom_{\min}(A) \\ \oplus \\ \C^{d''} \end{array} \to L_c^2(M;E)
\end{equation}
is invertible for large $|\lambda|$, and consider \eqref{DirichletErsatz} a ``Dirichlet problem'' for the operator $A - \lambda$.  Following Schulze's
viewpoint from the pseudodifferential edge-calculus \cite{Sz89, SzNorthHolland} we invert \eqref{DirichletErsatz} in the context of operator matrices by adding generalized Green remainders to the parametrix $B(\lambda)$. We then multiply the inverse
$\begin{pmatrix} B(\lambda) \\ T(\lambda) \end{pmatrix}$ from the left to the operator $A_{\Dom} - \lambda$, reducing the problem of inverting $A_\Dom - \lambda$ to the simpler problem of inverting the operator family
\begin{equation}\label{Randreduktion}
F(\lambda)= T(\lambda)(A-\lambda) : \Dom/\Dom_{\min} \to \C^{d''}.
\end{equation}
The operator $F(\lambda)$ can be interpreted as the reduction to the boundary of $A -\lambda$ under the boundary condition $\Dom$ by \eqref{DirichletErsatz}, and it plays a similar role as, e.g., the Dirichlet-to-Neumann map in classical boundary value problems.

It turns out that we may write the resolvent as
$$
(A_\Dom - \lambda)^{-1} = B(\lambda) + (A_\Dom - \lambda)^{-1}\Pi(\lambda)
$$
with $B(\lambda)$ from \eqref{Linksparametrix} and a finite
dimensional smoothing pseudodifferential projection $\Pi(\lambda)$ onto a
complement of the range of $A_{\Dom_{\min}} - \lambda$ in $L_c^2(M;E)$. The operators $B(\lambda)$ and $\Pi(\lambda)$ have complete asymptotic expansions as $|\lambda| \to \infty$ into homogeneous components in the interior and $\kappa$-homogeneous operator-valued components near the boundary, respectively. 

The structure of the paper is as follows: In Section
\ref{sec-Preliminaries} we recall basic facts about cone
operators and their symbols. Section \ref{sec-Domains} is devoted to
closed extensions in $L^2$ and in higher order Sobolev spaces. 
Section \ref{sec-RayConditions} concerns some relations between $A$ 
and its symbols regarding the discreteness of the spectrum and the
existence of sectors of minimal growth. In Section \ref{sec-Parametrix} we perform the construction of the parametrix \eqref{Linksparametrix} and establish the ``Dirichlet problem'' \eqref{DirichletErsatz}. Finally, in Section \ref{sec-Resolvents}, we prove the results about the existence and norm estimates of the resolvent by investigating the operator \eqref{Randreduktion}.

\section{Preliminaries}
\label{sec-Preliminaries}

Let $M$ be a smooth compact $n$-manifold with boundary and fix a defining function $x$ for $\partial M$ with $x>0$ in $\open M$.
If $E\to M$ is a complex vector bundle and $\Diff^m(M;E)$ is the space of differential operators on $C^\infty(M;E)$ of order $m$, then $\Diff^m_b(M;E)$ denotes the subspace consisting of totally characteristic differential operators on $C^\infty(M;E)$ of order $m$.

The elements of $x^{-m}\Diff^m_b(M;E)$, that is, differential
operators of the form $A=x^{-m}P$ with $P\in \Diff^m_b(M;E)$, are the
cone operators of order $m$.

According to \cite{GilKrainerMendoza1} we associate with $A$ an invariantly defined $c$-symbol
\begin{equation*}
\csym(A) \in C^\infty(\cT^*M\minus 0;\mathrm{End}(\cpi^*E))
\end{equation*}
on the $c$-cotangent bundle $\cT^*M \to M$, where $\cpi : \cT^*M \to M$ is the canonical projection map. Recall that $\cT^*M$ is the smooth vector bundle over $M$ whose space of smooth sections is
\begin{equation*}
C^\infty_{\cn}(M;T^*M)=\set{\eta\in C^\infty(M,T^*M): \iota^*\eta=0},
\end{equation*}
the space of $1$-forms on $M$ which are, over $\partial M$,
sections of the conormal bundle of $\partial M$ in $M$.

Let $\mathbf x^{-1}:\cT^*M\to\bT^*M$ be the natural isomorphism that is induced by the defining function $x$. Then the $c$-symbol of $A$ and the $b$-symbol of $x^mA$ are related as
\begin{equation*}
\csym(A)(\eta) =\bsym(x^mA)(\mathbf x^{-1}(\eta)).
\end{equation*}

\begin{definition}\label{ConeParameterElliptic}
The operator $A \in x^{-m}\Diff^m_b(M;E)$ is called $c$-elliptic if
\begin{equation*}
\csym(A) \in C^\infty(\cT^*M\minus 0;\mathrm{End}(\cpi^*E))
\end{equation*}
is an isomorphism. The family $\lambda\mapsto A-\lambda$ is called
$c$-elliptic with parameter in a set $\Lambda\subset \C$ if
\begin{equation*}
\csym(A)-\lambda \in C^\infty((\cT^*M\times \Lambda)\minus 0;\mathrm{End}((\cpi\times\mathit{id})^*E))
\end{equation*}
is an isomorphism. Here $\cpi\times\mathit{id}:(\cT^*M\times\Lambda)\minus 0\to M\times\Lambda$ is the canonical map.
\end{definition}

Let $E \to M$ be Hermitian, and $\mathfrak m$ be a positive $b$-density.
Recall that the Hilbert space $L^2_b(M;E)$ is the $L^2$ space of
sections of $E$ with respect to the Hermitian form on $E$ and the
density $\mathfrak m$. Thus the inner product is
\begin{equation*}
(u,v)_{L^2_b}=\int (u,v)_E\, \mathfrak m
\quad\text{ if }u,\ v \in L^2_b(M;E).
\end{equation*}
For non-negative integers $s$ the Sobolev spaces $H^s_b(M;E)$ are defined as
\begin{equation*}
 H^s_b(M;E)=\set{u\in L^2_b(M;E): Pu\in L^2_b(M;E)\
 \forall P\in \Diff^s_b(M;E)}.
\end{equation*}
The spaces $H^s_b(M;E)$ for general $s\in\R$ are
defined by interpolation and duality, and we set
\begin{equation*}
H^\infty_b(M;E)=\bigcap_s H^s_b(M;E),\quad
H^{-\infty}_b(M;E)=\bigcup_s H^s_b(M;E).
\end{equation*} The weighted spaces
\begin{equation*}
x^\mu H^s_b(M;E)=\set{x^\mu u: u\in H^s_b(M;E)}
\end{equation*}
are topologized so that $H^s_b(M;E)\ni u\mapsto x^\mu u\in x^\mu
H^s_b(M;E)$ is an isomorphism. In the case of $s=0$ one has
\begin{equation*}
x^\mu H^0_b(M;E)=x^\mu L^2_b(M;E)=L^2(M,x^{-2\mu}\mathfrak m;E),
\end{equation*}
and the Sobolev space based on $L^2(M,x^{-2\mu}\mathfrak m;E)$ and
$\Diff^s_b(M;E)$ is isomorphic to $x^\mu H^s_b(M;E)$.

To define a Mellin transform consistent with the density $\mathfrak
m$, pick a collar neighborhood $U_Y \cong Y\times[0,1)$ of the boundary $Y = \partial M$ in $M$,
and a defining function $x:M\to\R$ such that
\begin{equation}\label{ProductMeasure}
\mathfrak m=\frac{dx}{x}\otimes \pi_Y^*\mathfrak m_Y\text{ in } U_Y
\end{equation}
for some smooth density $\mathfrak m_Y$ on $\partial M$.
Let $\partial_x$ be the vector field tangent to the fibers of $U_Y\to Y$ such that
$\partial_x x=1$.

Fix $\omega\in C_0^\infty(-1,1)$ real valued, nonnegative and such
that $\omega=1$ in a neighborhood of $0$. Also fix a Hermitian
connection $\nabla$ on $E$. The Mellin transform of an element
$u\in C_0^\infty(\open M;E)$ is defined to be the entire function
$\hat u:\C\to C^\infty(Y;E|_Y)$ such that for any $v\in C^\infty(Y;E|_Y)$
\begin{equation*}
 (x^{-i\sigma}\omega u,\pi_Y^*v)_{L^2_b(M;E)}=
 (\hat u(\sigma),v)_{L^2(Y;E|_Y)}
\end{equation*}
By $\pi_Y^*v$ we mean the section of $E$ over $U_Y$ obtained by
parallel transport of $v$ along the fibers of $\pi_Y$. The Mellin transform thus defined extends to the spaces $x^\mu H^s_b(M;E)$ so as to give holomorphic functions on $\set{\sigma: \Im\sigma>-\mu}$ with values in $H^s(Y;E|_Y)$.
As is well known, the Mellin transform extends to the spaces $x^\mu
L^2_b(M;E)$ in such a way that if $u \in x^\mu L^2_b(M;E)$ then $\hat
u(\sigma)$ is holomorphic in $\set {\Im \sigma > -\mu}$ and in
$L^2(\set{\Im\sigma=-\mu}\times Y)$ with respect to $d\sigma\otimes
\mathfrak m_Y$.

Let $A = x^{-m}P$ with $P \in \Diff^m_b(M;E)$, and let
\begin{equation}\label{ConormalSymbol}
\C\ni \sigma\mapsto \hat P(\sigma) \in \Diff^m(Y;E|_Y)
\end{equation}
be the conormal symbol of $P$. Recall that $\hat{P}(\sigma)$ is elliptic for 
every $\sigma \in \C$ if $A$ is $c$-elliptic. The boundary spectrum of $A$ is
\begin{equation*}
\spec_b(A)=\set{\sigma\in \C: \hat P(\sigma)\text{ is not invertible}},
\end{equation*}
which is discrete if $A$ is $c$-elliptic,
and the conormal symbol of $A$ is defined to be that of the operator $P$.

Near $Y$ one can write
\begin{equation*}
P = \sum_{\ell=0}^m P_\ell' \circ (\nabla_{xD_x})^\ell
\end{equation*}
where the $P_\ell'$ are differential operators of order $m-\ell$
(defined on $U_Y$) such that for any smooth function $\phi(x)$ and
section $u$ of $E$ over $U_Y$, $P_\ell'(\phi(x) u) =
\phi(x)P_\ell'(u)$, in other words, of order zero in
$\nabla_{xD_x}$.
\begin{definition}\label{ConstantCoeff}
$P$ is said to have coefficients independent of $x$ near $Y$, 
or simply constant coefficients, if
\begin{equation*}
\nabla_{x\partial_x}P_k(u) = P_k(\nabla_{x\partial_x}u)
\end{equation*}
for any smooth section $u$ of $E$ supported in $U_Y$.
Correspondingly, $A$ is said to have coefficients independent of $x$ near $Y$ if
this holds for $P$.
\end{definition}
For any $N$ there are operators
$P_k$, $\tilde P_N\in \Diff^m_b(M;E)$ such that
\begin{equation}\label{bOperatorTaylor}
P = \sum_{k=0}^{N-1} P_k x^k + x^N \tilde P_N 
\end{equation}
where each $P_k$ has coefficients independent of $x$ near $Y$. If
$P_k$ has coefficients independent of $x$ near $Y$ then so does
its formal adjoint $P_k^\star$.

With $A=x^{-m}P\in x^{-m}\Diff^m_b(M;E)$ we associate on the model cone $Y^\wedge = \overline \R_+\times Y$ the operator
\begin{equation}\label{Awedge}
 A_\wedge = x^{-m} P_0,
\end{equation}
where $P_0\in \Diff^m_b(Y^\wedge;E)$ is the constant term in the
expansion \eqref{bOperatorTaylor} and has therefore coefficients
independent of $x$.

For $\varrho>0$ we consider the normalized dilation group action
from sections of $E$ to sections of $E$ on $\open Y^{\wedge}$ defined by
\begin{equation}\label{DilationGroup}
(\kappa_\varrho u)(x,y) = \varrho^{m/2}u(\varrho x,y).
\end{equation}
The normalizing factor $\varrho^{m/2}$ in the definition of $\kappa_\varrho$ is added only because it makes
\begin{equation*}
 \kappa_{\varrho}:x^{-m/2}L^2_b(Y^\wedge;E)\to x^{-m/2}L^2_b(Y^\wedge;E)
\end{equation*}
an isometry, where the measure on $L^2_b$ refers to the $b$-density
$\mathfrak m=\frac{dx}{x}\otimes \mathfrak m_Y$ on $Y^{\wedge}$.

Let $A^\star$ denote the formal adjoint of $A$ acting on $x^{-m/2}L^2_b(M;E)$. Then we have $(A_\wedge)^\star=(A^\star)_\wedge$.

The family $\lambda\mapsto A_\wedge -\lambda$ satisfies the homogeneity relation
\begin{equation}\label{HomogeneityOfA}
 A_\wedge - \varrho^m\lambda =
 \varrho^m \kappa_\varrho (A_\wedge -\lambda)\kappa_\varrho^{-1}
 \;\text{ for every } \varrho>0.
\end{equation}
\begin{definition}\label{kappaHomogeneous}
A family of operators $A(\lambda)$
acting on a $\kappa$-invariant space of distributions on $Y^\wedge$ will be called $\kappa$-homogeneous of degree $\nu$ if
\begin{equation*}
 A(\varrho^m\lambda) =
 \varrho^\nu \kappa_\varrho A(\lambda)\kappa_\varrho^{-1}
\end{equation*}
for every $\varrho>0$.
\end{definition}

This notion of homogeneity is systematically used in Schulze's
edge-calculus.

On $Y^\wedge$ it is convenient to introduce weighted Sobolev spaces
with a particular structure at infinity consistent with the structure
of the operators involved.
Let $\omega\in C^\infty_0(\R)$ be a nonnegative function with
$\omega(r)=1$ near $r=0$. We follow Schulze (cf. \cite{Sz89}) and
consider the space $H_{\rm cone}^s({Y^\wedge};E)$ consisting of
distributions $u$ such that given any coordinate patch $\Omega$ on
$Y$ diffeomorphic to an open subset of the sphere $S^{n-1}$, and
given any function $\vp \in C^\infty_0(\Omega)$, we have
$(1-\omega)\vp\,u \in H^s(\R^n;E)$ where $\R_+ \times S^{n-1}$ is
identified with $\R^n \minus \{0\}$ via polar coordinates.

For $s,\alpha\in\R$ we define $\mathcal{K}^{s,\alpha}(Y^\wedge;E)$ as
the space of distributions $u$ such that
\begin{equation*}
\omega u\in x^{\alpha}H^s_b(Y^\wedge;E) \text{ and }
(1-\omega) u\in x^{\frac{n-m}{2}}H_{\rm cone}^s({Y^\wedge};E)
\end{equation*}
for any cut-off function $\omega$. Note that
$H_{\rm cone}^0(Y^\wedge;E) = x^{-n/2}L^2_b(Y^\wedge;E)$.

It turns out that $C_0^\infty(\open Y^\wedge;E)$ is dense in
$\mathcal{K}^{s,\alpha}(Y^\wedge;E)$, and
\begin{equation}\label{BoundedAwedge}
 A_\wedge:\mathcal{K}^{s,\alpha}(Y^\wedge;E)
  \to \mathcal{K}^{s-m,\alpha-m}(Y^\wedge;E)
\end{equation}
is bounded for every $s$ and $\alpha$. The group
$\{\kappa_\varrho\}_{\varrho\in\R_+}$ is a strongly continuous
group of isomorphisms on $\mathcal{K}^{s,\alpha}$ for every
$s,\alpha\in\R$. As pointed out already, it defines an isometry
on the space
\begin{equation*}
\mathcal{K}^{0,-m/2}(Y^\wedge;E)=x^{-m/2}L^2_b(Y^\wedge;E)
\end{equation*}
which we will take as reference Hilbert space on $Y^\wedge$.

\section{Closed extensions}
\label{sec-Domains}

If $A\in x^{-m}\Diff^m_b(M;E)$, then for any $s$ and $\mu$,
\begin{equation*}
A:x^{\mu}H^{s}_b(M;E)\to x^{\mu-m}H^{s-m}_b(M;E)
\end{equation*}
is continuous. In order not to have to deal with the index $\mu$ we
normalize so that if our original interests are in $x^\mu
L^2_b(M;E)$, then we work with the operator
$x^{-\mu-m/2}Ax^{\mu+m/2}$ and base all the analysis on
$x^{-m/2}L^2_b(M;E)$. Clearly, $x^{-\mu-m/2}Ax^{\mu+m/2}\in
x^{-m}\Diff^m_b(M;E)$. The latter operator has the same $c$-symbol
as $A$, so is $c$-elliptic if and only if $A$ is so, and it has the
same spectral properties. This said, we assume that $\mu=-m/2$.

The closed extensions of elliptic cone operators on $x^{-m/2}L^2_b(M;E)$ have been studied by Lesch \cite{Le97} and by two of the authors of the present work in \cite{GiMe01}, among others. It is important for our purposes to admit arbitrary regularity. In analogy with the $x^{-m/2}L_b^2$-case,
two canonical closed extensions of the operator
\begin{equation}\label{InitialSobolevDomain}
A:C_0^\infty(\open M;E)\subset x^{-m/2}H^s_b(M;E)\to x^{-m/2}H^s_b(M;E),
\end{equation}
are singled out. Its closure
\begin{equation}\label{AminSobolev}
A:\Dom_{\min}^s(A)\subset x^{-m/2}H^s_b(M;E)\to x^{-m/2}H^s_b(M;E),
\end{equation}
and
\begin{equation}\label{AmaxSobolev}
A:\Dom_{\max}^s(A)\subset x^{-m/2}H^s_b(M;E)\to x^{-m/2}H^s_b(M;E),
\end{equation}
with
\begin{equation*}
\Dom_{\max}^s(A)=\set{u\in x^{-m/2}H^s_b(M;E) : Au\in x^{-m/2}H^s_b(M;E)}.
\end{equation*}
Both $\Dom_{\min}^s(A)$ and $\Dom_{\max}^s(A)$ are complete in the graph norm
\begin{equation}\label{GraphNorm}
\|u\|_{A,s}=\|u\|_{x^{-m/2}H^s_b}+\|Au\|_{x^{-m/2}H^s_b},
\end{equation}
and therefore $\Dom_{\min}^s(A)\subset\Dom_{\max}^s(A)$ is a
closed subspace. Clearly, for any closed extension
\begin{equation*}
A:\Dom\subset x^{-m/2}H^s_b(M;E)\to x^{-m/2}H^s_b(M;E)
\end{equation*}
of \eqref{InitialSobolevDomain} we have $\Dom_{\min}^s(A)\subset
\Dom$ and $\Dom\subset \Dom_{\max}^s(A)$ is closed (with respect
to the graph norm of $A$). These facts do not involve
$c$-ellipticity.

We will usually abbreviate $\Dom_{\min}^s(A)$ to $\Dom_{\min}^s$
and $\Dom_{\max}^s(A)$ to $\Dom_{\max}^s$ when the operator is
clear from the context. The operator $A$ with domain $\Dom$ will
be denoted by $A_{\Dom}$.

The proof of the following proposition characterizing
$\Dom_{\min}^s$ when $A$ is $c$-elliptic and $s$ is arbitrary is a
small variation of the characterization of $\Dom_{\min}^0$ as
given in Gil-Mendoza \cite{GiMe01}.

\begin{proposition}\label{CharacterizationDmin}
Let $A\in x^{-m}\Diff^m_b(M;E)$ be $c$-elliptic. Then
\begin{enumerate}[$(i)$]
\item $\Dom_{\min}^s=\Dom_{\max}^s\cap \big(\bigcap_{\eps>0} x^{m/2-\eps}H^{s+m}_b(M;E)\big)$.
\item $\Dom_{\min}^s = x^{m/2}H^{s+m}_b(M;E)$ if and only if $\spec_b(A)\cap \set{\Im \sigma=-m/2}=\varnothing$.
\end{enumerate}
\end{proposition}

Also the following theorem is a straightforward generalization of the
corresponding results for the case $s=0$, cf. Lesch \cite{Le97}, 
Gil-Mendoza \cite{GiMe01}.

\begin{theorem}\label{DomainCharakterisierung}
Let $A\in x^{-m}\Diff^m_b(M;E)$ be $c$-elliptic. Then:
\begin{enumerate}[$(i)$]
\item
With either of the domains $\Dom_{\min}^s$ or $\Dom_{\max}^s$, $A$
is Fredholm and so the former domain has finite codimension in the latter.
\item There is a one to one correspondence between the domains
$\Dom$ of closed extensions
\begin{equation*}
A: \Dom\subset x^{-m/2}H^s_b(M;E)\to x^{-m/2}H^{s}_b(M;E)
\end{equation*}
of \eqref{InitialSobolevDomain} and the subspaces of $\Dom_{\max}^s/\Dom_{\min}^s$.
\item For any sufficiently small $\eps>0$, the embeddings
\begin{equation*}\label{SobSpaceInclusions}
x^{m/2}H^{s+m}_b(M;E) \embed \Dom \embed x^{-m/2+\eps}H^{s+m}_b(M;E)
\end{equation*}
are continuous.
\item
Let $\Dom$ be such that $\Dom_{\min}^s\subset \Dom\subset \Dom_{\max}^s$. The operator $A : \Dom \to x^{-m/2}H^{s}_b(M;E)$
is Fredholm with index
\begin{equation}\label{WeakRelSInd}
\Ind A_{\Dom} = \Ind A_{\Dom_{\min}^s} + \dim \Dom/\Dom_{\min}^s.
\end{equation}
\end{enumerate}
\end{theorem}

Assuming that $A$ is $c$-elliptic, the space
$\Dom_{\max}^s/\Dom_{\min}^s$ can be identified with a (finite
dimensional) subspace $\Sing^s_{\max}\subset \Dom_{\max}^s$
complementary to $\Dom_{\min}^s$. Thus, the domains of the various
extensions of $A$ based on $x^{-m/2}H^s_b(M;E)$ are of the form
$\Dom_{\min}^s\oplus \Sing$ with $\Sing\subset \Sing_{\max}^s$. 
In fact, the complementary space can be chosen to be
independent of $s$, a subspace $\Sing_{\max}$ of $x^{-m/2}H^\infty_b(M;E)$,
\begin{equation*}
\Dom_{\max}^s=\Dom_{\min}^s\oplus \Sing_{\max}\quad \forall s\in \R.
\end{equation*}
A possible choice for $\Sing_{\max}$ is the orthogonal complement of
$\Dom_{\min}^0(A)$ in $\Dom_{\max}^0(A)$ with respect to the inner
product 
\begin{equation}\label{AInner}
(u,v)_A=(u,v)_{x^{-m/2}L^2_b}+(Au,Av)_{x^{-m/2}L^2_b},
\end{equation}
in other words, $\Sing_{\max}=\ker(A^\star A+\id)\cap \Dom^0_{\max}$. 
Another way to describe the complementary space is by means of
singular functions, see also Section~\ref{sec-Resolvents}.  

Granted this, one can then speak of the ``same'' extension of $A$
for different $s$; namely, if $\Sing\subset\Sing_{\max}$, let the
extension of $A$ based on $x^{-m/2}H^s_b(M;E)$ have domain
\begin{equation}\label{ComparableDomains}
\Dom^s=\Dom_{\min}^s\oplus\Sing.
\end{equation}
Then \eqref{WeakRelSInd} reads
\begin{equation*}
\Ind A_{\Dom^s} = \Ind A_{\Dom_{\min}^s} + \dim \Sing.
\end{equation*}
The index of $A_{\Dom_{\min}^s}$ is in fact independent of $s$. To
see this, we first observe that the kernel of $A$ in
$x^{-m/2}H^{-\infty}_b(M;E)$ is contained in
$x^{-m/2}H^{\infty}_b(M;E)$, and is therefore finite dimensional and
contained in each space $x^{-m/2}H^s_b(M;E)$. Next, using the
nonsingular sesquilinear pairing \begin{equation*}
x^{-m/2}H^s_b(M;E)\times x^{-m/2}H^{-s}_b(M;E)\ni
(u,v)\mapsto (u,v)_{x^{-m/2}L^2_b}\in \C,
\end{equation*}
we see that the annihilator of the range of $A_{\Dom_{\min}^s}$ is
the kernel $K^s$ of the formal adjoint $A^\star$ of $A$ acting on
$x^{-m/2}H^{-s}_b(M;E)$. Since $A^\star$ is also $c$-elliptic, its
kernel in $x^{-m/2}H^{-\infty}_b(M;E)$ is also a finite dimensional
subspace of $x^{-m/2}H^{\infty}_b(M;E)$. Thus $K^s$ is independent
of $s$. Since the range of $A_{\Dom_{\min}^s}$ is closed, this range
is the annihilator in $x^{-m/2}H^s(M;E)$ of $K^s$, so its
codimension is independent of $s$. Hence $\Ind A_{\Dom_{\min}^s}$ is
also independent of $s$. Thus:

\begin{proposition}\label{Spektralinvarianz}
Let $\Sing\subset \Sing_{\max}$ and define $\Dom^s$ as in
\eqref{ComparableDomains}. The index of
\begin{equation}\label{AmitDs}
A: \Dom^s\subset x^{-m/2}H^s_b(M;E)
\to x^{-m/2}H^{s}_b(M;E)
\end{equation}
is independent of $s$.
\end{proposition}

Let $P=x^m A$, an operator in $\Diff^m_b(M;E)$, and let $\lambda\in
\C$. Since $A-\lambda=x^{-m}(P-\lambda x^m)\in
x^{-m}\Diff^m_b(M;E)$, Proposition 4.1 of \cite{GiMe01} gives that
the minimal and maximal domains of $A-\lambda$ are those of $A$.
Since $A-\lambda\in x^{-m}\Diff^m_b(M;E)$ is $c$-elliptic if $A$ is
$c$-elliptic, also the kernel of
\begin{equation*}
A-\lambda:\Dom^s\subset x^{-m/2}H^s_b(M;E) \to x^{-m/2}H^s_b(M;E)
\end{equation*}
is independent of $s$ if $\Dom^s$ is the domain in \eqref{ComparableDomains}. Thus:
\begin{proposition}\label{SpectrumIndependent}
The spectrum of \eqref{AmitDs} is independent of $s$.
\end{proposition}

Sometimes it is useful to approximate a $c$-elliptic operator
$A\in x^{-m}\Diff_b^m(M;E)$ by operators having coefficients
independent of $x$ near the boundary $Y$ of $M$, 
see Definition~\ref{ConstantCoeff}. A simple and efficient 
approximation of $A$ can be obtained as follows.

Let $U_Y$ be a collar neighborhood of $Y$. For small $\tau>0$ let
\begin{equation*}
 \omega_\tau(x)=\omega(x/\tau)
\end{equation*}
where $\omega\in C_0^\infty(\overline\R_+)$ is a cut-off function with 
$\omega=1$ near $0$. Given $A$ let
\begin{equation}\label{eq:A_tau}
 A_{\tau} = \omega_\tau A_\wedge  + (1-\omega_\tau) A.
\end{equation}
For small enough $\tau > 0$ the operator $A_{\tau}$ is 
well defined, $c$-elliptic, and has the same conormal symbol and
therefore the same boundary spectrum as $A$.
Thus $\Dom_{\min}(A_{\tau})=\Dom_{\min}(A)$.
The following lemma was given in \cite{GiLoMe02}. Related results 
can also be found in \cite[Section 1.3]{Le97}.

\begin{lemma}\label{RedToConst}
As $\tau\to 0$, $A_{\tau}\to A$ in 
$\L(\Dom_{\min},x^{-m/2}L^2_b(M;E))$.
\end{lemma}
\begin{proof}
Since $A$ is $c$-elliptic, there is a bounded
parametrix $B:x^\gamma H^s_b\to x^{\gamma+m}H^{s+m}_b$ such that
\begin{equation*}
 R=I-BA:x^\gamma H^s_b\to x^{\gamma}H^\infty_b
\end{equation*}
is bounded for all $s$ and $\gamma$. 
Write $A=x^{-m}P$ and expand $P= P_0 + x\tilde P_1$ as in 
\eqref{bOperatorTaylor}. Then $x^{-m}P_0= A_\wedge$, and with 
$\tilde A=x^{-m}\tilde P_1$ we get
\begin{align*}
 A-A_{\tau} = x\omega_\tau \tilde A
 = x\omega_\tau \tilde A BA + x\omega_\tau \tilde A R
 = \tau \tilde\omega_\tau \tilde A BA + x\omega_\tau \tilde A R,
\end{align*}
where $\tilde\omega_\tau(x)=(x/\tau) \omega(x/\tau)$.
Now, $\tilde A B: x^{-m/2}L^2_b\to x^{-m/2}L^2_b$ is bounded, so if 
$u\in \Dom_{\min}(A)$, then
\begin{equation*}
 \|\tau \tilde\omega_\tau \tilde A BAu\|_{x^{-m/2}L^2_b} \leq
   c\,\tau\|Au\|_{x^{-m/2}L^2_b} \leq c\,\tau\|u\|_A.
\end{equation*}
Let $0<\alpha\ll 1$ and write $x\omega_\tau \tilde A R=
\tau^{1-\alpha}(\frac x{\tau})^{1-\alpha} \omega_\tau\,x^{\alpha}
\tilde A R$. The operator
\begin{equation*}
 x^{\alpha}\tilde A R: x^{m/2-\alpha}L^2_b\to x^{-m/2}L^2_b
\end{equation*}
and the embedding $(\Dom_{\min}(A),\|\cdot\|_A)\embed
x^{m/2-\alpha}L^2_b$ are both continuous, so
\begin{equation*}
 \|x \omega_\tau \tilde A Ru\|_{x^{-m/2}L^2_b}
 \leq \tilde c\,\tau^{1-\alpha}\|u\|_{x^{m/2-\alpha}L^2_b}
 \leq c\,\tau^{1-\alpha}\|u\|_A.
\end{equation*}
Altogether,
\begin{equation}\label{eq:DminEstimate}
 \|(A-A_{\tau})u\|_{x^{-m/2}L^2_b} \leq C\,\tau^{1-\alpha}\|u\|_A
\end{equation}
and thus $A_{\tau}\to A$ as $\tau\to 0$.
\end{proof}

In a similar way it can be shown that, for the formal adjoints, we
also have the convergence $A^\star_{\tau}\to A^\star$ as $\tau\to 0$. 

An immediate consequence of this lemma is the following result that 
was originally given in \cite[Section 1.3]{Le97}.
\begin{corollary} \label{AAtauSameDimension} 
For $A$ and $A_\tau$ as above, $\tau$ sufficiently small, we have
\[ \dim \Dom_{\max}(A)/\Dom_{\min}(A) =
   \dim \Dom_{\max}(A_\tau)/\Dom_{\min}(A_\tau). \]
\end{corollary}
\begin{proof}
We use the relative index formula \eqref{WeakRelSInd}
\begin{align*}
 \Ind A_{\tau,\Dom_{\max}} &= \Ind A_{\tau,\Dom_{\min}} +
 \dim\Dom_{\max}(A_\tau)/\Dom_{\min}(A_\tau), \\
 \Ind A_{\Dom_{\max}} &= \Ind A_{\Dom_{\min}} +
 \dim\Dom_{\max}(A)/\Dom_{\min}(A).
\end{align*}
By construction, $\Ind A_{\tau,\Dom_{\min}}=\Ind A_{\Dom_{\min}}$
and similarly $\Ind A^\star_{\tau,\Dom_{\min}}=
\Ind A^\star_{\Dom_{\min}}$ for $\tau$ sufficiently small.
This implies $\Ind A_{\tau,\Dom_{\max}}=\Ind A_{\Dom_{\max}}$
since $A^\star_{\tau,\Dom_{\min}}$ and $A^\star_{\Dom_{\min}}$ are
the Hilbert space adjoints of $A_{\tau,\Dom_{\max}}$ and
$A_{\Dom_{\max}}$, respectively. In conclusion, the quotient
spaces must have the same dimension.
\end{proof}

Similarly to the above, we consider extensions of the model operator
$$
A_\wedge : C_0^{\infty}(\open Y^{\wedge};E) \subset x^{-m/2}L_b^2(Y^{\wedge};E) \to
x^{-m/2}L_b^2(Y^{\wedge};E).
$$
Let $\Dom_{\wedge,\min}= \Dom_{\min}(A_{\wedge})$ be the completion of
$C_0^{\infty}(\open Y^{\wedge};E)$ with respect to the norm induced by the inner product
\begin{equation}\label{AWedgeInner}
(u,v)_{A_\wedge}=(u, v)_{x^{-m/2}L^2_b} + (A_\wedge u,A_\wedge v)_{x^{-m/2}L^2_b},
\end{equation}
and let
$$
\Dom_{\wedge,\max}= \Dom_{\max}(A_\wedge)=\set {u\in x^{-m/2} L^2_b(Y^\wedge;E): A_\wedge u\in x^{-m/2} L^2_b(Y^\wedge;E)}.
$$
Then
$$
A_{\wedge} : \Dom_{\wedge,\max} \subset x^{-m/2}L_b^2(Y^{\wedge};E) \to x^{-m/2}L_b^2(Y^{\wedge};E)
$$
is closed and densely defined, and $\Dom_{\wedge,\min} \subset \Dom_{\wedge,\max}$ is a closed
subspace with respect to the graph norm.
We have proved in \cite{GilKrainerMendoza1} that
$$
(1-\omega)\Dom_{\wedge,\max} = (1-\omega)\Dom_{\wedge,\min} = (1-\omega){\mathcal K}^{m,m/2}(Y^{\wedge};E)
$$
for all cut-off functions $\omega \in C_0^{\infty}(\overline{{\mathbb R}}_+)$ near zero, i.e.
$\omega = 1$ in a neighborhood of zero, and $\omega = 0$ near infinity.

Consequently, near infinity all domains $\Dom_{\wedge,\min} \subset
\Dom_{\wedge} \subset \Dom_{\wedge,\max}$ of $A_\wedge$ coincide
with $x^{\frac{n-m}{2}}H_{\rm cone}^m(Y^\wedge;E)$.
On the other hand, near the boundary, the closed extensions of
$A_\wedge$ are determined by its boundary spectrum which is the
same as the boundary spectrum of $A$. For this reason, many of
the results concerning the closed extensions of $A$ find their
analogs in the situation at hand. In fact, using an
approximation $A_\tau$ as in \eqref{eq:A_tau} with $\tau$ small,
one can easily describe the minimal and maximal extensions of
$A_\wedge$ on $Y^\wedge$ in terms of those of $A_\tau$ on the
manifold $M$. For instance, $u\in \Dom_{\max}(A_\wedge)$ if and
only if $(1-\omega)u\in x^{\frac{n-m}{2}}H_{\rm cone}^m(Y^\wedge;E)$
and $\omega u \in \Dom_{\max}(A_\tau)$ for some
cut-off function $\omega$ with small support and such that
$\omega=1$ near the boundary.

In particular, we have the embeddings
$$
\mathcal{K}^{m,m/2}(Y^\wedge;E) \embed \Dom_{\min}(A_\wedge) \embed
\Dom_{\max}(A_\wedge) \embed \mathcal{K}^{m,-m/2+\eps}(Y^\wedge;E).
$$
for some small $\eps>0$.

Because of \eqref{HomogeneityOfA} (with $\lambda=0$), both
$\Dom_{\min}(A_\wedge)$ and $\Dom_{\max}(A_\wedge)$ are
$\kappa$-invariant.

By the previous discussion, the following proposition is a direct
consequence of Proposition~\ref{CharacterizationDmin} and
Corollary~\ref{AAtauSameDimension}.

\begin{proposition}
Let $A\in x^{-m}\Diff^m_b(M;E)$ be $c$-elliptic. Then
\begin{enumerate}[$(i)$]
\item $\Dom_{\wedge,\min}=\Dom_{\wedge,\max}\cap
\big(\bigcap_{\eps>0} \mathcal K^{m,m/2-\eps}(Y^\wedge;E)\big)$.
\item $\Dom_{\wedge,\min} = \mathcal K^{m,m/2}(Y^\wedge;E)$ if and only
if $\spec_b(A)\cap \set{\Im \sigma=-m/2}=\varnothing$.
\item $\dim \Dom_{\wedge,\max}/\Dom_{\wedge,\min} =
\dim \Dom_{\max}(A)/\Dom_{\min}(A)$.
\end{enumerate}
\end{proposition}

Finally, we define the background spectrum of $A_\wedge$ as
\begin{align*}
\bgspec A_\wedge =\{\lambda\in \C : \;
& A_{\wedge,\Dom_{\wedge,\min}}-\lambda \text{ is not injective, or}\\
& A_{\wedge,\Dom_{\wedge,\max}}-\lambda\text{ is not surjective}\}.
\end{align*}
The complement $\bgres A_\wedge=\C\minus \bgspec A_\wedge$ is the background resolvent set.

\section{Ray conditions}
\label{sec-RayConditions}

The following theorem establishes the necessity of ray conditions
on the symbols of $A$ in order to have rays of minimal growth for
$A$ on some domain $\Dom$.

\begin{theorem}\label{NecessaryCondition}
Let $A\in x^{-m}\Diff^m(M;E)$ be $c$-elliptic. Suppose that there
is a domain $\Dom$, a ray
$$
\Gamma = \{z \in \C : z = re^{i\theta_0} \text{ for } r > 0\},
$$
and a number $R>0$ such that
$A - \lambda:\Dom\to x^{-m/2}L^2_b(M;E)$ is invertible for all
$\lambda\in\Gamma$ with $|\lambda|>R$. Suppose further that for
such $\lambda$, the resolvent
\begin{equation*}
(A_\Dom-\lambda)^{-1}:x^{-m/2}L^2_b(M;E)\to\Dom
\end{equation*}
is uniformly bounded in $\lambda$. Then
\begin{equation}\label{SymbolMinimalRay}
\bgspec A_\wedge\cap\Gamma = \varnothing\
\text{ and }\spec(\csym(A))\cap\bar\Gamma = 
\varnothing\text{ on }\cT^*M\minus 0.
\end{equation}
\end{theorem}

\begin{proof}
The hypotheses imply that $A-\lambda:\Dom_{\min}(A)\to
x^{-m/2}L^2_b(M;E)$ is injective for $\lambda\in \Gamma$ and
that, in fact, if $u\in \Dom_{\min}(A)$, then
\begin{equation}\label{InjectiveEstimate}
 \|(A-\lambda)u\| \ge C\|u\|_A
\end{equation}
for some constant $C>0$. Here $\|\cdot\|$ denotes the norm in
$x^{-m/2}L^2_b$ and $\|\cdot\|_A$ is the graph norm. We first
prove that
\begin{equation*}
A_\wedge-\lambda:\Dom_{\min}(A_\wedge)\to
   x^{-m/2}L^2_b(Y^\wedge;E) \;\text{ is injective}.
\end{equation*}
Note that $\Dom_{\min}(A_\wedge)$ and $\Dom_{\max}(A_\wedge)$
are invariant under the dilation $\kappa_\varrho$. 
If $v\in C_0^\infty(\open Y^\wedge;E)$, then for $\varrho>0$ small,
$\kappa_\varrho^{-1} v\in\Dom_{\min}(A_\wedge)$ is supported near
$Y$, the boundary of $Y^\wedge$, and gives an
element $\kappa_\varrho^{-1} v$ of $\Dom_{\min}(A)$. We have
\begin{align*}
 \|(\varrho^m\kappa_\varrho 
 A\kappa_\varrho^{-1}-\lambda)v\|
 &= \varrho^m\|\kappa_\varrho(A-\varrho^{-m}\lambda)
 \kappa_\varrho^{-1}v\| \\
 &= \varrho^m \|(A-\varrho^{-m}\lambda)\kappa_\varrho^{-1}v\|
\end{align*}
because $\kappa_\varrho$ is an isometry. Next, if $A-\lambda$ is injective,
then obviously so is $A-\varrho^{-m}\lambda$ for $\varrho\le 1$,
and by \eqref{InjectiveEstimate},
\begin{equation*}
\varrho^m \|(A-\varrho^{-m}\lambda)\kappa_\varrho^{-1}v\| \ge C \varrho^m \|\kappa_\varrho^{-1}v\|_A.
\end{equation*}
But
\begin{align*}
\varrho^m \|\kappa_\varrho^{-1}v\|_A
&=\varrho^m\|\kappa_\varrho^{-1}v\|
    +\varrho^m\|A\kappa_\varrho^{-1}v\|\\
&= \varrho^m\|v\| +\|\varrho^m \kappa_\varrho 
    A\kappa_\varrho^{-1}v\|
\end{align*}
using again that $\kappa_\varrho$ is an isometry. Thus
\begin{equation*}
\|(\varrho^m\kappa_\varrho A\kappa_\varrho^{-1}-\lambda)v\|
\geq C\bigl(\varrho^m\|v\|
    +\|\varrho^m \kappa_\varrho A\kappa_\varrho^{-1}v\|\bigr)
\end{equation*}
for some $C>0$ and all small $\varrho$. In view of the definition
of $A_\wedge$, taking the limit as $\varrho\to 0$ we arrive at
\begin{equation}\label{AwedgeEstimate}
 \| (A_\wedge -\lambda)v\| \ge C\|A_\wedge v\|
\end{equation}
for all $v\in C_0^\infty(\open Y^\wedge;E)$. Now,
for an arbitrary $v\in\Dom_{\min}(A_\wedge)$ there exist a
sequence $\{v_k\}\subset C_0^\infty(\open Y^\wedge;E)$ such that
$v_k\to v$ and $A_\wedge v_k\to A_\wedge v$ in $x^{-m/2}L^2_b$ as
$k\to \infty$, so $(A_\wedge-\lambda)v_k\to (A_\wedge-\lambda) v$
in $x^{-m/2}L^2_b$.  Thus, since \eqref{AwedgeEstimate} holds for
the $v_k$, it holds for any $v\in \Dom_{\min}(A_\wedge)$.

The estimate \eqref{AwedgeEstimate} implies the injectivity of
$A_\wedge-\lambda$ on $\Dom_{\min}(A_\wedge)$ for $\lambda\ne 0$.
Indeed, if $(A_\wedge-\lambda)v=0$, then $A_\wedge v = 0$, so
$\lambda v= 0$. Thus $v=0$ since $\lambda\ne 0$.

The surjectivity of $A_\wedge-\lambda:\Dom_{\max}(A_\wedge)\to
x^{-m/2}L^2_b(Y^\wedge;E)$ follows from the injectivity of
$A_\wedge^\star-\overline\lambda:\Dom_{\min}(A_\wedge^\star)\to
x^{-m/2}L^2_b(Y^\wedge;E)$. The latter is a consequence of the
injectivity of $(A^\star-\overline \lambda)$ on
$\Dom_{\min}(A^\star)$ for $\lambda \in \Gamma$ and the above
argument. This proves the first assertion in
\eqref{SymbolMinimalRay}.

We now prove the second assertion. Since $A$ is $c$-elliptic,
$A_\wedge$ is elliptic in the usual sense in the interior of
$Y^\wedge$. So the usual elliptic \emph{a priori} estimate holds
in compact subsets of $\open Y^\wedge$. Thus there is a
constant $C>0$ such that
\begin{equation*}
 \|v\|_{\K^{m,m/2}}\leq C\big(\|A_\wedge v\| + \|v\|\big)
\end{equation*}
for $v\in \K^{m,m/2}(Y^\wedge;E)$, $\supp v\subset \set{1\leq x\leq
2}\times Y$. The inequality \eqref{AwedgeEstimate} now gives
\begin{equation}\label{PellipticEstimate}
 \|v\|_{\K^{m,m/2}}\leq C\big(\|(A_\wedge -\lambda)v\| + \|v\|\big)
\end{equation}
for $v\in \K^{m,m/2}(Y^\wedge;E)$, $\supp v\subset \set{1\leq x\leq
2}\times Y$, with some $C$ independent of $\lambda$. By standard
arguments (see e.g. Seeley \cite{SeeleyBook}) this gives that
$\sym(A_\wedge)-\lambda$
is invertible for $\lambda\in \Gamma$ when $1\leq x\leq 2$. But
\begin{equation*}
\sym(A_\wedge)(x,y;\xi,\eta)-\lambda=x^{-m}\big( 
\csym(A_\wedge)(y;x\xi,\eta)-x^m \lambda\big).
\end{equation*}
In this formula we made use of the fact that the $c$-symbol of
$A_\wedge$ is independent of $x$. Replacing $x\xi$ by $\xi$ and
$x^m \lambda$ by $\lambda$, and using that
$\csym(A_\wedge)=\csym(A)|_Y$ we reach the conclusion that
\begin{equation*}
\csym(A)-\lambda
\end{equation*}
is invertible over $Y$, and therefore over a neighborhood of $Y$
in $M$, when $\lambda\in \Gamma$. The hypothesis on $A$ also
implies estimates like \eqref{PellipticEstimate} for $A$ on
compact subsets of the interior of $M$. Thus also
$\sym(A)-\lambda$ is invertible over compact subsets of the
interior of $M$ when $\lambda\in \Gamma$. This gives the second
statement in \eqref{SymbolMinimalRay}.
\end{proof}

The following is a partial converse of Theorem~\ref{NecessaryCondition}.

\begin{theorem}\label{SufficientCondition}
Let $A\in x^{-m}\Diff_b^m(M;E)$ be $c$-elliptic. If
\eqref{SymbolMinimalRay} holds, then there exists a domain $\Dom$
such that $\spec A_{\Dom}$ is discrete.
\end{theorem}

\begin{proof}
We will use the parametrix from Section~\ref{sec-Parametrix} to prove
the statement. First of all, the compactness of $M$ and the
spectral condition on the symbol $\csym(A)$ imply that there exists
some closed sector $\Lambda$ with $\Gamma\subset\open \Lambda$ such that
$\spec(\csym(A))\cap \Lambda =\varnothing$ on $\cT^*M \minus 0$.
Consequently, $A-\lambda$ is $c$-elliptic with parameter $\lambda \in
\Lambda$, cf. Definition~\ref{ConeParameterElliptic}.  We choose
$\Lambda$ in such a way that $\Lambda\minus\{0\}\subset \bgres A_\wedge$ 
also holds; this is possible because $\bgres A_\wedge$ is a union of
open sectors, see \cite{GilKrainerMendoza1}.  Then, for
$\lambda\in\Lambda\minus\{0\}$, we also have that $A_\wedge
-\lambda:\Dom_{\min}(A_\wedge)\to x^{-m/2}L^2_b(Y^\wedge;E)$ is
injective and therefore, by Theorem~\ref{PDrei},
\begin{equation*}
 A-\lambda:\Dom_{\min}(A)\to x^{-m/2}L^2_b(M;E)
\end{equation*}
is injective for $\lambda$ sufficiently large.

On the other hand, the surjectivity of $A_\wedge
-\lambda:\Dom_{\max}(A_\wedge)\to x^{-m/2}L^2_b(Y^\wedge;E)$ implies
the injectivity of $A_\wedge^\star -\bar \lambda$ on
$\Dom_{\min}(A_\wedge^\star)$. Since $A^\star - \bar\lambda$ is also
$c$-elliptic with parameter $\bar \lambda$ in the complex conjugate of
$\Lambda$, we can use Theorem~\ref{PDrei} with $A^\star$ instead of
$A$ to conclude that $A^\star-\bar\lambda:\Dom_{\min}(A^\star)\to
x^{-m/2}L^2_b(M;E)$ is injective for $\bar\lambda$ sufficiently large.
Thus, for such $\lambda$, we get the surjectivity of
\begin{equation*}
 A-\lambda:\Dom_{\max}(A)\to x^{-m/2}L^2_b(M;E).
\end{equation*}
Consequently, for $\lambda$ large, $A - \lambda$ is injective on $\Dom_{\min}$ and
surjective on $\Dom_{\max}$, and hence there exists a domain $\Dom$ such that
\begin{equation*}
 A_{\Dom}-\lambda:\Dom\to x^{-m/2}L^2_b(M;E)
\end{equation*}
is invertible. Thus $\spec A_{\Dom} \not= \C$, so it must be discrete.
\end{proof}

Observe that for $\lambda\in\Gamma$, $|\lambda|>R>0$, the norm
$\|(A_\Dom-\lambda)^{-1}\|_{\L(x^{-m/2}L^2_b(M;E),\Dom)}$ is
uniformly bounded if and only if
\begin{equation*}
 \|(A_\Dom-\lambda)^{-1}\|_{\L(x^{-m/2}L^2_b(M;E))} =O(|\lambda|^{-1}) 
\;\text{ as } |\lambda|\to\infty. 
\end{equation*}
Stronger and more precise statements about resolvents of elliptic
cone operators will be given in Section~\ref{sec-Resolvents}.

\section{Parametrix construction}
\label{sec-Parametrix}

In this section we assume $\Lambda$ to be a closed sector in $\C$ of
the form
\begin{equation*}
 \Lambda = \set{z\in\C : z=re^{i\theta} \text{ for } r\ge 0,\
 \theta\in \R, \ |\theta-\theta_0|\le a}
\end{equation*}
for some real $\theta_0$ and $a>0$, and assume that $A-\lambda$ is
$c$-elliptic with parameter $\lambda\in\Lambda$ according to
Definition~\ref{ConeParameterElliptic}, and that
\begin{equation}\label{Adachinjektiv1}
 A_\wedge - \lambda : \Dom_{\min}(A_{\wedge}) \to
 x^{-m/2}L^2_b(Y^{\wedge};E) \text{ is injective if }\lambda \in \Lambda \minus
\set{0}.
\end{equation}
Our goal is to construct a parameter-dependent parametrix of
\begin{equation}\label{Aminuslambdamin}
A - \lambda : \Dom_{\min}(A) \to x^{-m/2}L^2_b(M;E)
\end{equation}
by means of three crucial steps that we proceed to outline.

\medskip 
{\sc Step 1:} The first step is concerned with the construction of a
pseudodifferential parametrix $B_1(\lambda)$ of
$A-\lambda: C_0^{\infty}(\open M;E) \to C_0^{\infty}(\open M;E)$
taking care of the degeneracy of the complete symbol
of $A - \lambda$ near the boundary of $M$. The parametrix
$B_1(\lambda)$ is constructed within a corresponding (sub)calculus
of parameter-dependent pseudodifferential operators that are built
upon degenerate symbols.

\medskip 
{\sc Step 2:} In the second step the parametrix $B_1(\lambda)$ is
refined to a parametrix
\begin{equation*}
B_2(\lambda) : x^{-m/2}L^2_b(M;E) \to \Dom_{\min}(A)
\end{equation*}
which is continuous and pointwise a Fredholm inverse of $A - \lambda$.
The remainders
\begin{align}\label{Reste2}
B_2(\lambda)(A - \lambda) - 1 &: \Dom_{\min}(A) \to \Dom_{\min}(A), \\
\label{Reste2b} (A - \lambda)B_2(\lambda) - 1 &:
x^{-m/2}L^2_b(M;E) \to x^{-m/2}L^2_b(M;E)
\end{align}
are parameter-dependent smoothing pseudodifferential operators in
\begin{equation*}
C_0^{\infty}(\open M;E) \to C^{\infty}(\open M;E)
\end{equation*}
since $B_2(\lambda)$ is a refinement of $B_1(\lambda)$, but the
operator norms in the spaces \eqref{Reste2} and \eqref{Reste2b} are
not decreasing as $|\lambda| \to \infty$.

\medskip 
{\sc Step 3:} While in the first two steps we only make use of the
$c$-ellipticity with parameter, we now need the additional
requirement that \eqref{Adachinjektiv1} holds. In view of the
$\kappa$-homogeneity of $A_\wedge-\lambda$,
$$
A_{\wedge} - \varrho^m\lambda =
\varrho^m\kappa_{\varrho}(A_{\wedge} -
\lambda)\kappa_{\varrho}^{-1}
\;\text{ for } \lambda\ne 0,\; \varrho > 0,
$$
 we only need to require
\eqref{Adachinjektiv1} for $|\lambda| = 1$. Recall that the minimal
domain $\Dom_{\min}(A_{\wedge})$ is invariant under the action of
$\kappa_\varrho$.

Under the additional assumption \eqref{Adachinjektiv1} we will refine 
$B_2(\lambda)$ to obtain a parameter-dependent parametrix 
$B(\lambda)$ such that
$$
B(\lambda)(A - \lambda) - 1 : \Dom_{\min}(A) \to \Dom_{\min}(A)
$$
is compactly supported in $\lambda \in \Lambda$. In particular,
for $\lambda$ sufficiently large the operator family $A -
\lambda : \Dom_{\min}(A) \to x^{-m/2}L^2_b(M;E)$ is injective, and
the parametrix $B(\lambda)$ is a left-inverse. 
Moreover, for $\lambda$ large, the smoothing remainder
$$
\Pi(\lambda)= 1 - (A - \lambda)B(\lambda)
$$
is a projection on $x^{-m/2}L^2_b(M;E)$ to a complement of the
range of $A - \lambda$ on $\Dom_{\min}(A)$, i.e., $(A -
\lambda)B(\lambda)$ is a projection onto $\rg(A_{\min} - \lambda)$.

For the final construction of $B(\lambda)$ we adopt Schulze's viewpoint 
from the pseudo\-differential edge-calculus, see e.g. \cite{SzNorthHolland,
SzWiley98}, and add extra conditions of trace and potential type 
within a suitably defined class of Green remainders.

\bigskip 
We now proceed to construct a suitable parametrix of $A-\lambda$ as
outlined above.  The first step is the parametrix construction in the
interior of the manifold, assuming only that $A-\lambda$ is
$c$-elliptic with parameter in a closed sector $\Lambda\subset\C$.

On $M$ we fix a collar neighborhood diffeomorphic to $[0,1)\times
Y$, $Y= \partial{M}$, and consider local coordinates of the form
$[0,1)\times\Omega \subset \overline{\R}_+\times\R^{n-1}$ near the
boundary, where $\Omega \subset \R^{n-1}$ corresponds to a chart
on $Y$. Moreover, these coordinates are chosen in such a way that
the push-forward of the vector bundle $E$ is trivial on
$[0,1)\times\Omega$ (e.g., choose $\Omega$ contractible).

In these coordinates the operator $A - \lambda$ takes the form
\begin{equation}\label{Aminuslambdalokal}
A - \lambda = x^{-m}\bigg(\sum_{k+|\alpha| \leq m}
a_{k\alpha}(x,y)D_y^{\alpha}(xD_x)^k - x^{m}\lambda\bigg),
\end{equation}
where the $a_{k\alpha}$ are smooth matrix-valued coefficients on
$[0,1)\times\Omega$. 
The $c$-ellipticity with parameter of the family $A-\lambda$ implies
that, in the interior of $M$, it is elliptic with parameter in the usual 
sense,  and in local coordinates near the boundary, 
\begin{equation*}
 \sum_{k+|\alpha| = m}a_{k\alpha}(x,y)\eta^{\alpha}\xi^k -\lambda
\end{equation*}
is invertible for all $(\xi,\eta,\lambda) \in
\bigl(\R\times\R^{n-1}\times\Lambda\bigr)\minus \{0\}$ and
$(x,y) \in [0,1)\times\Omega$.

 From equation \eqref{Aminuslambdalokal} we deduce that the
complete symbol of $A - \lambda$ in $(0,1)\times\Omega$ is of the
form $x^{-m}a(x,y,x\xi,\eta,x^{m}\lambda)$ for some
parameter-dependent classical symbol $a(x,y,\xi,\eta,\lambda)$ of
order $m$, and the $c$-ellipticity condition near the boundary is 
equivalent to the invertibility of the principal component
$a_{(m)}(x,y,\xi,\eta,\lambda)$ of $a$.
These observations give rise to the class of parameter-dependent
pseudodifferential operators that we will consider below.

For the rest of this section we will work (without loss of
generality) with scalar symbols; the general case of matrix-valued
symbols is straightforward.

Sometimes we will denote the variables in $(0,1)\times\Omega\,$ 
by $z=(x,y)$ and $z'=(x',y')$, and the corresponding covariables in
$\R^n$ by $\zeta=(\xi,\eta)\in \R\times\R^{n-1}$. 

\begin{definition}\label{DegClass}
For $\mu \in \R$ let $\Psi^{\mu}(\Lambda)$ denote the
space of all pseudodifferential operators
\begin{equation*}
A(\lambda) : C_0^{\infty}((0,1)\times\Omega) \to
C^{\infty}((0,1)\times\Omega)
\end{equation*}
depending on the parameter $\lambda \in \Lambda$ of the form
\begin{equation}\label{AllgemeineParameterOps}
A(\lambda) u(z) = \frac{1}{(2\pi)^n}\iint
e^{i(z-z')\cdot\zeta}\,\tilde{a}(z,\zeta,\lambda)
u(z')\,dz'\,d\zeta + C(\lambda)u(z)
\end{equation}
for $z,z'\in (0,1)\times\Omega$, $\zeta\in\R^n$, 
where the family $C(\lambda) \in \Psi^{-\infty}(\Lambda)$ is a
parameter-dependent smoothing operator of the form
\begin{equation*}
C(\lambda)u(z) =
\int k(z,z',\lambda)u(z')\,dz'
\end{equation*}
with rapidly decreasing integral kernel $k(z,z',\lambda) \in
\S(\Lambda,C^{\infty}((0,1)\times\Omega\times (0,1)\times\Omega))$, 
and where the symbol 
$\tilde a(z,\zeta,\lambda)=\tilde{a}(x,y,\xi,\eta,\lambda)$ satisfies
$$
\tilde{a}(x,y,\xi,\eta,\lambda) =
x^{-\mu}a(x,y,x\xi,\eta,x^{d}\lambda)
$$
with $a(x,y,\xi,\eta,\lambda) \in
C^{\infty}([0,1)\times\Omega\times\R\times \R^{n-1}\times\Lambda)$
satisfying for all multi-indices $\alpha$, $\beta$, and $\gamma$,
the symbol estimates
$$
|\partial_{(x,y)}^{\alpha}\partial_{(\xi,\eta)}^{\beta}
\partial_{\lambda}^{\gamma}a(x,y,\xi,\eta,\lambda)|
= O\bigl(\bigl(1 + |\xi| + |\eta| +
|\lambda|^{1/d}\bigr)^{\mu-|\beta|-d|\gamma|}\bigr) $$ as
$|(\xi,\eta,\lambda)| \to \infty$, locally uniformly for $(x,y)
\in [0,1)\times\Omega$. Here $d \in {\N}$ is a fixed parameter for
the class $\Psi^{\infty}(\Lambda)$ which refers to the anisotropy;
in the case of the operator $A-\lambda$ we have $d=m=\textup{ord}(A)$.
Moreover, the symbol $a(x,y,\xi,\eta,\lambda)$ is assumed to be
classical: It admits an asymptotic expansion
\begin{equation}\label{ClassicalExpansion}
a \sim \sum_{j=0}^{\infty}\chi(\xi,\eta,\lambda)
a_{(\mu-j)}(x,y,\xi,\eta,\lambda),
\end{equation} 
where $\chi \in C^{\infty}({\R}\times\R^{n-1}\times\Lambda)$ is a
function such that $\chi = 0$ near the origin and $\chi = 1$ for
$|(\xi,\eta,\lambda)|$ large, and the components 
$a_{(\mu-j)}(x,y,\xi,\eta,\lambda)$ satisfy the homogeneity relation
\begin{equation*}
a_{(\mu-j)}(x,y,\varrho\xi,\varrho\eta,\varrho^d\lambda) =
\varrho^{\mu-j}a_{(\mu-j)}(x,y,\xi,\eta,\lambda)
\end{equation*}
for $\varrho > 0$ and $(\xi,\eta,\lambda) \in
({\R}\times\R^{n-1}\times\Lambda) \minus \{0\}$. The
parameter-dependent principal symbol of $A(\lambda)$ is then given
by $x^{-\mu}a_{(\mu)}(x,y,x\xi,\eta,x^d\lambda)$.
\end{definition}

Note that the symbol $a(x,y,\xi,\eta,\lambda)$ is smooth in $x$ up to $x=0$.

\begin{proposition}\label{DegClassComp}
Let $A(\lambda) \in \Psi^{\mu_1}(\Lambda)$ and $B(\lambda) \in
\Psi^{\mu_2}(\Lambda)$ with either $A(\lambda)$ or $B(\lambda)$
being properly supported, uniformly in $\lambda \in \Lambda$. Then
the composition
$$
A(\lambda)B(\lambda) : C_0^{\infty}((0,1)\times\Omega) \to
C^{\infty}((0,1)\times\Omega)
$$
belongs to $\Psi^{\mu_1+\mu_2}(\Lambda)$.
\end{proposition}
\begin{proof}
Let $\tilde{a}(x,y,\xi,\eta,\lambda)$ and
$\tilde{b}(x,y,\xi,\eta,\lambda)$ be complete symbols associated
with $A(\lambda)$ and $B(\lambda)$ according to
\eqref{AllgemeineParameterOps}. Then the corresponding complete symbol 
of the composition has the asymptotic expansion
\begin{equation}\label{LeibnizProdukt}\tag{\ref{DegClassComp}a}
\sum\limits_{k+|\alpha|=0}^{\infty}\frac{1}{k!\alpha!}
\partial_{\xi}^k\partial_{\eta}^{\alpha}\tilde{a}(x,y,\xi,\eta,\lambda)
D_x^kD_y^{\alpha}\tilde{b}(x,y,\xi,\eta,\lambda).
\end{equation}
Now write
\begin{align*}
\tilde{a}(x,y,\xi,\eta,\lambda)
  &= x^{-\mu_1}a(x,y,x\xi,\eta,x^{d}\lambda), \\
\tilde{b}(x,y,\xi,\eta,\lambda)
  &= x^{-\mu_2}b(x,y,x\xi,\eta,x^{d}\lambda)
\end{align*}
with $a$ and $b$ as in Definition~\ref{DegClass}. This gives
$$
\partial_{\xi}^k\partial_{\eta}^{\alpha}\tilde{a}(x,y,\xi,\eta,\lambda)
= x^{-\mu_1}\bigl(\partial_{\xi}^k\partial_{\eta}^{\alpha}
a\bigr)(x,y,x\xi,\eta,x^{d}\lambda) x^k.
$$
Since $(xD_x)D_y^{\alpha}\tilde{b}(x,y,\xi,\eta,\lambda)$ equals
$$
x^{-\mu_2} \bigl((-\mu_2 + xD_x + {\xi}D_{\xi} +
d\lambda_1D_{\lambda_1} + d\lambda_2D_{\lambda_2})
D_y^{\alpha}b\bigr)(x,y,x\xi,\eta,x^{d}\lambda),
$$
and since $x^kD_x^k = \sum\limits_{j=0}^k c_{kj}(xD_x)^j$ with
some universal constants $c_{kj}$, we see that each term in the
asymptotic expansion \eqref{LeibnizProdukt} is of the form
$$
\frac{1}{k!\alpha!}\partial_{\xi}^k\partial_{\eta}^{\alpha}
\tilde{a}(x,y,\xi,\eta,\lambda)D_x^kD_y^{\alpha}
\tilde{b}(x,y,\xi,\eta,\lambda) =
x^{-(\mu_1+\mu_2)}p_{k,\alpha}(x,y,x\xi,\eta,x^{d}\lambda)
$$
with a parameter-dependent symbol $p_{k,\alpha}$ of order
$\mu_1+\mu_2-k-|\alpha|$ which satisfies the conditions of
Definition \ref{DegClass}. In conclusion, if $p$ is such that
\begin{equation*}
p(x,y,\xi,\eta,\lambda) \sim \sum\limits_{k+|\alpha|=0}^{\infty}
p_{k,\alpha}(x,y,\xi,\eta,\lambda),
\end{equation*}
then $x^{-(\mu_1+\mu_2)}p(x,y,x\xi,\eta,x^{d}\lambda)$ is a complete
symbol of the composition $A(\lambda)B(\lambda)$ and the proposition
follows.
\end{proof}

\begin{definition}\label{Lokalcellipt}
Let $A(\lambda)\! \in \Psi^{\mu}(\Lambda)$ with 
principal symbol $x^{-\mu}a_{(\mu)}(x,y,x\xi,\eta,x^d\lambda)$. 
The family $A(\lambda)$ is said to be $c$-elliptic with parameter
$\lambda \in \Lambda$ if $a_{(\mu)}(x,y,\xi,\eta,\lambda)$ is
invertible for all $(x,y) \in [0,1)\times\Omega$ and
$(\xi,\eta,\lambda) \in ({\R}\times\R^{n-1}\times\Lambda)\minus
\{0\}$.
\end{definition}

\begin{proposition}\label{LokalcelliptParametrix}
For $A(\lambda) \in \Psi^{\mu}(\Lambda)$ the following are
equivalent:
\begin{enumerate}[$(i)$]
\item $A(\lambda)$ is $c$-elliptic with parameter $\lambda \in
\Lambda$.
\item There exists a parametrix $Q(\lambda)\in \Psi^{-\mu}(\Lambda)$,
properly supported $($uniformly in $\lambda)$, such that
$A(\lambda)Q(\lambda) - 1$  and $Q(\lambda)A(\lambda)-1$
both belong to $\Psi^{-\infty}(\Lambda)$. 
\end{enumerate}
\end{proposition}
\begin{proof}
For the proof we need the auxiliary operator class
$\Psi^{\mu,0}(\Lambda)= x^{\mu}\Psi^{\mu}(\Lambda)$. It is easy to
see from the proof of Proposition \ref{DegClassComp} that the
composition gives rise to
$$
\Psi^{\mu_1,0}(\Lambda)\times\Psi^{\mu_2,0}(\Lambda) \to
\Psi^{\mu_1+\mu_2,0}(\Lambda)
$$
provided that one of the factors is properly supported (uniformly
in $\lambda$). Actually, it is not necessary to couple the weight
factor and the order of the operators as it is done for the elements 
of $\Psi^{\mu}(\Lambda)$.

Let $A(\lambda) \in \Psi^{\mu}(\Lambda)$ be $c$-elliptic with
parameter. Without loss of generality assume that $A(\lambda)$ is
properly supported, uniformly in $\lambda$. Let $Q'(\lambda) \in
\Psi^{-\mu}(\Lambda)$ be properly supported with complete symbol
$x^{\mu}(\chi{\cdot}a_{(\mu)}^{-1})(x,y,x\xi,\eta,x^d\lambda)$,
where $\chi$ is as in \eqref{ClassicalExpansion}. Thus
$R_r(\lambda) = A(\lambda)Q'(\lambda) - 1$ and 
$R_l(\lambda) = Q'(\lambda)A(\lambda) - 1$ both belong to 
$\Psi^{-1,0}(\Lambda)$, and are properly supported,
uniformly in $\lambda$. For $k \in \N$ let
$r_k(x,y,\xi,\eta,\lambda)$ be of order $-k$ such that
$r_k(x,y,x\xi,\eta,x^{d}\lambda)$ is a complete symbol of
$R_l^k(\lambda) \in \Psi^{-k,0}(\Lambda)$. Let
$r(x,y,\xi,\eta,\lambda)$ be of order $-1$ such that
$$
r(x,y,\xi,\eta,\lambda) \sim
\sum\limits_{k=1}^{\infty}(-1)^kr_k(x,y,\xi,\eta,\lambda),
$$
and let $R'(\lambda) \in \Psi^{-1,0}(\Lambda)$ be properly supported
having $r(x,y,x\xi,\eta,x^d\lambda)$ as complete symbol. Then 
$$
(1 + R'(\lambda))Q'(\lambda)A(\lambda) - 1 \in \bigcap\limits_{k
\in \N}\Psi^{-k,0}(\Lambda) = \Psi^{-\infty}(\Lambda),
$$
so $(1 + R'(\lambda))Q'(\lambda)\in \Psi^{-\mu}(\Lambda)$ is a left
parametrix of $A(\lambda)$. In the same way we obtain a right parametrix.  
The other direction of the proposition is immediate.
\end{proof}

We now pass to the collar neighborhood $[0,1){\times}Y \subset M$:
The restriction of the bundle $E$ to $[0,1){\times}Y$ is
isomorphic to the pull-back of a bundle on $Y$. For simplicity, we
denote this bundle by the same letter $E$, and the sections of the
bundle $E$ on $[0,1){\times}Y$ are then represented as
$C^{\infty}([0,1),C^{\infty}(Y;E))$. We consider families of
pseudodifferential operators
$$
A(\lambda) : C_0^{\infty}((0,1),C^{\infty}(Y;E)) \to
C^{\infty}((0,1),C^{\infty}(Y;E))
$$
on $(0,1){\times}Y$ acting in sections of the bundle $E$ which
depend anisotropically on the parameter $\lambda \in \Lambda$.
With respect to the fixed splitting of variables these operators 
can be written as follows:
\begin{equation}\label{AllgemeineOpsimKragen}
A(\lambda)u(x) = \frac{1}{2\pi}\iint
e^{i(x-x')\xi}\tilde{a}(x,\xi,\lambda)u(x')\,dx'\,d\xi + C(\lambda)u(x)
\end{equation}
for $x, x'\in (0,1)$, $\xi\in\R$, where $C(\lambda) \in 
\Psi^{-\infty}(\Lambda)$ is a parameter-dependent smoothing operator
\begin{equation*}
C(\lambda)u(x) = \int k(x,x',\lambda)u(x')\,dx'
\end{equation*}
with integral kernel $k(x,x',\lambda) \in 
\S(\Lambda,C^{\infty}((0,1){\times}(0,1),L^{-\infty}(Y)))$.
As in the local case, cf. Definition~\ref{DegClass}, we use here the  
notation $\Psi^{-\infty}(\Lambda)$ for the remainder class.

Moreover, the symbol $\tilde{a}(x,\xi,\lambda)$ is a smooth
function of $x \in (0,1)$ taking values in the space
$L^{\mu,(1,d)}(Y;{\mathbb R}\times\Lambda)$ of pseudodifferential
operators of order $\mu \in {\mathbb R}$ on $Y$ depending on the
parameters $(\xi,\lambda) \in {\mathbb R}\times\Lambda$. Recall
that a family of operators
\begin{equation*}
B(\xi,\lambda) : C^{\infty}(Y;E) \to C^{\infty}(Y;E)
\end{equation*}
belongs to $L^{\mu,(1,d)}(Y;{\mathbb R}\times\Lambda)$ if, in a local
patch $\Omega$, it is of the form
\begin{equation*}
B(\xi,\lambda)u(y) = \frac{1}{(2\pi)^{n-1}}\iint e^{i(y-y')\eta} 
b(y,\xi,\eta,\lambda)u(y')\,dy'\,d\eta + D(\xi,\lambda)u(y)
\end{equation*}
for $y, y'\in\Omega$, $\eta\in\R^{n-1}$, where
\begin{equation*}
D(\xi,\lambda)u(y) = \int c(y,y',\xi,\lambda)u(y')\,dy'
\end{equation*}
with integral kernel $c(y,y',\xi,\lambda) \in 
\S(\R\times\Lambda,C^{\infty}(\Omega\times\Omega))$,
and where the symbol $b(y,\xi,\eta,\lambda)$ satisfies the
symbol estimates of Definition~\ref{DegClass}, but here in the
$x$-independent case.

As before, we do not consider general families of pseudodifferential
operators on $(0,1){\times}Y$ and restrict ourselves to operators in
$\Psi^{\mu}(\Lambda)$ where the symbol $\tilde{a}(x,\xi,\lambda)$ in
\eqref{AllgemeineOpsimKragen} is required to be of the form
\begin{equation*}
\tilde{a}(x,\xi,\lambda) = x^{-\mu}a(x,x\xi,x^{d}\lambda),
\end{equation*}
where $a(x,\xi,\lambda)$ is smooth in $x \in [0,1)$ with values in
$L^{\mu,(1,d)}(Y;\R\times\Lambda)$. Observe that this is
precisely the class of operators that is obtained via globalizing
the local classes from Definition~\ref{DegClass} to the collar
neighborhood $(0,1){\times}Y$.

The parameter-dependent homogeneous principal symbol of an
operator in $\Psi^{\mu}(\Lambda)$ extends to an anisotropic
homogeneous section on
$(\cT^*([0,1){\times}Y)\times\Lambda)\minus 0$, and the global 
meaning of the $c$-ellipticity from Definition~\ref{Lokalcellipt} 
is the invertibility of the principal symbol there. From 
Proposition~\ref{LokalcelliptParametrix} we get the following:

\begin{proposition}\label{KragencelliptParametrix}
There exists a parametrix $Q(\lambda)\in \Psi^{-m}(\Lambda)$ of $A
- \lambda$ which is properly supported $($uniformly in $\lambda)$
and has the form
\begin{equation*}
Q(\lambda)u(x) = \frac{1}{2\pi}\iint e^{i(x-x')\xi}
\tilde{p}(x,\xi,\lambda)u(x') \,dx'\,d\xi
\end{equation*}
for $x,x'\in (0,1)$, $\xi\in\R$, with 
$\tilde{p}(x,\xi,\lambda) = x^mp(x,x\xi,x^m\lambda)$.
\end{proposition}
\begin{proof}
The existence of a properly supported parametrix in
$\Psi^{-m}(\Lambda)$ follows immediately from Proposition
\ref{LokalcelliptParametrix}. We only need to verify that the
remainder term $C(\lambda)$ from equation
\eqref{AllgemeineOpsimKragen} can be arranged to vanish. Let first
\begin{equation*}
\tilde{Q}(\lambda)u(x) = \frac{1}{2\pi}\iint e^{i(x-x')\xi}
\tilde{q}(x,\xi,\lambda)u(x')\,dx'\,d\xi + C(\lambda)u(x)
\end{equation*}
be a parametrix of $A - \lambda$ in
$\Psi^{-m}(\Lambda)$, obtained by patching together local
parametrices from Proposition~\ref{LokalcelliptParametrix}, where
$\tilde{q}(x,\xi,\lambda) = x^m q(x,x\xi,x^m\lambda)$. We get the
desired $Q(\lambda)$ by setting
\begin{equation*}
p(x,\xi,\lambda) = \bigl({\mathcal F}_{x' \to \xi}\vp(x'){\mathcal
F}_{\xi \to x'}^{-1}q\bigr)(x,\xi,\lambda),
\end{equation*}
where $\mathcal F$ denotes the Fourier transform, and $\varphi \in
C_0^{\infty}(\R)$ is a function with $\varphi = 1$ in a 
neighborhood of the origin.
\end{proof}

\medskip 
We are finally ready to construct a parameter-dependent parametrix
$B_1(\lambda)$ of $A - \lambda$ on $M$. The
important aspect of the following theorem is the structure of the
complete symbol of $B_1(\lambda)$ close to the boundary of $M$.

\begin{theorem}\label{PEins}
Let $Q_{\mathrm{int}}(\lambda)$ be a standard parameter-dependent 
parametrix of $A - \lambda$ on $\open M$ which is properly
supported $($uniformly in $\lambda)$, and let
$Q(\lambda) \in \Psi^{-m}(\Lambda)$ be the parametrix of $A -
\lambda$ on $(0,1){\times}Y$ from
Proposition~\ref{KragencelliptParametrix}. Then for any cut-off
functions $\omega,\omega_0,\omega_1 \in C_0^{\infty}([0,1))$ with
$\omega_1 \prec \omega \prec \omega_0$, the properly supported
pseudodifferential operator
\begin{equation*}
B_1(\lambda)= {\omega}Q(\lambda){\omega_0} + 
 (1 -\omega)Q_{\mathrm{int}}(\lambda)(1 - {\omega_1})
\end{equation*}
is a parametrix of $A - \lambda$ on $M$.
\end{theorem}

Recall that a cut-off function $\omega \in C_0^{\infty}([0,1))$ is
a function which equals $1$ in a neighborhood of the origin. 
Observe that these functions can also be considered as functions on
$M$ supported in the collar neighborhood $[0,1){\times}Y$ of the
boundary. Moreover, we use the notation $\vp \prec \psi$ to
indicate that the function $\psi$ equals $1$ in a neighborhood of
the support of the function $\vp$, in particular, $\vp\psi =\vp$.

\medskip 
The second step in our parametrix construction concerns the refinement 
of $B_1(\lambda)$ from Theorem~\ref{PEins} to a Fredholm inverse.
First of all, we want to modify $B_1(\lambda)$ in order to get a family 
of bounded operators
\begin{equation*}
B_1(\lambda) : x^{-m/2}H^s_b(M;E) \to \Dom_{\min}^{s}(A)
\end{equation*}
for any $s \in \R$, where $\Dom_{\min}^{s}(A)$ denotes the
minimal domain of $A$ in $x^{-m/2}H^s_b(M;E)$, cf.
Section~\ref{sec-Domains}. Recall that for every $t\in\R$,
\begin{equation*}
 x^{m/2}H^{t+m}_b(M;E) \embed \Dom_{\min}^t \embed
 x^{-m/2+\eps}H^{t+m}_b(M;E).
\end{equation*}
Also, we use the notation $\Dom_{\min}(A) = \Dom_{\min}^0(A)$.

By Mellin quantization, one can easily modify $B_1(\lambda)$ in
such a way that
\begin{equation*}
B_1(\lambda) : x^{-m/2}H^s_b(M;E) \to x^{m/2}H^{s+m}_b(M;E)
\end{equation*}
is bounded for every $s \in \R$. Mellin representations
of pseudodifferential operators are standard. The following
proposition is a direct consequence of known results about the
Mellin quantization that can be found for instance in \cite{GSSOsaka}.

\begin{proposition}\label{MellinQuantisierung}
Let $Q(\lambda)$ be the parametrix of $A - \lambda$ from
Proposition~\ref{KragencelliptParametrix} defined via the symbol
$p(x,\xi,\lambda)$. Let
\begin{equation*}
h(x,\sigma,\lambda)= \frac{1}{2\pi}\iint
e^{-i(r-1)\xi} r^{i\sigma}\vp(r) p(x,\xi,\lambda)\,dr\,d\xi
\end{equation*}
for $r,x,\xi\in\R$, $\sigma\in\C$, where $\vp \in C_0^{\infty}(\R_+)$ 
is a function such that $\vp=1$ near $r = 1$. 
If we redefine $Q(\lambda)$ as
\begin{equation*}
Q(\lambda)u(x) = \frac{1}{2{\pi}i}\int\limits_{\Im\sigma= m/2\;}
\int\limits_{(0,1)}\Bigl(\frac{x}{x'}\Bigr)^{i\sigma}
x^m h(x,\sigma,x^m\lambda)u(x')\,\frac{dx'}{x'}\,d\sigma,
\end{equation*}
then the corresponding family $B_1(\lambda)$ from Theorem~\ref{PEins} 
is again a properly supported parametrix of $A-\lambda$
such that, in addition,  
\begin{equation*}
 B_1(\lambda): x^{-m/2}H^s_b(M;E) \to x^{m/2}H^{s+m}_b(M;E)
 \embed \Dom_{\min}^{s}(A)
\end{equation*}
is bounded for every $s \in \R$.
\end{proposition}

Our goal in this second step is to refine this parameter-dependent
parametrix in such a way that the remainders are elements of order zero
in a suitable class of Green operators that will be defined below.
To this end we consider scales of Hilbert spaces
$\{{\mathcal E}^s\}_{s \in \R}$ on $M$ and associated scales
$\{{\mathcal E}^{s,\delta}_{\wedge}\}_{s,\delta \in \R}$ on
$Y^{\wedge}$ as follows: Either ${\mathcal E}^s=
x^{\gamma}H^s_b(M;E)$ for some weight $\gamma \in \R$, or
${\mathcal E}^s= \Dom_{\min}^{s-m}(A)$. With the Sobolev spaces
${\mathcal E} = x^{\gamma}H$ we associate
$$
{\mathcal E}^{s,\delta}_{\wedge}=
{\omega}\bigl(x^{\gamma}H^s_b(Y^{\wedge};E)\bigr) +
(1 -\omega)\bigl(x^{\frac{n-m}{2}-\delta}
H^s_{\textup{cone}}(Y^{\wedge};E)\bigr),
$$
and for the scale of minimal domains ${\mathcal E} = \Dom_{\min}$
we define
$$
{\mathcal E}^{s,\delta}_{\wedge}=
{\omega}\Dom_{\min}^{s-m}(A_{\wedge}) +
(1-\omega)\bigl(x^{\frac{n-m}{2}-\delta}
H^s_{\textup{cone}}(Y^{\wedge};E)\bigr).
$$
Here $\omega \in C_0^{\infty}([0,1))$ denotes, as usual, a cut-off function near the origin. Note that in the latter case we have 
${\mathcal E}^{m,0}_{\wedge} = \Dom_{\min}(A_{\wedge})$. 
Recall that $n = \dim M$.

\begin{definition}\label{AdmissibleRemainders}
An operator family $G(\lambda) : C_0^{\infty}(\open M;E) \to
C^{\infty}(\open M;E)$ is called a \emph{Green remainder} of order
$\mu \in \R$ with respect to the scales $(\mathcal E,\mathcal F)$ if
for all cut-off functions $\omega,\tilde{\omega} \in C_0^{\infty}([0,1))$
the following holds:
\begin{enumerate}[$(i)$]
\item $(1-\omega)G(\lambda),\; G(\lambda)(1-\tilde{\omega}) \in
      \bigcap\limits_{s,t \in \R}
      \S (\Lambda,{\mathcal K}({\mathcal E}^s,{\mathcal F}^t))$,
\item $g(\lambda)= \omega{G(\lambda)}\tilde{\omega}:
  C_0^{\infty}(\open Y^{\wedge};E) \to C^{\infty}(\open Y^{\wedge};E)$
is a \emph{Green symbol}, i.e., a classical operator-valued symbol of
order $\mu \in \R$ in the following sense:
\begin{equation*}
g(\lambda) \in \bigcap\limits_{s,t,\delta,\delta' \in
{\R}}C^{\infty}(\Lambda,\K({\mathcal E}_{\wedge}^{s,\delta},
{\mathcal F}_{\wedge}^{t,\delta'})),
\end{equation*}
and for all multi-indices $\alpha \in \N_0^2$,
\begin{equation}\label{AdmissibleSymbol}
\Bigl\|{\kappa}_{[\lambda]^{1/m}}^{-1}\partial_{\lambda}^{\alpha}g(\lambda)
 {\kappa_{[\lambda]^{1/m}}}\Bigr\|_{
 \K({\mathcal E}_{\wedge}^{s,\delta},{\mathcal
 F}_{\wedge}^{t,\delta'})} = O(|\lambda|^{\mu/m - |\alpha|})
\end{equation}
as $|\lambda| \to \infty$. Here 
$\mathcal K({\mathcal E}^s,{\mathcal F}^t)$ denotes the space of
compact operators from ${\mathcal E}^s$ to ${\mathcal F}^t$, and
$[\cdot]$ is a strictly positive smoothing of the absolute value 
$|\cdot|$ near the origin. Without loss of generality we may assume 
$[\lambda]> 1$ for every $\lambda$.

Moreover, for $j \in \N_0$ there exist
\begin{equation*}
 g_{(\mu-j)}(\lambda) \in \bigcap\limits_{s,t,\delta,\delta' 
 \in \R}C^{\infty}(\Lambda \minus \{0\},
 \K({\mathcal E}_{\wedge}^{s,\delta},{\mathcal F}_{\wedge}^{t,\delta'}))
\end{equation*}
such that
\begin{equation*}
g_{(\mu-j)}(\varrho^m\lambda) =
\varrho^{\mu-j}\kappa_{\varrho}g_{(\mu-j)}(\lambda)\kappa_{\varrho}^{-1}
\quad\text{for } \varrho > 0,
\end{equation*}
and for some function $\chi \in C^{\infty}(\Lambda)$ with $\chi=0$
near zero and $\chi=1$ near $\infty$, and all $j \in {\mathbb
N}_0$, the symbol estimates \eqref{AdmissibleSymbol} hold for
$g(\lambda) - \sum_{k=0}^{j-1}\chi(\lambda)g_{(\mu-k)}(\lambda)$
with $\mu$ replaced by $\mu-j$.
\end{enumerate}
\end{definition}

As usual, the cut-off functions in $C_0^{\infty}([0,1))$ are
considered as functions on both $M$ and $Y^{\wedge}$, and 
$\{\kappa_{\varrho}\}_{\varrho \in \R_+}$ is the dilation group from
\eqref{DilationGroup}. 
The $\kappa$-homogeneous components $g_{(\mu-j)}(\lambda)$ are
well-defined for the Green remainder $G(\lambda)$, i.e., they do
not depend on the particular choice of cut-off functions (see also Lemma
\ref{Greensymbol} below). Hence a Green remainder is determined by
an asymptotic expansion
\begin{equation}\label{GreenAsymptotik}
G(\lambda) \sim \sum\limits_{j=0}^{\infty}G_{(\mu-j)}(\lambda)
\end{equation}
up to Green remainders of order $-\infty$, where
$G_{(\mu-j)}(\lambda) = g_{(\mu-j)}(\lambda)$.  The principal component
of $G(\lambda)$ in this expansion will be denoted by
\begin{equation*}
 G_{\wedge}(\lambda)= G_{(\mu)}(\lambda). 
\end{equation*}
Note that in view of Definition~\ref{AdmissibleRemainders}(i) every 
Green remainder $G(\lambda)$ is a parameter-dependent smoothing
pseudodifferential operator over the manifold $\open M$.

It should be pointed out that the choice of the compact operators as
operator ideal for the Green remainders is just for convenience; we
could also pass to the Schatten classes $\ell^p(\mathcal
E_{\wedge}^s,\mathcal F_{\wedge}^t)$ for arbitrary $p > 0$, or even to
s-nuclear operators in $\bigcap_{p >0}\ell^p(\mathcal
E_{\wedge}^s,\mathcal F_{\wedge}^t)$.  This is useful for applications
to index theory, especially the case of trace class remainders.

\begin{lemma}\label{Greensymbol}
Let $g(\lambda)$ be a Green symbol of order $\mu \in {\mathbb R}$,
and $\omega \in C_0^{\infty}(\overline{{\mathbb R}}_+)$ a cut-off
function near zero. Then $(1-\omega)g(\lambda)$ and
$g(\lambda)(1-\omega)$ are Green symbols of order $-\infty$, i.e.,
\begin{equation*}
(1-\omega)g(\lambda), \;\; g(\lambda)(1-\omega) \in 
\S(\Lambda,\K({\mathcal E}_{\wedge}^{s,\delta},
{\mathcal F}_{\wedge}^{t,\delta'})).
\end{equation*}
\end{lemma}
\begin{proof}
We only need to prove that
\begin{equation*}
(1-\omega)g(\lambda) = O([\lambda]^{-L})
\;\text{ as } |\lambda| \to \infty,
\;\text{ for all } L \in \R. 
\end{equation*}
The argument for higher derivatives and for $g(\lambda)(1-\omega)$ 
is analogous. 

Write $(1-\omega(x)) = \vp_k(x)x^k$ for every $k \in \N_0$. 
Note that $\vp_k \in C^{\infty}(\R_+)$ is
supported away from the origin, and $\vp_k(x) = \frac{1}{x^k}$ for
sufficiently large $x$. Then, for any given $s,t,\delta,\delta'\in\R$, 
and denoting the norms in $\L({\mathcal E}_{\wedge}^{s,\delta},{\mathcal
F}_{\wedge}^{t,\delta'})$ and $\L({\mathcal F}_{\wedge}^{t,\delta'})$
by $\|\cdot\|_{\delta,\delta'}$ and $\|\cdot\|_{\delta'}$,
respectively, we have
\begin{align*}
\Bigl\|{\kappa}^{-1}_{[\lambda]^{1/m}}
&(1-\omega)g(\lambda){\kappa}_{[\lambda]^{1/m}}\Bigr\|_{\delta,\delta'}\\
&=
\Bigl\|\vp_k\bigl(\tfrac{x}{[\lambda]^{1/m}}\bigr)[\lambda]^{-k/m}x^k
{\kappa}^{-1}_{[\lambda]^{1/m}}g(\lambda){\kappa}_{[\lambda]^{1/m}}
\Bigr\|_{\delta,\delta'} \\
&\leq C \Bigl\|\vp_k\bigl(\tfrac{x}{[\lambda]^{1/m}}\bigr)
\Bigr\|_{\delta'-k} \cdot
\Bigl\|{\kappa}^{-1}_{[\lambda]^{1/m}}g(\lambda)
{\kappa}_{[\lambda]^{1/m}}\Bigr\|_{\delta,\delta'-k}\cdot [\lambda]^{-k/m} \\
&\leq \tilde C \Bigl\|\vp_k\bigl(\tfrac{x}{[\lambda]^{1/m}}\bigr)
\Bigr\|_{\delta'-k} \cdot [\lambda]^{\frac{\mu-k}{m}}
\end{align*}
for some constants $C$ and $\tilde C$. As the norm of $\vp_k(x/
[\lambda]^{1/m})$ is $O(1)$ as $|\lambda| \to \infty$, the
assertion follows for $(1-\omega)g(\lambda)$. 
\end{proof}

A direct consequence from Lemma \ref{Greensymbol} is that the Green
remainders form an algebra. The homogeneous components of the product
of two Green remainders are determined by formally multiplying the
asymptotic sums \eqref{GreenAsymptotik}. In particular, 
\begin{equation*}
 (G_1G_2)_{\wedge}(\lambda) =
 G_{1,\wedge}(\lambda)G_{2,\wedge}(\lambda).
\end{equation*}

\begin{lemma}\label{GreenIdeal}
Let $G(\lambda)$ be a Green remainder of order $\mu \in \R$. Then
\begin{enumerate}[$(i)$]
\item $(A-\lambda)G(\lambda)$ and $G(\lambda)(A-\lambda)$ are
Green remainders of order $\mu+m$.
\item $B_1(\lambda)G(\lambda)$ and $G(\lambda)B_1(\lambda)$ are
Green remainders of order $\mu-m$.
\end{enumerate}
In all four cases the principal components are the composition of the principal components of the factors..
\end{lemma}
Recall that the principal component of $A-\lambda$
is $A_\wedge-\lambda$. On the other hand, the principal component
of $B_1(\lambda)$ is given by
\begin{equation}\label{PEinsPrinc}
B_{1,\wedge}(\lambda)u(x) =
 x^m\Bigl(\frac{1}{2{\pi}i}\Bigr)\int\limits_{\Im\sigma =
 m/2\;}\int\limits_{\R_+}
 \Bigl(\frac{x}{x'}\Bigr)^{i\sigma}
 h(0,\sigma,x^m\lambda)u(x')\,\frac{dx'}{x'}\,d\sigma
\end{equation}
for $u \in C_0^{\infty}(\R_+,C^{\infty}(Y;E))$, where 
$h(x,\sigma,\lambda)$ 
is the symbol from Proposition~\ref{MellinQuantisierung}.
For the above compositions to make sense, we are tacitly assuming 
that $G(\lambda)$ acts on corresponding scales.
\begin{proof}
Let us consider $(A - \lambda)G(\lambda)$. The product 
$G(\lambda)(A-\lambda)$ can be treated in a similar way.
In the collar neighborhood $(0,1){\times}Y$ we have
$$
A = x^{-m}\sum\limits_{j=0}^ma_j(x)(xD_x)^j,
$$
where $a_j(x) \in C^{\infty}([0,1),\Diff^{m-j}(Y;E))$. We set
$A_{(m)}(\lambda)= A_{\wedge} - \lambda$, and for $k \in \N$,
$$
A_{(m-k)}(\lambda)=
x^{-m+k}\sum_{j=0}^m\frac{1}{k!}\bigl(\partial_x^ka_j\bigr)(0)(xD_x)^j.
$$
Observe that for each $j$, $A_{(j)}(\lambda):
C_0^{\infty}(\open Y^{\wedge};E) \to C^{\infty}(\open Y^{\wedge};E)$, and
$$
\omega\Bigl((A-\lambda) -
\sum\limits_{k=0}^{N-1}A_{(m-k)}(\lambda)\Bigr)\tilde{\omega} \in
x^{-m+N}\Diff_b^m(Y^{\wedge};E)
$$
for any cut-off functions $\omega$, $\tilde{\omega} \in
C_0^{\infty}([0,1))$.

Let $\omega \in C_0^{\infty}([0,1))$ be an arbitrary cut-off
function. Then, as the operator norm of $A - \lambda$ grows
polynomially, it follows immediately
that $(A-\lambda)G(\lambda)(1-\omega)$ is rapidly decreasing in
$\Lambda$. On the other hand, using a suitable cut-off function
$\omega' \in C_0^{\infty}([0,1))$, we may write
$$
(1-\omega)(A - \lambda)G(\lambda) =
(1-\omega)(A-\lambda)(1-\omega')G(\lambda).
$$
Thus also $(1-\omega)(A-\lambda)G(\lambda)$ is rapidly decreasing
in $\Lambda$.

It remains to consider $\omega(A-\lambda)G(\lambda)\tilde{\omega}$
for cut-off functions $\omega,\;\tilde{\omega} \in
C_0^{\infty}([0,1))$. Choose cut-off functions $\omega_0$ and
$\omega_1$ such that $\omega \prec \omega_1 \prec \omega_0$. Then
\begin{align*}
   \omega(A-\lambda)G(\lambda)\tilde{\omega}
   &= \omega(A-\lambda)\omega_1\omega_0 G(\lambda)\tilde{\omega}\\
   &=\omega\Bigl(\sum_{k=0}^{N-1}A_{(m-k)}(\lambda)\Bigr)
     \omega_1\omega_0 G(\lambda)\tilde{\omega}
     + \omega{\tilde A_{N}}\omega_1\omega_0 G(\lambda)\tilde{\omega}
\end{align*}
for $N \in \N_0$, where $\tilde A_{N} \in
x^{-m+N}\Diff_b^{m}(Y^{\wedge};E)$. Since $g(\lambda) = \omega_0
G(\lambda)\tilde{\omega}$ is a Green symbol, it is easy to see
that $\omega{\tilde A_{N}}\omega_1 g(\lambda)$ is an operator-valued
symbol of order $\mu+m-N$, i.e., the estimates
\eqref{AdmissibleSymbol} hold with $\mu+m-N$ instead of $\mu$. The
argument here is to consider separately the terms
$\omega(x)\omega(x[\lambda]^{1/m }){\tilde A_{N}}\omega_1 g(\lambda)$
and $\omega(x)(1-\omega(x[\lambda]^{1/m})){\tilde A_{N}}\omega_1
g(\lambda)$.

Now, using the $\kappa$-homogeneity
$$
A_{(m-k)}(\varrho^m\lambda) =
\varrho^{m-k}{\kappa}_{\varrho}A_{(m-k)}(\lambda){\kappa}_{\varrho}^{-1}
$$
for $\varrho > 0$ and $\lambda \in \Lambda \minus \{0\}$, and because
of Lemma \ref{Greensymbol}, we finally conclude that $(A-\lambda)G(\lambda)$ 
is a Green remainder of order $\mu+m$.  Moreover, the homogeneous 
components of $(A-\lambda)G(\lambda)$ are given by
$$
\bigl((A-\lambda)G(\lambda)\bigr)_{(\mu+m-j)} =
\sum\limits_{k+l=j}A_{(m-k)}(\lambda)G_{(\mu-l)}(\lambda).
$$

The analysis for the products $G(\lambda)B_1(\lambda)$ and
$B_1(\lambda)G(\lambda)$ follows the same lines. At the
places where the locality of $(A - \lambda)$ was used, we can still
draw the desired conclusions for $B_1(\lambda)$, noting that for 
cut-off functions $\omega \prec \tilde{\omega}$ in $C_0^{\infty}([0,1))$, 
the operator families $\omega{B_1(\lambda)}(1-\tilde{\omega})$ and
$(1-\tilde{\omega}){B_1(\lambda)}\omega$ are Green remainders of
order $-\infty$. Moreover, on $Y^{\wedge}$ we expand
$B_1(\lambda)$ into components given by
\begin{equation*}
u \mapsto x^{m+k}\frac{1}{2{\pi}i}\int\limits_{\Im\sigma = m/2\,}
\int\limits_{\R_+} \Bigl(\frac{x}{x'}\Bigr)^{i\sigma}\frac{1}{k!}
\bigl(\partial_x^k h\bigr)(0,\sigma,x^m\lambda)u(x')\,\frac{dx'}{x'}\,d\sigma,
\quad k\in\N_0,
\end{equation*}
for $u \in C_0^{\infty}(\R_+,C^{\infty}(Y;E))$, and proceed as above.
\end{proof}

\begin{proposition}\label{GreenDminCharakt}
For an operator family
\begin{equation*}
 G(\lambda) : C_0^{\infty}(\open M;E) \to C^{\infty}(\open M;E)
\end{equation*}
the following are equivalent:
\begin{enumerate}[$(i)$]
\item $G(\lambda)$ is a Green remainder of order $\mu \in {\R}$ in 
the scales $({\mathcal E},\Dom_{\min})$.
\item $G(\lambda)$ is a Green remainder of order $\mu \in {\R}$ in 
the scales $({\mathcal E},x^{m/2-\varepsilon}H)$ for
every $\varepsilon > 0$, and $(A-\lambda)G(\lambda)$ is Green of
order $\mu + m$ in $({\mathcal E},x^{-m/2}H)$.
\end{enumerate}
\end{proposition}
\begin{proof}
The direction (i) $\Rightarrow$ (ii) follows from
Lemma~\ref{GreenIdeal} noting that
\begin{equation*}
 \Dom_{\min}^t(A) = \Dom_{\max}^t(A) \cap
 \biggl(\bigcap_{\eps > 0}x^{m/2-\eps}H^{t+m}_b(M;E)\biggr).
\end{equation*}
Let us now assume (ii). Then it is evident that for every cut-off
function $\omega \in C_0^{\infty}([0,1))$ the operator families
$(1-\omega)G(\lambda)$ and $G(\lambda)(1-\omega)$ are rapidly
decreasing in $\Lambda$ with values in the scale $\Dom_{\min}$ of
minimal domains. Hence it remains to consider
$\omega{G(\lambda)}\tilde{\omega}$ for cut-off functions
$\omega,\tilde{\omega} \in C_0^{\infty}([0,1))$.

Note first that the assertion of the proposition is obviously
valid at the level of Green symbols, i.e., $g(\lambda)$ is a Green
symbol of order $\mu \in \R$ with values in the
$\Dom_{\min}$-scale on $Y^{\wedge}$ if and only if $g(\lambda)$ is
a Green symbol of order $\mu \in \R$ with values in the
scale $x^{m/2-\varepsilon}H$ of Sobolev spaces on $Y^{\wedge}$ for
every $\varepsilon > 0$, and $(A_{\wedge}-\lambda)g(\lambda)$ is a
Green symbol of order $\mu + m$ with values in the scale
$x^{-m/2}H$ on $Y^{\wedge}$ (note that we are concerned with the
associated scales on $Y^{\wedge}$ in the sense of Definition
\ref{AdmissibleRemainders}).

Now let $\omega_0$ be another cut-off function such that 
$\omega \prec \omega_0$. Thus $\omega_0 \omega = \omega$ and so
\begin{equation*}
 (A_{\wedge} - \lambda)\bigl(\omega{G(\lambda)}\tilde{\omega}\bigr)
 = {\omega_0}(A-\lambda){\omega_0}
 \bigl(\omega{G(\lambda)}\tilde{\omega}\bigr) +
 \omega_0\tilde A \omega_0\bigl(\omega G(\lambda)\tilde{\omega}\bigr)
\end{equation*}
for some $\tilde A \in x^{-m+1}\Diff_b^m(Y^{\wedge};E)$. Hence 
$\omega_0 \tilde A\omega_0 \bigl(\omega{G(\lambda)}\tilde{\omega}\bigr)$ 
is a Green symbol of order $\mu+m-1$ with values in the scale
$x^{-m/2}H$ on $Y^{\wedge}$. Observe that this argument makes use
of our assumption that $G(\lambda)$ is a Green remainder of order
$\mu \in \R$ in the scales $({\mathcal E},x^{m/2-\eps}H)$ 
for every $\eps > 0$.

On the other hand, we may write
\begin{align*}
 {\omega_0}(A-\lambda){\omega_0}
 \bigl(\omega{G(\lambda)}\tilde{\omega}\bigr) 
 &= {\omega_0}(A-\lambda)\omega{G(\lambda)}\tilde{\omega} \\
 &= {\omega_0\omega}(A-\lambda)G(\lambda)\tilde{\omega} +
 {\omega_0}[(A-\lambda),\omega]G(\lambda)\tilde{\omega} \\ 
 &= {\omega}(A-\lambda)G(\lambda)\tilde{\omega} +
 {\omega_0}[(A-\lambda),\omega]G(\lambda)\tilde{\omega},
\end{align*}
where ${\omega_0}[(A-\lambda),\omega]G(\lambda)\tilde{\omega}$ is
rapidly decreasing in $\Lambda$. Thus we have proved
\begin{equation*}
 (A_{\wedge} - \lambda)\bigl(\omega{G(\lambda)}\tilde{\omega}\bigr)
 \equiv {\omega}(A-\lambda)G(\lambda)\tilde{\omega}
\end{equation*}
modulo a Green symbol of order $\mu+m-1$ with values in the scale
of Sobolev spaces $x^{-m/2}H$ on $Y^{\wedge}$, and as
${\omega}(A-\lambda)G(\lambda)\tilde{\omega}$ is a Green symbol
of order $\mu + m$ by our assumption (ii), the proposition follows.
\end{proof}

Let $\hat P_0(\sigma): C^{\infty}(Y;E|_Y) \to C^{\infty}(Y;E|_Y)$
be the conormal symbol of $A=x^{-m}P$, cf. \eqref{ConormalSymbol}. 
Since $A$ is assumed to be
$c$-elliptic, we know that the inverse $\hat{P}_0^{-1}(\sigma)$ of
$\hat P_0(\sigma)$ is a finitely meromorphic Fredholm function on $\C$,
and there exists a sufficiently small $\eps_0>0$
such that $\hat{P}_0(\sigma)$ is invertible in
$$
 \set{\sigma \in \C : -m/2-\eps_0 < \Im\sigma <
 -m/2+\eps_0,\; \Im\sigma \neq -m/2},
$$
with a holomorphic inverse there. Define
\begin{equation}\label{SmoothMellRest}
h_0(\sigma)= \hat{P}_0^{-1}(\sigma-im) - h(0,\sigma,0),
\end{equation}
where $h$ is the holomorphic Mellin symbol from Proposition
\ref{MellinQuantisierung}. Then $h_0(\sigma)$ is
finitely meromorphic in $\C$ taking values in
$L^{-\infty}(Y)$ and it is rapidly decreasing as $|\Re\sigma| 
\to \infty$, uniformly for $\Im\sigma$ in compact intervals.
Moreover, the strip
$$
 \set{\sigma \in \C : m/2-\eps_0 < \Im\sigma <
 m/2+\eps_0,\; \Im\sigma \neq m/2},
$$
is free of poles of $h_0(\sigma)$.

For arbitrary $0 < \eps < \eps_0$ and cut-off
function $\omega \in C_0^{\infty}([0,1))$ we define
\begin{equation*}
 M(\lambda): C_0^{\infty}(\open M;E)\to C^{\infty}(\open M;E)
\end{equation*} 
via 
\begin{equation*}
u \mapsto x^m\omega(x[\lambda]^{1/m})\biggl(\frac{1}{2{\pi}i}\!
\int\limits_{\Im\sigma = m/2+\eps\,} 
\int\limits_{\R_+}\Bigl(\frac{x}{x'}\Bigr)^{i\sigma} h_0(\sigma)
\omega(x'[\lambda]^{1/m})u(x')\, \frac{dx'}{x'}\,d\sigma\biggr) 
\end{equation*}
with the Mellin symbol $h_0(\sigma)$ from \eqref{SmoothMellRest}.
$M(\lambda)$ is a parameter-dependent smoothing
operator, and since the function $\omega(x[\lambda]^{1/m})$ is supported
in the collar $[0,1){\times}Y$, $M(\lambda)$ can be
regarded as an operator on both $M$ and $Y^\wedge$.

For $\lambda \neq 0$ we also define 
\begin{equation*}
 M_{\wedge}(\lambda): C_0^{\infty}(\open Y^{\wedge};E)
 \to C^{\infty}(\open Y^{\wedge};E)
\end{equation*} 
via
\begin{equation*}
 u \mapsto x^m\omega(x|\lambda|^{1/m})\biggl(\frac{1}{2{\pi}i}\!
 \int\limits_{\Im\sigma = m/2+\eps\,}\int\limits_{\R_+}
 \Bigl(\frac{x}{x'}\Bigr)^{i\sigma} h_0(\sigma)
 \omega(x'|\lambda|^{1/m})u(x')\,\frac{dx'}{x'}\,d\sigma\biggr).
\end{equation*}
Observe that $M_\wedge(\lambda)$ is $\kappa$-homogeneous of degree $-m$.

\begin{theorem}\label{PZwei}
Set $B_2(\lambda)= B_1(\lambda) + M(\lambda)$. Then
$$
B_2(\lambda): x^{-m/2}H^s_b(M;E) \to \Dom_{\min}^{s}(A)
$$
is a parameter-dependent parametrix of $A-\lambda$, and the remainders
\begin{align*}
G_1(\lambda) &= (A-\lambda)B_2(\lambda) - 1 : 
x^{-m/2}H^s_b(M;E) \to x^{-m/2}H^t_b(M;E),  \\
G_2(\lambda) &= B_2(\lambda)(A-\lambda) - 1 : 
\Dom_{\min}^{s}(A) \to \Dom_{\min}^{t}(A)
\end{align*}
are Green families of order zero in the sense of Definition
\ref{AdmissibleRemainders} with principal components 
given by
\begin{align*}
 G_{1,\wedge}(\lambda) = (A_{\wedge} -
 \lambda)B_{2,\wedge}(\lambda) - 1 \;\text{ and }\;
 G_{2,\wedge}(\lambda) = B_{2,\wedge}(\lambda)(A_{\wedge}-\lambda) - 1,
\end{align*}
where
\begin{equation}\label{PZweiPrinc}
B_{2,\wedge}(\lambda) =B_{1,\wedge}(\lambda)+M_{\wedge}(\lambda)
\end{equation}
with $B_{1,\wedge}(\lambda)$ as in \eqref{PEinsPrinc}.
\end{theorem}
\begin{proof}
Let us begin by noting that
$$
 B_2(\lambda) : x^{-m/2}H_b^s(M;E) \to 
 \bigcap_{\eps> 0}x^{m/2-\eps}H_b^{s+m}(M;E)
$$
is continuous. Hence, in order to show that $B_2(\lambda)$ maps
indeed into $\Dom_{\min}^{s}(A)$, it suffices to check that
$$
(A-\lambda)B_2(\lambda): x^{-m/2}H_b^s(M;E) \to x^{-m/2}H_b^s(M;E).
$$
We will prove that this operator is in fact of the form $1+G_1(\lambda)$.

By the standard composition rules for (parameter-dependent) cone
operators in cone Sobolev spaces (see e.g. \cite{GiHeat01},
\cite{GSSOsaka}, and \cite{SzWiley98}), we know that
\begin{equation*}
 (A-\lambda)B_1(\lambda) = 1 + \tilde{M}(\lambda) + G(\lambda),
\end{equation*}
where $G(\lambda)$ is a Green remainder of order zero in the
scales $(x^{-m/2}H,x^{-m/2}H)$, and $\tilde{M}(\lambda)$ is a
smoothing Mellin operator given by
$$
\tilde{M}(\lambda)u(x) =
\omega(x[\lambda]^{1/m})\biggl(\frac{1}{2{\pi}i}\!
\int\limits_{\Im\sigma = m/2\,}\int\limits_{\R_+}
\Bigl(\frac{x}{x'}\Bigr)^{i\sigma} \tilde{h}_0(\sigma)
\omega(x'[\lambda]^{1/m})u(x')\,\frac{dx'}{x'}\,d\sigma\biggr)
$$
with a holomorphic Mellin symbol 
\begin{equation}\label{SmoothMellinRest2}\tag{\ref{PZwei}a}
 \tilde{h}_0(\sigma) = \hat{P}_0(\sigma-im)h(0,\sigma,0) - 1
 = -\hat{P}_0(\sigma-im)h_0(\sigma)
\end{equation}
with $h_0$ as in \eqref{SmoothMellRest}.
Moreover, the principal components satisfy the identity
$$
(A_{\wedge}-\lambda)B_{1,\wedge}(\lambda) = 1 +
\tilde{M}_{\wedge}(\lambda) + G_{\wedge}(\lambda),
$$
where $\tilde{M}_{\wedge}(\lambda)$ is defined by
replacing $[\lambda]$ by $|\lambda|$ in $\tilde M(\lambda)$. 

Next we consider the composition $(A-\lambda)M(\lambda)$. As
$M(\lambda)$ is a Green remainder of order $-m$ in the scales
$(x^{-m/2}H,x^{m/2-\varepsilon}H)$ for every $\varepsilon > 0$, 
we conclude that up to a Green remainder of order $0$ in
$(x^{-m/2}H,x^{-m/2}H)$ we may write
\begin{align*}
(A-\lambda)M(\lambda) 
&\equiv {\omega_0}(x[\lambda]^{1/m})A_{\wedge}
 {\omega_0}(x[\lambda]^{1/m})M(\lambda) - \lambda{M(\lambda)}\\
&\equiv {\omega_0}(x[\lambda]^{1/m})A_{\wedge}
 {\omega_0}(x[\lambda]^{1/m})M(\lambda),
\end{align*}
where $\omega_0$ is a cut-off function with $\omega \prec \omega_0$,
so $\omega_0\omega=\omega$.  
Because of the relation \eqref{SmoothMellinRest2}, and since the 
commutator $[A_{\wedge},\omega(x[\lambda]^{1/m})] =
[A_{\wedge},\omega(x[\lambda]^{1/m})]{\omega_0}(x[\lambda]^{1/m})$
produces arbitrary flatness near the origin, we have 
\begin{equation*}
\omega_0 (x[\lambda]^{1/m})A_{\wedge}\omega_0(x[\lambda]^{1/m})
M(\lambda) \equiv -\tilde{M}(\lambda)
\end{equation*}
modulo a Green remainder of order zero in $(x^{-m/2}H,x^{-m/2}H)$. 

Hence we have proved that $(A-\lambda)M(\lambda) =
-\tilde{M}(\lambda) + \tilde{G}(\lambda)$ for some Green remainder 
$\tilde{G}(\lambda)$ of order zero in $(x^{-m/2}H,x^{-m/2}H)$. 
Consequently, 
\begin{equation*}
 (A-\lambda)B_2(\lambda) = 1 + G_1(\lambda) 
\end{equation*}
with $G_1(\lambda)=G(\lambda)+\tilde{G}(\lambda)$, and by
$\kappa$-homogeneity the principal components necessarily satisfy
$(A_{\wedge}-\lambda)B_{2,\wedge}(\lambda) =
1 + G_{1,\wedge}(\lambda)$.  Thus the assertion of the theorem regarding
the composition $(A-\lambda)B_2(\lambda)$ is proved.

It remains to investigate the composition
$B_2(\lambda)(A-\lambda)$. Again, we first apply the standard
composition rules of (parameter-dependent) cone operators in cone
Sobolev spaces to see that $B_2(\lambda)(A-\lambda) = 1 +
G_2(\lambda)$, where $G_2(\lambda)$ is a Green remainder of order
zero in the scales $(\Dom_{\min},x^{m/2-\eps}H)$ for
arbitrary $\eps > 0$. Moreover, the principal components
satisfy the desired identity
$B_{2,\wedge}(\lambda)(A_{\wedge}-\lambda) = 1 +
G_{2,\wedge}(\lambda)$. As $(A-\lambda)G_2(\lambda) =
G_1(\lambda)(A - \lambda)$, we obtain from Lemma \ref{GreenIdeal}
that $(A-\lambda)G_2(\lambda)$ is a Green remainder of order $m$
in $(\Dom_{\min},x^{-m/2}H)$. Proposition~\ref{GreenDminCharakt} 
now implies that $G_2(\lambda)$ is a Green
remainder of order zero in $(\Dom_{\min},\Dom_{\min})$. 
\end{proof}

\begin{remark}\label{PZweiEigenschaften}
The parametrix $B_2(\lambda)$ has the following properties.
\begin{enumerate}[(i)]
\item As a consequence of Theorem~\ref{PZwei}, for
$\lambda \in \Lambda \minus \{0\}$,
$$
A_{\wedge} - \lambda : \Dom_{\min}(A_{\wedge}) \to
x^{-m/2}L^2_b(Y^{\wedge};E)
$$
is Fredholm and $B_{2,\wedge}(\lambda)$ is a Fredholm inverse.
\item The principal component $B_{2,\wedge}(\lambda)$ is
$\kappa$-homogeneous of degree $-m$, i.e.,
$$
B_{2,\wedge}(\varrho^m\lambda) =
\varrho^{-m}{\kappa}_{\varrho}B_{2,\wedge}(\lambda){\kappa}_{\varrho}^{-1}
: C_0^{\infty}(\open Y^{\wedge};E) \to C^{\infty}(\open Y^{\wedge};E)
$$
for $\varrho > 0$ and $\lambda\in\Lambda\minus \{0\}$.
\item Let $G(\lambda)$ be a Green remainder of order $\mu \in \R$.
Then $B_2(\lambda)G(\lambda)$ and $G(\lambda)B_2(\lambda)$ are
Green remainders of order $\mu - m$, and the principal components
are given as $B_{2,\wedge}(\lambda)G_{\wedge}(\lambda)$ and
$G_{\wedge}(\lambda)B_{2,\wedge}(\lambda)$, respectively.
\item For every $s \in \R$ the following equivalent norm estimates hold:
\begin{align}
\label{PZweiSobolev}
\|B_2(\lambda)\|_{\L(x^{-m/2}H_b^s)}
&\leq \textup{const}\cdot[\lambda]^{2|s|/m-1}, \\
\label{PZweiDomain}
\|B_2(\lambda)\|_{\L(x^{-m/2}H_b^s,\Dom_{\min}^{s}(A))}
&\leq \textup{const}\cdot[\lambda]^{2|s|/m}.
\end{align}
If $G(\lambda)$ is an arbitrary Green remainder of order
$-m$, then $B_2(\lambda) + G(\lambda)$ is also an
admissible parameter-dependent parametrix of $A - \lambda$
satisfying the same norm estimates as $B_2(\lambda)$.
\end{enumerate}
\end{remark}
\begin{proof}
Let us prove (iii): By Lemma \ref{GreenIdeal} we only have to deal
with the terms $M(\lambda)G(\lambda)$ and $G(\lambda)M(\lambda)$.
But, since $M(\lambda): C_0^{\infty}(\open Y^{\wedge};E) \to
C^{\infty}(\open Y^{\wedge};E)$ satisfies
$$
M(\varrho^m\lambda) =
\varrho^{-m}{\kappa}_{\varrho}M(\lambda){\kappa}_{\varrho}^{-1}
$$
for $|\lambda| \gg 0$ and $\varrho \geq 1$, the assertion for
these terms is evident.

We now prove (iv): The group action $\{\kappa_{\varrho}\}_{\varrho
\in \R_+}$ satisfies the estimate
\begin{equation*}
 \bigl\|\kappa_{[\lambda]^{1/m}}\bigr\|_{\L(\K^{s,-m/2})}
  \leq \textup{const}\cdot[\lambda]^{|s|/m} 
\end{equation*}
on the space $\K^{s,-m/2}(Y^\wedge;E)$.
Recall that $\{\kappa_{\varrho}\}_{\varrho \in \R_+}$ is
defined to be unitary in $x^{-m/2}L^2_b(Y^{\wedge};E)$. Hence every
Green remainder $G(\lambda)$ of order zero in the scales
$(x^{-m/2}H,x^{-m/2}H)$ satisfies the norm estimate
$$
 \|G(\lambda)\|_{\L(x^{-m/2}H_b^s)} \leq
 \textup{const}\cdot[\lambda]^{2|s|/m}.
$$
Together with Theorem \ref{PZwei} this implies that the asserted
estimates are actually equivalent. Moreover, \eqref{PZweiSobolev}
follows from the estimates for the group action and the standard
estimates for parameter-dependent pseudodifferential operators in
Sobolev spaces, cf. Shubin \cite[Section 9]{Shubin}.
\end{proof}

\medskip 
As outlined at the beginning of this section, our goal is the
construction of a parametrix $B(\lambda)$ of $A-\lambda$ that is
a left-inverse for $\lambda$ sufficiently large. To achieve this,
we additionally require that the family
\begin{equation*}
 A_{\wedge} - \lambda : \Dom_{\min}(A_{\wedge}) \to
 x^{-m/2}L^2_b(Y^{\wedge};E)
\end{equation*}
be injective for all $\lambda \in \Lambda \minus \{0\}$.

In the remaining part of this section we will prove the following 
theorem:

\begin{theorem}\label{PDrei}
Let $B_2(\lambda)$ be the parametrix from Theorem~\ref{PZwei}.
Then there exists a Green remainder $G(\lambda)$ of order $-m$ in
the scales $(x^{-m/2}H,\Dom_{\min})$ such that
\[ B(\lambda)= B_2(\lambda) + G(\lambda) \]
is a parameter-dependent parametrix of $A - \lambda$ with
$B(\lambda)(A-\lambda) = 1$ for $\lambda$ sufficiently
large. In particular, for these values of $\lambda$, 
$(A-\lambda)B(\lambda)$ is a projection onto $\rg(A-\lambda)$, the 
range of 
\begin{equation*}
 A - \lambda : \Dom_{\min}^s(A) \to x^{-m/2}H_b^{s}(M;E),
\end{equation*}
and thus the Green remainder
\begin{equation*}
\Pi(\lambda)= 1 - (A-\lambda)B(\lambda)
\end{equation*}
is a projection onto some complement of $\rg(A-\lambda)$ in
$x^{-m/2}H_b^{s}(M;E)$ which is finite dimensional, is contained 
in $x^{-m/2}H_b^{\infty}(M;E)$, and is independent of $s$.
\end{theorem}

For the proof of this theorem we first introduce the following
class of generalized Green remainders.

\begin{definition}\label{GenGreenRem}
We consider scales of Hilbert spaces $\{{\mathcal E}^s\}_{s \in
\R}$ on $M$ and associated scales $\{{\mathcal
E}^{s,\delta}_{\wedge}\}_{s,\delta \in \R}$ on
$\open Y^{\wedge}$ as in Definition \ref{AdmissibleRemainders}.
Moreover, let $N_-,N_+ \in \N_0$.

An operator family
$$
G(\lambda) :
\begin{array}{c}
C_0^{\infty}(\open M;E) \\ \oplus \\ \C^{N_-}
\end{array} \to
\begin{array}{c}
C^{\infty}(\open M;E) \\ \oplus \\ \C^{N_+}
\end{array}
$$
is called a \emph{generalized Green remainder} of order $\mu \in
\R$ in the scales of spaces $\bigl({\mathcal E} \oplus
\C^{N_-},{\mathcal F} \oplus \C^{N_+}\bigr)$, if
for any cut-off functions $\omega,\tilde{\omega} \in
C_0^{\infty}([0,1))$ it holds:
\begin{enumerate}[(i)]
\item  For every $s,t\in\R$ the families
$$
\begin{pmatrix} (1-\omega) & 0 \\ 0 & 0 \end{pmatrix}G(\lambda)
\;\text{ and }\; G(\lambda)\begin{pmatrix} (1-\tilde{\omega}) & 0 \\
0 & 0
\end{pmatrix}
$$
are rapidly decreasing in $\Lambda$ with values in the compact
operators mapping
\begin{equation*}
\begin{array}{c}
\mathcal E^s \\ \oplus \\ \C^{N_-}
\end{array}
\to
\begin{array}{c}
\mathcal F^t \\ \oplus \\ \C^{N_+}
\end{array}\!.
\end{equation*}
\item The family $g(\lambda)$ given by
$$
g(\lambda)= \begin{pmatrix} \omega & 0 \\ 0 & 1
\end{pmatrix}{G(\lambda)}
\begin{pmatrix} \tilde{\omega} & 0 \\ 0 & 1 \end{pmatrix} :
\begin{array}{c}
C_0^{\infty}(\open Y^{\wedge};E) \\ \oplus \\ \C^{N_-}
\end{array} \to
\begin{array}{c}
C^{\infty}(\open Y^{\wedge};E) \\ \oplus \\ \C^{N_+}
\end{array}
$$
is a \emph{generalized Green symbol}, i.e., it is a classical
operator-valued symbol of order $\mu \in \R$ in the
following sense:
\begin{equation*}
g(\lambda) \in \bigcap\limits_{s,t,\delta,\delta' \in \R}
C^{\infty}\bigl(\Lambda,\K\bigl(
{\mathcal E}_{\wedge}^{s,\delta} \oplus \C^{N_-},
{\mathcal F}_{\wedge}^{t,\delta'} \oplus \C^{N_+} \bigr)\bigr),
\end{equation*}
and for all multi-indices $\alpha \in \N_0^2$,
\begin{equation}\label{GenGreenSymbol}
\biggl\|\begin{pmatrix} {\kappa}_{[\lambda]^{1/m}} & 0 \\ 0 & 1
\end{pmatrix}^{-1}
\partial_{\lambda}^{\alpha}g(\lambda)
\begin{pmatrix} \kappa_{[\lambda]^{1/m}} & 0 \\ 0 & 1 \end{pmatrix}\biggr\|
= O(|\lambda|^{\mu/m - |\alpha|})
\end{equation}
as $|\lambda| \to \infty$. Moreover, for $j \in \N_0$
there exist
\begin{equation*}
g_{(\mu-j)}(\lambda) \in \bigcap\limits_{s,t,\delta,\delta' \in \R}
C^{\infty}\bigl(\Lambda \minus \{0\},\K\bigl(
{\mathcal E}_{\wedge}^{s,\delta} \oplus \C^{N_-},
{\mathcal F}_{\wedge}^{t,\delta'} \oplus \C^{N_+} \bigr)\bigr),
\end{equation*}
such that
$$
g_{(\mu-j)}(\varrho^m\lambda) = \varrho^{\mu-j}
\begin{pmatrix} \kappa_{\varrho} & 0 \\ 0 & 1 \end{pmatrix}
g_{(\mu-j)}(\lambda)
\begin{pmatrix} \kappa_{\varrho} & 0 \\ 0 & 1 \end{pmatrix}^{-1}
$$
for every $\varrho > 0$, and for some function $\chi\in
C^{\infty}(\Lambda)$ with $\chi=0$ near zero and $\chi=1$ near
$\infty$, the symbol estimates \eqref{GenGreenSymbol} hold for
$g(\lambda) -
\sum\limits_{k=0}^{j-1}\chi(\lambda)g_{(\mu-k)}(\lambda)$ with
$\mu$ replaced by $\mu-j$.
\end{enumerate}
\end{definition}

Note that when $N_- = N_+ = 0$, we recover the class of
Green remainders from Definition~\ref{AdmissibleRemainders}.
Also for generalized Green remainders, the 
$\kappa$-homogeneous components $g_{(\mu-j)}(\lambda)$ are 
well-defined for $G(\lambda)$, i.e., they do not depend on the choice 
of the cut-off functions. Thus a generalized Green remainder is 
determined by an asymptotic expansion
\begin{equation}\label{GreenAsymptotikgen}
G(\lambda) \sim \sum\limits_{j=0}^{\infty}G_{(\mu-j)}(\lambda)
\end{equation}
up to generalized Green remainders of order $-\infty$, where
$G_{(\mu-j)}(\lambda) = g_{(\mu-j)}(\lambda)$. The principal
component will again be denoted by 
$G_{\wedge}(\lambda)= G_{(\mu)}(\lambda)$.

We will be particularly concerned with the operators
$$
\begin{pmatrix} A - \lambda & 0 \\ 0 & 0 \end{pmatrix} + G(\lambda), \quad
\begin{pmatrix} B_2(\lambda) & 0 \\ 0 & 0 \end{pmatrix} + G'(\lambda)
$$
for generalized Green remainders $G(\lambda)$ and $G'(\lambda)$ of
order $m$ and $-m$, respectively. We will also need their
$\kappa$-homogeneous principal components
$$
\begin{pmatrix} A_{\wedge} - \lambda & 0 \\ 0 & 0 \end{pmatrix} 
+ G_{\wedge}(\lambda), \quad
\begin{pmatrix} B_{2,\wedge}(\lambda) & 0 \\ 0 & 0 \end{pmatrix}
+ G'_{\wedge}(\lambda).
$$
Lemma~\ref{GreenIdeal} (as well as (iii) in Remark
\ref{PZweiEigenschaften}) continues to hold in this more general
framework, and Theorem~\ref{PZwei} implies
\begin{align*}
\biggl(\begin{pmatrix} A - \lambda & 0 \\ 0 & 0 \end{pmatrix} +
G(\lambda)\biggr) \biggl(\begin{pmatrix} B_2(\lambda) & 0 \\ 0 & 0
\end{pmatrix} + G'(\lambda)\biggr)
&= 1 + G_1(\lambda), \\
\biggl(\begin{pmatrix} B_2(\lambda) & 0 \\ 0 & 0 \end{pmatrix} +
G'(\lambda)\biggr) \biggl(\begin{pmatrix} A - \lambda & 0 \\ 0 & 0
\end{pmatrix} + G(\lambda)\biggr) &= 1 + G_2(\lambda)
\end{align*}
with generalized Green remainders $G_1(\lambda)$ and
$G_2(\lambda)$ of order zero, provided the scales are such that
the composition makes sense. Moreover, the principal components
satisfy the same relations.

\begin{lemma}\label{EinsPlusGreen}
Let $G(\lambda)$ be a generalized Green remainder of order zero in
the scales $({\mathcal E} \oplus \C^N,{\mathcal E} \oplus
\C^N)$ for some $N \in \N_0$, and assume that
$$
1 + G_{\wedge}(\lambda) : \begin{array}{c} {\mathcal
E}_{\wedge}^{s,\delta} \\ \oplus \\ \C^N \end{array} \to
\begin{array}{c} {\mathcal E}_{\wedge}^{s,\delta} \\ \oplus \\
\C^N \end{array}
$$
is invertible for all $\lambda \in \Lambda \minus \{0\}$ and
some $s,\delta \in \R$. Then there exists a generalized
Green remainder $\tilde{G}(\lambda)$ of order zero such that
\[(1+ G(\lambda))(1 + \tilde{G}(\lambda)) - 1 \;\text{ and }\; (1 +
\tilde{G}(\lambda))(1 + G(\lambda)) - 1 \] are generalized Green
remainders of order $-\infty$. Moreover, $\tilde{G}(\lambda)$ can
be arranged in such a way that these remainders are compactly
supported in $\Lambda$, thus $(1 + \tilde{G}(\lambda))$ inverts
$(1 + G(\lambda))$ for every $\lambda$ sufficiently large.
\end{lemma}
\begin{proof}
The inverse of $1 + G_{\wedge}(\lambda)$ can be written as
$$
(1 + G_{\wedge}(\lambda))^{-1} = 1 + \tilde{G}_{\wedge}(\lambda)
$$
where $\tilde{G}_{\wedge}(\lambda)=
G_{\wedge}(\lambda)(1+G_{\wedge}(\lambda))^{-1}
G_{\wedge}(\lambda)- G_{\wedge}(\lambda)$ is a homogeneous Green
symbol of order zero. For $\lambda \in \Lambda$ set
$$
G'(\lambda)=
\begin{pmatrix} \omega & 0 \\ 0 & 1
\end{pmatrix}\chi(\lambda)\tilde{G}_{\wedge}(\lambda)
\begin{pmatrix} \omega & 0 \\ 0 & 1 \end{pmatrix},
$$
where $\omega \in C_0^{\infty}([0,1))$ is a cut-off function and
$\chi\in C^\infty(\Lambda)$ is a function with $\chi=0$ near 0 and
$\chi=1$ near $\infty$. Hence $G'(\lambda)$ is a generalized Green
remainder of order zero, and by construction we obtain
$$
(1 + G(\lambda))(1 + G'(\lambda)) = 1 + \tilde{G}_1(\lambda),
\quad (1 + G'(\lambda))(1 + G(\lambda)) = 1 + \tilde{G}_2(\lambda)
$$
with generalized Green remainders $\tilde{G}_1(\lambda)$ and
$\tilde{G}_2(\lambda)$ of order $-1$.

As the class of generalized Green remainders is asymptotically
complete, there exists a generalized Green remainder
$\tilde{G}_R(\lambda)$ of order $-1$ with
$$
\tilde{G}_R(\lambda) \sim
\sum\limits_{k=1}^{\infty}(-1)^k\tilde{G}_1^k(\lambda).
$$
This asymptotic expansion holds up to generalized Green remainders
of order $-\infty$. Hence
$$
(1 + G(\lambda))(1 + G'(\lambda))(1 + \tilde{G}_R(\lambda)) = 1 +
\tilde{G}_{(-\infty)}(\lambda)
$$
with a generalized Green remainder
$\tilde{G}_{(-\infty)}(\lambda)$ of order $-\infty$. In particular, 
the operator norm of $\tilde{G}_{(-\infty)}(\lambda)$
is decreasing as $|\lambda| \to \infty$ and therefore 
$1 + \tilde{G}_{(-\infty)}(\lambda)$ is invertible for $\lambda$ large. 
Moreover, the inverse can be written as
$$
\bigl(1 + \tilde{G}_{(-\infty)}(\lambda)\bigr)^{-1} = 1 +
\tilde{G}^{(-\infty)}(\lambda),
$$
where $ \tilde{G}^{(-\infty)}(\lambda)=
\tilde{G}_{(-\infty)}(\lambda)\bigl(1 +
\tilde{G}_{(-\infty)}(\lambda)\bigr)^{-1}\tilde{G}_{(-\infty)}(\lambda)
- \tilde{G}_{(-\infty)}(\lambda)$. Note that if $\chi\in
C^\infty(\Lambda)$ is a suitable function with $\chi=0$ near 0 and
$\chi=1$ near $\infty$, then
$\chi(\lambda)\tilde{G}^{(-\infty)}(\lambda)$ is a generalized
Green remainder of order $-\infty$. Summing up, we have proved
that
$$
(1 + G(\lambda))(1 + G'(\lambda))(1 + \tilde{G}_R(\lambda))(1 +
\chi(\lambda)\tilde{G}^{(-\infty)}(\lambda))- 1
$$
is compactly supported in $\Lambda$. Finally, we define
$\tilde{G}(\lambda)$ by
$$
1 + \tilde{G}(\lambda) = (1 + G'(\lambda))(1 +
\tilde{G}_R(\lambda))(1 +
\chi(\lambda)\tilde{G}^{(-\infty)}(\lambda)).
$$
By construction, $\tilde{G}(\lambda)$ is a generalized Green
remainder  of order zero, and $1 + \tilde{G}(\lambda)$ inverts $1
+ G(\lambda)$ from the right for large values of $\lambda$.

In the same way, we can prove that $1 + G(\lambda)$ has a left-inverse for $\lambda$ sufficiently large. This inverse must be necessarily $1 + \tilde{G}(\lambda)$ and the lemma is proved.
\end{proof}

The following theorem implies Theorem~\ref{PDrei}.

\begin{theorem}\label{AminuslambdaplusK}
For $\lambda\in\Lambda\minus\{0\}$
let $d''=-\Ind (A_{\wedge,\Dom_{\min}} -\lambda)$,
There exists a generalized Green remainder $\begin{pmatrix} 0 &
K(\lambda)\end{pmatrix}$ of order $m$ in the scales $(\Dom_{\min}
\oplus \C^{d''},x^{-m/2}H)$ such that
$$
\begin{pmatrix} A - \lambda & K(\lambda) \end{pmatrix} :
\begin{array}{c}
\Dom_{\min}^s(A) \\ \oplus \\ \C^{d''}
\end{array} \to x^{-m/2}H_b^{s}(M;E)
$$
is invertible for $\lambda$ sufficiently large. Moreover,
the inverse can be written as
$$
\begin{pmatrix} A - \lambda & K(\lambda) \end{pmatrix}^{-1} =
\begin{pmatrix} B_2(\lambda) + G(\lambda) \\ T(\lambda) \end{pmatrix},
$$
where $\begin{pmatrix} G(\lambda) \\ T(\lambda) \end{pmatrix}$ is
a generalized Green remainder of order $-m$ in the corresponding
scales $(x^{-m/2}H,\Dom_{\min} \oplus \C^{d''})$. In particular,
the parameter-dependent parametrix
\begin{equation*}
 B(\lambda)= B_2(\lambda) +  G(\lambda)
\end{equation*}
satisfies the conditions of Theorem \ref{PDrei}.
\end{theorem}
\begin{proof}
 From Theorem~\ref{FredFam} (see also Remark~\ref{DenseSubspace}
and Corollary~\ref{ExtraCondSector}) we conclude that there exists
$k_{\wedge}(\lambda)$ such that
$$
\begin{pmatrix} A_{\wedge} - \lambda & k_{\wedge}(\lambda) \end{pmatrix} :
\begin{array}{c} \Dom_{\min}(A_{\wedge}) \\ \oplus \\
\C^{d''} \end{array} \to x^{-m/2}L_b^2(Y^{\wedge};E)
$$
is invertible for $\lambda \in \Lambda \minus \{0\}$, and
$k_{\wedge}(\lambda)$ can be arranged to be a homogeneous
principal Green symbol of order $m$.

Let $\omega \in C_0^{\infty}([0,1))$ be a cut-off function and let
$\chi\in C^\infty(\Lambda)$ be a function with $\chi=0$ near 0 and
$\chi=1$ near $\infty$. If we set
$K(\lambda)= \omega\chi(\lambda)k_{\wedge}(\lambda)$, then
$\begin{pmatrix} 0 & K(\lambda) \end{pmatrix}$ is a generalized Green 
remainder of order $m$.  We will prove that the theorem holds with this 
particular choice for $K(\lambda)$.

As $B_{2,\wedge}(\lambda)$ is a Fredholm inverse of $A_{\wedge}
- \lambda$ for $\lambda\in\Lambda \minus \{0\}$, we may apply
once again the results from Appendix~\ref{sec-ExtraConditions} to
conclude the existence of families $\tilde{k}_{\wedge}(\lambda)$,
$\tilde{t}_{\wedge}(\lambda)$, and $\tilde{q}_{\wedge}(\lambda)$
such that
\begin{equation*}
\begin{pmatrix}
  B_{2,\wedge}(\lambda)   
  & \tilde{k}_{\wedge}(\lambda) \\ 
  \tilde{t}_{\wedge}(\lambda) & \tilde{q}_{\wedge}(\lambda)
\end{pmatrix}:
\begin{array}{c} 
  x^{-m/2}L_b^2(Y^{\wedge};E) \\ \oplus \\ \C^{N_-} 
\end{array} \to
\begin{array}{c} 
  \Dom_{\min}(A_{\wedge}) \\ \oplus \\ \C^{N_+}  
\end{array}
\end{equation*}
is invertible for $\lambda \in \Lambda \minus \{0\}$, and
$\begin{pmatrix} 
  {\bf 0} & \tilde{k}_{\wedge}(\lambda) \\ 
  \tilde{t}_{\wedge}(\lambda) & \tilde{q}_{\wedge}(\lambda)
\end{pmatrix}$
is a homogeneous principal Green symbol of order $-m$. Note that
by construction $N_+ - N_- = \Ind B_{2,\wedge}(\lambda) = d''$. According to
$\C^{N_+}=\C^{d''}\oplus \C^{N_-}$ we
arbitrarily decompose 
\begin{equation*}
\tilde t_\wedge(\lambda)=
\begin{pmatrix} \tilde t_{\wedge,1}(\lambda)\\ 
\tilde t_{\wedge,2}(\lambda) \end{pmatrix} \;\text{ and }\;
\tilde q_\wedge(\lambda)=
\begin{pmatrix} \tilde q_{\wedge,1}(\lambda)\\ 
\tilde q_{\wedge,2}(\lambda) \end{pmatrix},
\end{equation*}
and let
\begin{equation*}
G'(\lambda)=
\begin{pmatrix} \omega &0 &0 \\ 0 &1 &0\\ 0 &0 &1\end{pmatrix}
\chi(\lambda)
\begin{pmatrix}
  {\bf 0} & \tilde{k}_{\wedge}(\lambda) \\ 
  \tilde{t}_{\wedge,1}(\lambda) & \tilde{q}_{\wedge,1}(\lambda) \\
  \tilde{t}_{\wedge,2}(\lambda) & \tilde{q}_{\wedge,2}(\lambda)
\end{pmatrix}
\begin{pmatrix} \omega & 0 \\ 0 & 1 \end{pmatrix},
\end{equation*}
where $\omega$ and $\chi$ are as above. Then $G'(\lambda)$ is a
generalized Green remainder of order $-m$ in the scales
$(x^{-m/2}H \oplus \C^{N_-},\Dom_{\min} \oplus \C^{N_+})$.
We now let
\begin{equation*}
\A(\lambda) =
\left(\!\begin{array}{cc|c}
\! A - \lambda & K(\lambda)\! & 0 \\ \hline
\Vsp 0 & 0 & [\lambda] 
\end{array}\!\right)
\;\text{ and }\;
\B(\lambda) =
\left(\!\begin{array}{c|c} 
  B_2(\lambda) & 0 \\ 0 & 0 \\ \hline \Vsp 0 & 0 
\end{array}\!\right) 
+ G'(\lambda),
\end{equation*}
and consider the compositions
\begin{align*}
\A(\lambda)\B(\lambda) &= 1 + G_1(\lambda) \;\text{ on } 
x^{-m/2}L_b^2(M;E)\oplus \C^{N_-},  \\
\B(\lambda)\A(\lambda) &= 1 + G_2(\lambda) \;\text{ on }
\bigl(\Dom_{\min}(A)\oplus \C^{d''}\bigr)\oplus \C^{N_-}.
\end{align*}
Note that 
$\left(\!\begin{array}{cc|c}
0 & K(\lambda)\! & 0 \\ \hline \Vsp 0 & 0 & [\lambda] 
\end{array}\!\right)$
is a generalized Green remainder of order $m$
with principal component
\begin{equation*}
\left(\!\begin{array}{cc|c}
0 & k_\wedge(\lambda)\! & 0 \\ \hline \Vsp 0 & 0 & |\lambda| 
\end{array}\!\right).
\end{equation*}
Hence $G_1(\lambda)$ and $G_2(\lambda)$ are generalized Green
remainders of order zero, and by construction both $1 +
G_{1,\wedge}(\lambda)$ and $1 + G_{2,\wedge}(\lambda)$ are
invertible for $\lambda \in \Lambda \minus \{0\}$.

Lemma~\ref{EinsPlusGreen} now implies the invertibility of 
$\A(\lambda)$ for $\lambda$ large. Consequently, the 
diagonal matrix structure of $\A(\lambda)$ gives the invertibility 
of $\begin{pmatrix} A - \lambda & K(\lambda) \end{pmatrix}$.
Moreover, 
\begin{equation*}
\A(\lambda)^{-1}= 
\left(\!\begin{array}{cc|c}
\! A - \lambda & K(\lambda)\! & 0 \\ \hline
\Vsp 0 & 0 & [\lambda]
\end{array}\!\right)^{\!-1} =
\B(\lambda)(1+ \tilde G(\lambda))
\end{equation*}
for some generalized Green remainder $\tilde{G}(\lambda)$ of order 
$-m$. Thus  
$\begin{pmatrix} A - \lambda & K(\lambda) \end{pmatrix}^{-1}$
must be of the form 
\begin{equation*}
\begin{pmatrix} B_2(\lambda)+ G(\lambda)\\ T(\lambda) \end{pmatrix}
\end{equation*}
which proves the theorem.
\end{proof}

\begin{corollary}
 For $\lambda\in\Lambda\minus \{0\}$ we have
 $\Ind(A_{\wedge,\Dom_{\min}}-\lambda) = \Ind A_{\Dom_{\min}}$.
\end{corollary}

As stated above, the parameter-dependent family 
$B(\lambda)= B_2(\lambda) +  G(\lambda)$ is a parametrix of 
$(A-\lambda)$ satisfying the conditions of Theorem~\ref{PDrei}. 
Let us draw some consequences of that theorem.

\begin{corollary}\label{MeromorphLeftInverse}
There exists a discrete set $\Delta \subset \C$ such that
$$
A - \lambda : \Dom_{\min}^s(A) \to x^{-m/2}H_b^{s}(M;E)
$$
is injective for $\lambda \in \C \minus \Delta$, and there
exists a finitely meromorphic left-inverse.
\end{corollary}
\begin{proof}
Due to Theorem~\ref{PDrei},
$$
A - \lambda : \Dom_{\min}^s(A) \to x^{-m/2}H_b^{s}(M;E)
$$
is injective for $\lambda \in \Lambda$ sufficiently large,
and the parametrix $B(\lambda)$ is a left-inverse.

Fix some large $\lambda_0 \in \Lambda$ and consider
the operator function
$$
F: \C\ni\lambda\mapsto B(\lambda_0)(A-\lambda)\in \L(\Dom_{\min}^s(A)).
$$
Then $F$ is a holomorphic Fredholm family on $\C$, and
$F(\lambda_0) = 1$ is invertible. The well known theorem on the
inversion of holomorphic Fredholm families now implies that the
inverse $\C \minus \Delta \ni \lambda \mapsto
F(\lambda)^{-1}$ is a finitely meromorphic operator function,
where $\Delta \subset \C$ is discrete. Hence $A - \lambda$ is
injective for $\lambda \in \C \minus \Delta$, and
$F(\lambda)^{-1}B(\lambda_0)$ is a finitely meromorphic left-inverse.
\end{proof}

\begin{corollary}\label{InverseDarstellung}
Let $\lambda_0 \in \Lambda$ and assume there exists some
domain $\Dom^s$ such that
\begin{equation*}\label{Aminuslambda0}
A - \lambda_0 : \Dom^s \to x^{-m/2}H_b^{s}(M;E)
\end{equation*}
is invertible. Then it is invertible for all $s\in\R$, and we have
\begin{equation*}
(A-\lambda_0)^{-1} = B(\lambda_0) + (A-\lambda_0)^{-1}\Pi(\lambda_0)
\end{equation*}
with the parametrix $B(\lambda)$ and the projection $\Pi(\lambda)$
from Theorem~\ref{PDrei}.
\end{corollary}

\section{Resolvents}
\label{sec-Resolvents}

The elements of the quotient 
\begin{equation*}
	\tilde\Sing_{\max}=\Dom_{\max}/\Dom_{\min}
\end{equation*}
can be conveniently identified with singular functions as follows.
Let $u \in \Dom_{\max}$. Then there is a finite sum of the form
\begin{equation}\label{Singulaerfunktion}
\tilde u = \!\!\sum_{-\frac{m}{2} < \Im(\sigma) < \frac{m}{2}}
\Bigl(\sum_{k=0}^{m_{\sigma}}c_{\sigma,k}(y)\log^kx\Bigr)x^{i\sigma}
\end{equation}
with $c_{\sigma,k}(y) \in C^{\infty}(Y;E)$ such that
$u - {\omega}\tilde u \in \Dom_{\min}$, where 
$\omega \in C_0^{\infty}([0,1))$ is a cut-off function near zero.
The function $\tilde u$ is uniquely determined by the equivalence class 
$u + \Dom_{\min}$, and in this way we may identify $\tilde\Sing_{\max}$ 
with a finite dimensional subspace of $C^{\infty}(\open Y^\wedge;E)$ 
consisting of singular functions \eqref{Singulaerfunktion}.
Analogously, we also obtain an identification of 
\begin{equation*}
 \tilde\Sing_{\wedge,\max}=\Dom_{\wedge,\max}/\Dom_{\wedge,\min}
\end{equation*}
with a finite dimensional space of functions of the form 
\eqref{Singulaerfunktion}.

In order to prove the existence of sectors of minimal growth for a given
extension $A_{\Dom}$, we are led to consider a particular extension 
$A_{\wedge,\Dom_{\wedge}}$ of the model operator. Thereby, the domain 
$\Dom_{\wedge}$ is associated to $\Dom$ via
\begin{equation}\label{DwedgeIdentitaet}
\Dom_{\wedge} / \Dom_{\wedge,\min} = \theta\bigl(\Dom / \Dom_{\min}\bigr),
\end{equation}
where
$$
\theta: \tilde\Sing_{\max} \to \tilde\Sing_{\wedge,\max}
$$
is the natural isomorphism introduced in \cite{GilKrainerMendoza1}.

Using the identification of the quotients with spaces of singular functions,
we briefly recall the definition of $\theta$. To this end, we split
\begin{equation}\label{APhiDarstellung}
A = x^{-m}\sum_{k=0}^{m-1}P_kx^k + \tilde P_m
\end{equation}
near $Y$, where each $P_k \in \Diff_b^m(Y^{\wedge};E)$ has coefficients 
independent of $x$, and $\tilde P_m \in \Diff_b^m(Y^{\wedge};E)$. Let 
$\hat{P}_k(\sigma)$ be the conormal symbol associated with $P_k$.
In this section, all arguments involving \eqref{APhiDarstellung} will refer 
to functions that are supported near $Y$, so we may assume that the 
coefficients of $\tilde P_m$ vanish near infinity.
In slight abuse of the notation from \cite{GilKrainerMendoza1} we now write
\begin{align*}
\tilde\Sing_{\max}= \bigoplus_{\sigma_0 \in \Sigma}\tilde\Sing_{\sigma_0}
\;\text{ and }\;
\tilde\Sing_{\wedge,\max}= \bigoplus_{\sigma_0 \in \Sigma}
\tilde\Sing_{\wedge,\sigma_0},
\end{align*}
where
\begin{equation}\label{SigmaDefinition}
\Sigma= \spec_b(A) \cap \{\sigma \in \C : -m/2 < \Im(\sigma) < m/2\}.
\end{equation}

The space $\tilde\Sing_{\wedge,\sigma_0}$ consists of all singular functions 
of the form
$$
\Bigl(\sum_{k=0}^{m_{\sigma_0}}c_{\sigma_0,k}(y)\log^kx\Bigr)x^{i\sigma_0}
$$
that are associated with elements of $\tilde\Sing_{\wedge,\max}$.
The operator $\theta$ acts isomorphically between $\tilde\Sing_{\sigma_0} \to
\tilde\Sing_{\wedge,\sigma_0}$. Both, the space $\tilde\Sing_{\sigma_0}$
and the operator itself, are easiest understood from its inverse
\begin{equation}\label{thetainvers}
\theta^{-1}|_{\tilde\Sing_{\wedge,\sigma_0}} = \sum_{k=0}^{N(\sigma_0)}
\e_{\sigma_0,k} : \tilde\Sing_{\wedge,\sigma_0} \to \tilde\Sing_{\sigma_0},
\end{equation}
where $N(\sigma_0) \in {\mathbb N}_0$ is the largest integer such that 
$\Im \sigma_0 - N(\sigma_0) > -m/2$, and the operators 
\begin{equation*}
\e_{\sigma_0,k} : \tilde\Sing_{\wedge,\sigma_0} \to 
C^{\infty}(\open Y^\wedge;E)
\end{equation*}
are inductively defined as follows:
\begin{itemize}
\item $\e_{\sigma_0,0} = \id$, the identity map.
\item Given $\e_{\sigma_0,0},\ldots,\e_{\sigma_0,\vartheta-1}$ for some
$\vartheta \in \{1,\ldots,N(\sigma_0)-1\}$, we define 
$\e_{\sigma_0,\vartheta}(\psi)$ for $\psi \in \tilde\Sing_{\wedge,\sigma_0}$ 
to be the unique singular function of the form
$$
\Bigl(\sum_{k=0}^{m_{\sigma_0-i\vartheta}}c_{\sigma_0-i\vartheta,k}(y)
\log^kx\Bigr)x^{i(\sigma_0-i\vartheta)}
$$
such that
\begin{equation*}
\qquad
(\omega \e_{\sigma_0,\vartheta}(\psi))^{\wedge}(\sigma) + 
\hat{P}_0(\sigma)^{-1}\Bigl(\sum_{k=1}^{\vartheta}\hat{P}_k(\sigma)
\s_{\sigma_0-i\vartheta}(\omega \e_{\sigma_0,\vartheta-k}(\psi))^{\wedge}
(\sigma+ik)\Bigr)
\end{equation*}
is holomorphic at $\sigma = \sigma_0 - i\vartheta$, where 
$(\omega \e_{\sigma_0,\vartheta-k}(\psi))^{\wedge}(\sigma)$ is
the Mellin transform of the function $\omega \e_{\sigma_0,\vartheta-k}(\psi)$, 
and $\s_{\sigma_0-i\vartheta}(\omega \e_{\sigma_0,\vartheta-k}(\psi))^{\wedge}
(\sigma+ik)$ is the singular part of the Laurent expansion at $\sigma_0 - 
i\vartheta$. Here, $\omega \in C_0^{\infty}(\overline{\R}_+)$ is an arbitrary 
cut-off function near zero.
Recall that the Mellin transform of $\omega \e_{\sigma_0,\vartheta-k}(\psi)$ 
is meromorphic in $\C$ with only one pole at $\sigma_0 - i(\vartheta-k)$.
\end{itemize}

It is of interest to note that this construction yields
$$
\sum_{k=0}^{\vartheta}\bigl(P_kx^k\bigr)(\e_{\sigma_0,\vartheta-k}(\psi))=0
$$
for every $\psi \in \tilde\Sing_{\wedge,\sigma_0}$ and every
$\vartheta = 0,\ldots,N(\sigma_0)$.

In conclusion, every space $\tilde\Sing_{\sigma_0}$ consists of singular 
functions of the form
\begin{equation*}
\tilde u = \sum_{\vartheta=0}^{N(\sigma_0)}
\Bigl(\sum_{k=0}^{m_{\sigma_0-i\vartheta}}c_{\sigma_0-i\vartheta,k}(y)
\log^kx\Bigr)x^{i(\sigma_0-i\vartheta)},
\end{equation*}
and we have
\begin{equation} \label{ThetaOperator}
\theta\tilde u = \Bigl(\sum_{k=0}^{m_{\sigma_0}}c_{\sigma_0,k}(y)
\log^kx\Bigr)x^{i\sigma_0}.
\end{equation}

The main result of this section concerns
the existence of sectors of minimal growth for closed extensions of a 
$c$-elliptic cone operator $A$. Recall that a sector
\begin{equation*}
\Lambda = \{\lambda \in {\C}: \lambda = re^{i\theta} \; \text{for} \;
r \geq 0, \; \theta \in {\R}, \; |\theta - \theta_0| \leq a\},
\end{equation*}
with $\theta_0 \in {\R}$ and $a > 0$,
is called a sector of minimal growth for the extension
$$
A_{\Dom} : \Dom \subset x^{-m/2}L_b^2(M;E) \to x^{-m/2}L^2_b(M;E)
$$
if for $\lambda \in \Lambda$ with $|\lambda| > R$ sufficiently large
$$
A_\Dom - \lambda : \Dom \to x^{-m/2}L_b^2(M;E)
$$
is invertible, and the resolvent $(A_\Dom-\lambda)^{-1}$ satisfies the 
equivalent norm estimates  
\begin{equation}\label{ResolventEstimates}
\begin{split} 
 &\bigl\|(A_\Dom-\lambda)^{-1}\bigr\|_{\L(x^{-m/2}L^2_b)}
 = O(|\lambda|^{-1}) \;\text{ as } |\lambda| \to \infty,\\
 &\bigl\|(A_\Dom-\lambda)^{-1}\bigr\|_{\L(x^{-m/2}L^2_b,
 \Dom_{\max})} = O(1)
 \;\text{ as } |\lambda| \to \infty.
 \end{split}
\end{equation}

Analogously, we call $\Lambda$ a sector of minimal growth for 
$A_{\wedge,\Dom_{\wedge}}$ if 
$$
A_{\wedge,\Dom_{\wedge}} - \lambda : \Dom_{\wedge} \to 
x^{-m/2}L_b^2(Y^{\wedge};E)
$$
is invertible for large $|\lambda| > 0$ in $\Lambda$, and the inverse
satisfies the equivalent estimates
\begin{equation}\label{AWedgeRMG}
\begin{split}
 &\bigl\|(A_{\wedge,\Dom_\wedge}-\lambda)^{-1}\bigr\|_{\L(x^{-m/2}L^2_b)}
 = O(|\lambda|^{-1}) \;\text{ as } |\lambda| \to \infty,\\
 &\bigl\|(A_{\wedge,\Dom_\wedge}-\lambda)^{-1}\bigr\|_{\L(x^{-m/2}L^2_b,
 \Dom_{\max})} = O(1) \;\text{ as } |\lambda| \to \infty.
\end{split} 
\end{equation}

\begin{theorem}\label{RayMinimalGrowth} 
Let $A\in x^{-m}\Diff_b^m(M;E)$ be $c$-elliptic with parameter in $\Lambda$. 
Let $\Dom \subset x^{-m/2}L^2_b(M;E)$ be a domain such that $A_\Dom$ is closed
and let $\Dom_{\wedge}$ be the associated domain defined via 
\eqref{DwedgeIdentitaet}. Assume that $\Lambda$ is a sector of minimal 
growth for the extension $A_{\wedge,\Dom_\wedge}$.
Then $\Lambda$ is a sector of minimal growth for the operator $A_{\Dom}$.
Moreover, the resolvent of $A_\Dom$ satisfies the equation
\begin{equation}\label{ResolventeDarstellung}
 (A_\Dom-\lambda)^{-1} = B(\lambda) + (A_\Dom-\lambda)^{-1}\Pi(\lambda)
\end{equation}
with the parametrix $B(\lambda)$ and the projection $\Pi(\lambda)$
from Theorem~\ref{PDrei}.
\end{theorem}

Before we prove this theorem, we discuss some interesting properties
of the resolvent conditions on $A_\wedge$. For more details see 
\cite{GilKrainerMendoza1}.

\begin{proposition} 
If $\Dom_\wedge$ is $\kappa$-invariant, then the invertibility of
$A_{\wedge,\Dom_\wedge}-\lambda$ for $\lambda \in \Lambda$ with
$|\lambda| > R$ implies the invertibility of $A_{\wedge,\Dom_\wedge}
-\lambda$ for all $\lambda\in\Lambda\minus \{0\}$, and $\Lambda$ is
a sector of minimal growth for $A_{\wedge,\Dom_{\wedge}}$.
\end{proposition}
\begin{proposition} 
If $\Lambda$ is a sector of minimal growth for the operator 
$A_{\wedge}$ with domain $\Dom_{\wedge}$, then $\Lambda$ is also a sector of 
minimal growth for $A_{\wedge}$
with domain $\kappa_{\varrho}\Dom_{\wedge}$ for any $\varrho > 0$.
In particular, the resolvent $B_{\varrho,\wedge}(\lambda)$ of 
$A_{\wedge,\kappa_{\varrho}\Dom_{\wedge}}$ satisfies
\begin{equation*}
 B_{\varrho,\wedge}(\lambda) = \varrho^{-m}\kappa_{\varrho}
(A_{\wedge,\Dom_\wedge} - \varrho^{-m}\lambda)^{-1}\kappa^{-1}_{\varrho}.
\end{equation*}
\end{proposition}

In general, the norm estimates \eqref{AWedgeRMG} are not easy to
check. However, the following proposition shows that this resolvent condition
only needs to be verified for the projection of 
$(A_{\wedge,\Dom_\wedge}-\lambda)^{-1}$ onto the finite dimensional space 
$\tilde\Sing_{\wedge,\max}=\Dom_{\wedge,\max}/\Dom_{\wedge,\min}$.

\begin{proposition}\label{PropSmaxCond}
Let $A$ be c-elliptic with parameter in $\Lambda$. The sector $\Lambda$ 
is a sector of minimal growth for $A_{\wedge,\Dom_\wedge}$ if and only if
$$
A_{\wedge,\Dom_\wedge}-\lambda: \Dom_{\wedge}\to x^{-m/2}L_b^2(Y^{\wedge};E)
$$
is invertible for large $|\lambda|>0$, and the inverse satisfies the estimate
\begin{equation}\label{SmaxCondition}
 \bigl\|\kappa_{|\lambda|^{1/m}}^{-1}q_{\wedge}
 (A_{\wedge,\Dom_\wedge}-\lambda)^{-1}\bigr\|_{\L(x^{-m/2}L^2_b,
 \tilde\Sing_{\wedge,\max})} = O(|\lambda|^{-1}) 
\;\text{ as } |\lambda|\to\infty.
\end{equation}
Here $q_{\wedge}: \Dom_{\wedge,\max} \to \tilde\Sing_{\wedge,\max}$ denotes 
the canonical projection.
\end{proposition}
\begin{proof}
We first observe that the $\kappa$-homogeneity of $A_{\wedge}$ implies
$$
A_{\wedge}\kappa_{|\lambda|^{1/m}}^{-1}(A_{\wedge,\Dom_\wedge} - \lambda)^{-1} 
=\kappa_{|\lambda|^{1/m}}^{-1}|\lambda|^{-1}A_{\wedge}(A_{\wedge,\Dom_\wedge} 
- \lambda)^{-1}
$$
as operators in $\L(x^{-m/2}L^2_b)$. Using this identity and the fact that
$\kappa_{\varrho}$ is an isometry in $\L(x^{-m/2}L^2_b)$, one can easily
see that the estimates \eqref{AWedgeRMG} are equivalent to
\begin{equation}\label{KappaResolventEstimate}
\|\kappa_{|\lambda|^{1/m}}^{-1}
(A_{\wedge,\Dom_\wedge}-\lambda)^{-1}\|_{\L(x^{-m/2}L^2_b,\Dom_{\wedge,\max})} 
= O(|\lambda|^{-1}) \;\text{ as } |\lambda|\to\infty,
\end{equation}
and therefore \eqref{SmaxCondition} holds. Note that 
$\kappa_\varrho q_\wedge = q_\wedge\kappa_\varrho$.

Conversely, assume that we have \eqref{SmaxCondition}. 
Let $B_{\wedge}(\lambda)$ be the principal part of
the parametrix $B(\lambda)$ from Theorem~\ref{PDrei}. Then, for $\lambda
\in\Lambda\minus\{0\}$, we have
\begin{equation*}
1-B_{\wedge}(\lambda)(A_{\wedge}-\lambda)=0 \;\text{ on } \Dom_{\wedge,\min},
\end{equation*}
and we may write
\begin{equation*}
(A_{\wedge,\Dom_\wedge} - \lambda)^{-1} = B_{\wedge}(\lambda) +
(1 - B_{\wedge}(\lambda)(A_{\wedge} - \lambda))
q_\wedge(A_{\wedge,\Dom_\wedge} - \lambda)^{-1}
\end{equation*}
as operators in $\L(x^{-m/2}L_b^2,\Dom_{\wedge,\max})$. 
Since $B_{\wedge}(\lambda)$ and $(A_\wedge-\lambda)$ are $\kappa$-homogeneous 
of degree $-m$ and $m$, respectively, we have the identities
\begin{align*}
\kappa_{|\lambda|^{1/m}}^{-1} B_\wedge(\lambda) &=
|\lambda|^{-1}B_\wedge\bigl(\tfrac{\lambda}{|\lambda|}\bigr)
\kappa_{|\lambda|^{1/m}}^{-1},\\
\kappa_{|\lambda|^{1/m}}^{-1} (A_\wedge -\lambda) &=
|\lambda| \bigl(A_\wedge-\tfrac{\lambda}{|\lambda|}\bigr)
\kappa_{|\lambda|^{1/m}}^{-1},
\end{align*}	
which imply
\begin{multline*}
\kappa_{|\lambda|^{1/m}}^{-1} (A_{\wedge,\Dom_\wedge} -\lambda)^{-1} 
= |\lambda|^{-1}B_\wedge\bigl(\tfrac{\lambda}{|\lambda|}\bigr)
\kappa_{|\lambda|^{1/m}}^{-1} \\ 
+ \Bigl(1 - B_{\wedge}\bigl(\tfrac{\lambda}{|\lambda|}\bigr)
\bigl(A_{\wedge} - \tfrac{\lambda}{|\lambda|}\bigr)\Bigr)
\kappa_{|\lambda|^{1/m}}^{-1}q_\wedge(A_{\wedge,\Dom_\wedge} - \lambda)^{-1}.
\end{multline*}
Passing to the norm in $\L(x^{-m/2}L_b^2,\Dom_{\wedge,\max})$ and using 
\eqref{SmaxCondition} we obtain \eqref{KappaResolventEstimate} which is 
equivalent to the estimates \eqref{AWedgeRMG}. 
\end{proof}

For the proof of Theorem \ref{RayMinimalGrowth} we need further
ingredients.  First of all, using the operator $\theta$ defined via
\eqref{thetainvers} and \eqref{ThetaOperator}, we now define on 
$\tilde{\Sing}_{\max}$ the group action
\begin{equation}\label{kappaSchlange}
\tilde{\kappa}_{\varrho}= \theta^{-1}\kappa_{\varrho}\theta \quad \text{for $\varrho > 0$.}
\end{equation}
We may write $\tilde{\kappa}_{\varrho} = \kappa_{\varrho}L_{\varrho}$, where
$$
L_{\varrho} = \kappa_{\varrho}^{-1}\theta^{-1}\kappa_{\varrho}\theta :
\tilde{\Sing}_{\max} \to C^{\infty}(\open Y^{\wedge};E)
$$
is the direct sum of the operators $L_{\varrho}|_{\tilde{\Sing}_{\sigma_0}}$ which act
as follows:

For $\tilde{u} \in \tilde{\Sing}_{\sigma_0}$ we have
\begin{equation}\label{AltLRho}
L_{\varrho}\tilde{u} = \sum\limits_{\vartheta=0}^{N(\sigma_0)}\varrho^{-\vartheta}\e_{\sigma_0,\vartheta}(\varrho)(\theta\tilde{u}),
\end{equation}
where $\e_{\sigma_0,\vartheta}(\varrho)$ is defined as
\begin{equation*}
\e_{\sigma_0,\vartheta}(\varrho) = \varrho^{\vartheta}\kappa_{\varrho}^{-1}
\e_{\sigma_0,\vartheta}\kappa_{\varrho}:
\tilde{\Sing}_{\wedge,\sigma_0}\to C^{\infty}(\open Y^{\wedge};E).
\end{equation*}
In particular, $\e_{\sigma_0,0}(\varrho)(\tilde u) = \tilde u$ for
all $\varrho\in\R_+$ and $\tilde u\in\tilde\Sing_{\wedge,\sigma_0}$.

\begin{lemma}\label{Ltheta} \
\begin{enumerate}[$(i)$]
\item For every $\psi \in \tilde{\Sing}_{\wedge,\sigma_0}$ and every 
$\vartheta \in \{0,\ldots,N(\sigma_0)\}$ there exists a polynomial 
$q_\vartheta(y,\log x,\log \varrho)$ in $(\log x,\log \varrho)$ with
coefficients in $C^{\infty}(Y;E)$ such that
\begin{equation}\label{esigmavarrho}
\e_{\sigma_0,\vartheta}(\varrho)(\psi) = q_\vartheta(y,\log x,\log \varrho)
x^{i(\sigma_0 - i\vartheta)},
\end{equation}
and the degree of $q_\vartheta$ with respect to $(\log x,\log \varrho)$ is 
bounded by some $\mu \in \N_0$ which is independent of $\sigma_0\in\Sigma$, 
$\psi \in \tilde{\Sing}_{\wedge,\sigma_0}$, and
$\vartheta \in \{0,\ldots,N(\sigma_0)\}$.

\item Let $\omega \in C_0^{\infty}(\overline{\R}_+)$ be any cut-off
function near the origin, i.e.,  $\omega = 1$ near zero and $\omega =
0$ near infinity. Then the operator family
$$
\omega\bigl(L_{\varrho} - \theta\bigr) : \tilde{\Sing}_{\max} \to \K^{\infty,-m/2}(Y^{\wedge};E)
$$
satisfies for every $s \in \R$ the norm estimate
$$
\|\omega\bigl(L_{\varrho} - \theta\bigr)\|_{\L(\tilde{\Sing}_{\max},
\K^{s,-m/2})} = O(\varrho^{-1}\log^{\mu}\varrho)
\;\text{ as } \varrho \to \infty, 
$$
where $\mu \in \N_0$ is the bound for the degrees of the polynomials
$q_\vartheta$ in $(i)$, and $\K^{s,-m/2}(Y^{\wedge};E)$ is the weighted
Sobolev space defined in Section \ref{sec-Preliminaries}.
\end{enumerate}
\end{lemma}
\begin{proof}
As $\Sigma$ is a finite set and all spaces
$\tilde{\Sing}_{\wedge,\sigma_0}$ are finite dimensional, it
suffices to show that \eqref{esigmavarrho} holds for a
basis of $\tilde{\Sing}_{\wedge,\sigma_0}$.  
We pick a basis $\{\psi_0,\ldots,\psi_K\}
\subset \tilde{\Sing}_{\wedge,\sigma_0}$ which is a Jordan basis for 
the infinitesimal generator $(\frac{m}2+x\partial_x)$ of the group
$\kappa_{\varrho}|_{\tilde{\Sing}_{\wedge,\sigma_0}} \in
\L(\tilde{\Sing}_{\wedge,\sigma_0})$.
Recall that $\tilde{\Sing}_{\wedge,\max}$ is $\kappa$-invariant, and so 
are necessarily all the spaces $\tilde{\Sing}_{\wedge,\sigma_0}$.
Note that the only eigenvalue of $(\frac{m}2+x\partial_x)$ on
$\tilde{\Sing}_{\wedge,\sigma_0}$ is ${m/2+i\sigma_0}$.

Consequently, for each $j$ we may write
$$
\kappa_{\varrho}\psi_j = \varrho^{m/2+i\sigma_0}
 \sum_{k=0}^{K}p_{jk}(\log\varrho) \psi_k,
$$
where $p_{jk}$ is a polynomial, and thus
$$
\e_{\sigma_0,\vartheta}(\varrho)(\psi_j) = 
\varrho^{\vartheta}\kappa_{\varrho}^{-1}\e_{\sigma_0,\vartheta}
(\kappa_{\varrho}\psi_j) = \sum_{k=0}^{K}p_{jk}(\log\varrho) 
\varrho^{i(\sigma_0-i\vartheta)}\varrho^{m/2}\kappa_{\varrho}^{-1}
\e_{\sigma_0,\vartheta}(\psi_k).
$$
Every $\e_{\sigma_0,\vartheta}(\psi_k)$ is a singular function of the form
$$
\Bigl(\sum\limits_{\nu=0}^{m^{(k)}_{\sigma_0-i\vartheta}}
c^{(k)}_{\sigma_0-i\vartheta,\nu}(y)\log^{\nu}x\Bigr)
x^{i(\sigma_0-i\vartheta)},
$$
and so
$$
\varrho^{i(\sigma_0-i\vartheta)}\varrho^{m/2}\kappa_{\varrho}^{-1}
\e_{\sigma_0,\vartheta}(\psi_k) =
\Bigl(\sum\limits_{\nu=0}^{m^{(k)}_{\sigma_0-i\vartheta}}
c^{(k)}_{\sigma_0-i\vartheta,\nu}(y)(\log x - \log \varrho)^{\nu}\Bigr)
x^{i(\sigma_0-i\vartheta)}.
$$
Hence (i) is proved.

For the proof of (ii) note that according to \eqref{AltLRho} and
(i), we have for $\tilde{u} \in \tilde{\Sing}_{\sigma_0}$ 
\begin{align*}
\omega\bigl(L_{\varrho} - \theta\bigr)\tilde{u} &= \omega\sum\limits_{\vartheta=1}^{N(\sigma_0)}
\varrho^{-\vartheta}\e_{\sigma_0,\vartheta}(\varrho)(\theta\tilde{u}) \\
&= \varrho^{-1}\sum\limits_{\vartheta=1}^{N(\sigma_0)}\varrho^{1-\vartheta}\omega q_{\vartheta}(y,\log x,\log \varrho)x^{i(\sigma_0-i\vartheta)},
\end{align*}
and consequently
$$
\|\omega\bigl(L_{\varrho} - \theta\bigr)\tilde{u}\|_{\K^{s,-m/2}} 
\leq \textup{const}\cdot \bigl(\varrho^{-1}\log^{\mu}\varrho\bigr)
$$
for $\varrho \geq 1$, which then in fact holds for all $\tilde{u} \in \tilde{\Sing}_{\max}$. As
$$
\omega\bigl(L_{\varrho} - \theta\bigr) : \tilde{\Sing}_{\max} \to \K^{s,-m/2}(Y^{\wedge};E)
$$
is continuous for every $\varrho > 0$, we obtain (ii) from the Banach-Steinhaus theorem.
\end{proof}

\begin{lemma}\label{kappaSchlangelift}
Fix a cut-off function $\omega \in C_0^{\infty}([0,1))$ near $0$. 
For $\varrho \geq 1$ consider the operator family
\begin{equation*}
\tilde{K}(\varrho)= \omega_{\varrho} \tilde{\kappa}_{\varrho} :
\tilde\Sing_{\max} \to \Dom^{\infty}_{\max}(A) 
= \bigcap_{t \in {\R}} \Dom^t_{\max}(A),
\end{equation*}
where $\omega_{\varrho}(x)= \omega(\varrho x)$. If $q: \Dom_{\max}(A) \to 
\tilde\Sing_{\max}$ is the canonical projection, then 
\begin{equation*}
q{\circ}\tilde{K}(\varrho) = \tilde{\kappa}_{\varrho},
\end{equation*}
and we have the norm estimates 
\begin{align} \label{Kestimate1}
&\|\tilde{K}(\varrho)\|_{\L(\tilde\Sing_{\max},x^{-m/2}L^2_b)} = O(1)
\;\text{ as } \varrho\to\infty, \\ \label{Kestimate2}
&\|\tilde{K}(\varrho)\|_{\L(\tilde\Sing_{\max},\Dom_{\max})} = O(\varrho^m)
\;\text{ as } \varrho\to\infty.
\end{align}
Moreover, for every $t \in \R$ there exists $M_t \in\R$ such that
\begin{equation}\label{Kestimate3}
\|\tilde{K}(\varrho)\|_{\L(\tilde\Sing_{\max},\Dom^t_{\max})} = O(\varrho^{M_t})
\;\text{ as } \varrho\to\infty.
\end{equation}
\end{lemma}
\begin{proof}
That $\tilde{K}(\varrho)$ is a lift of $\tilde{\kappa}_{\varrho}$ to 
$\Dom_{\max}^{\infty}(A)$ is evident from the definition. 
In order to show the norm estimates, it is sufficient 
to consider for each $\sigma_0 \in \Sigma$ the restriction
$$
\tilde{K}_{\sigma_0}(\varrho)= \tilde{K}(\varrho)|_{\tilde\Sing_{\sigma_0}}: 
\tilde\Sing_{\sigma_0} \to \Dom_{\max}^{\infty}(A)
$$
and prove the estimates for this operator.
Recall that $\tilde\kappa_\varrho=\kappa_\varrho L_\varrho$ so that for 
$\tilde u\in \tilde\Sing_{\sigma_0}$ we have $\tilde K_{\sigma_0}(\varrho)
\tilde u = \kappa_\varrho(\omega L_\varrho\tilde u)$. On the other hand, 
by Lemma~\ref{Ltheta}, $\omega L_{\varrho} \to \omega \theta$ in 
$\L(\tilde\Sing_{\max},x^{-m/2}L_b^2)$ as $\varrho \to \infty$, so the 
family $\omega L_{\varrho}$ is uniformly bounded for $\varrho \ge 1$. Thus
\begin{equation*}
\|\tilde{K}_{\sigma_0}(\varrho)\tilde u\|_{x^{-m/2}L_b^2(M;E)} 
\leq \textup{const} \|\kappa_{\varrho}\bigl(\omega 
L_{\varrho}\tilde u\bigr)\|_{x^{-m/2}L_b^2(Y^{\wedge};E)} 
\leq \textup{const} \|\omega \tilde u\|_{\Dom_{\max}}
\end{equation*}
since the norm $\|\omega \tilde u\|_{\Dom_{\max}}$ is an admissible norm 
on the finite dimensional space $\tilde\Sing_{\sigma_0}$.
Recall that $\kappa_{\varrho}$ is an isometry in $x^{-m/2}L_b^2$. 
Finally, the above estimate gives \eqref{Kestimate1}.

For proving \eqref{Kestimate2} we only need to show that
$$
\|A\tilde{K}_{\sigma_0}(\varrho)\|_{\L(\tilde\Sing_{\sigma_0},x^{-m/2}L_b^2)} =
O(\varrho^m) \;\text{ as } \varrho \to \infty. 
$$
Thus we will prove that there exists a constant $C > 0$, independent of 
$\tilde u\in \tilde\Sing_{\sigma_0}$ and $\varrho \geq 1$, such that
\begin{equation*}
\|A(\kappa_{\varrho}(\omega L_{\varrho}\tilde u))\|_{x^{-m/2}
L_b^2(Y^{\wedge};E)} \leq C\varrho^m \|\omega \tilde u\|_{\Dom_{\max}}.
\end{equation*}
To this end we split $A$ near the boundary as in \eqref{APhiDarstellung} 
and use \eqref{AltLRho} to obtain
\begin{equation}\label{PhiAPhiL}
\begin{split}
A(\kappa_{\varrho} & (\omega L_{\varrho}\tilde u)) \\
&= \Bigl(x^{-m}\sum_{k=0}^{m-1}P_kx^k\Bigr)\kappa_{\varrho}
\bigl({\omega}L_{\varrho}\tilde u\bigr) +
\tilde P_m\kappa_{\varrho}\bigl({\omega}L_{\varrho}\tilde u\bigr) \\
&= \varrho^m\kappa_{\varrho}
\Bigl(x^{-m}\sum_{k=0}^{m-1}\varrho^{-k}P_kx^k\Bigr)
\Bigl({\omega}\sum_{j=0}^{N(\sigma_0)}\varrho^{-j}
\e_{\sigma_0,j}(\varrho)(\theta \tilde u)\Bigr) +
\tilde P_m\kappa_{\varrho}\bigl({\omega}L_{\varrho}\tilde u\bigr) \\
&= \sum_{\vartheta=0}^{2m-2} \varrho^{m-\vartheta}\kappa_{\varrho}\Bigl(
x^{-m}\!\!\!\sum_{\substack{k+j=\vartheta \\ 0 \leq k,j \leq m-1}}\!\!\!
\bigl(P_kx^k\bigr)\bigl({\omega}\e_{\sigma_0,j}(\varrho)(\theta \tilde u)\bigr)
\Bigr) +\tilde P_m\kappa_{\varrho}\bigl({\omega}L_{\varrho}\tilde u\bigr)
\end{split}
\end{equation}
with the convention that $e_{\sigma_0,j}(\varrho)= 0$ for $j > N(\sigma_0)$.

For every $\vartheta\in\{0,\dots,2m-2\}$ we consider the family of linear maps 
\begin{equation}\label{KreisSumme}
\tilde u\mapsto 
x^{-m}\!\!\sum_{\substack{k+j=\vartheta \\ 0 \leq k,j \leq m-1}}\!\!
\bigl(P_kx^k\bigr)\bigl({\omega}\e_{\sigma_0,j}(\varrho)(\theta \tilde u)\bigr)
:\tilde\Sing_{\sigma_0}\to x^{-m/2}L_b^2(Y^{\wedge};E).
\end{equation}
We will prove that \eqref{KreisSumme} is well-defined, i.e., every $\tilde u \in \tilde{\Sing}_{\sigma_0}$ is 
indeed mapped into $x^{-m/2}L_b^2(Y^{\wedge};E)$, and that the norms are 
bounded by a constant times $\log^\mu \varrho$ as $\varrho\to\infty$ with
$\mu$ as in Lemma~\ref{Ltheta}. 
Thus for every $\vartheta \in \{0,\ldots,2m-2\}$ we have
\begin{multline*}
\Bigl\|\varrho^{m-\vartheta}\kappa_{\varrho}\Bigl(
x^{-m}\!\!\!\sum_{\substack{k+j=\vartheta \\ 0 \leq k,j \leq m-1}}
\!\!\bigl(P_kx^k\bigr)\bigl({\omega}\e_{\sigma_0,j}(\varrho)(\theta \tilde u)
\bigr)\Bigr)\Bigr\|_{x^{-m/2}L^2_b} \\
\leq \textup{const}\cdot\bigl(\varrho^{m-\vartheta} \log^\mu
\varrho\bigr) \|\omega \tilde u\|_{\Dom_{\max}},
\end{multline*}
while for $\vartheta=0$, 
\begin{equation}\label{Theta0Estimate}
\varrho^m \kappa_\varrho x^{-m} P_0\omega 
\e_{\sigma_0,0}(\varrho)(\theta \tilde u) =
\varrho^m \kappa_\varrho A_\wedge\omega (\theta \tilde u) =
A_\wedge \kappa_\varrho\bigl(\omega \theta \tilde u\bigr),
\end{equation}
so for this term we have a norm estimate without $\log$.

Let $\tilde\omega\in C_0^\infty(\overline\R_+)$ be a cut-off function near 
$0$ with $\omega\prec\tilde\omega$. Then there exist suitable 
$\vp,\tilde \vp\in C_0^\infty(\R_+)$ such that for all 
$\tilde u\in\tilde\Sing_{\sigma_0}$,
\begin{multline}\label{SumSplitting}
x^{-m}\!\!\!\sum_{\substack{k+j=\vartheta \\ 0 \leq k,j \leq m-1}}\!\!
\big(P_kx^k\big)\big({\omega}\e_{\sigma_0,j}(\varrho)(\theta \tilde u)\big)\\
= \tilde\omega x^{-m}\!\!\!\sum_{\substack{k+j=\vartheta \\ 0 \leq k,j 
\leq m-1}}\!\! \bigl(P_kx^k\bigr)\e_{\sigma_0,j}(\varrho)(\theta \tilde u) +
\tilde\vp x^{-m}\!\!\!\sum_{\substack{k+j=\vartheta \\ 0 \leq k,j \leq m-1}}
\!\!  \bigl(P_kx^k\bigr)\vp\e_{\sigma_0,j}(\varrho)(\theta \tilde u).
\end{multline}

According to Lemma \ref{Ltheta} the second sum in \eqref{SumSplitting} is a
polynomial in $\log \varrho$ of degree at most $\mu$ with coefficients in
$x^{-m/2}L^2_b(Y^{\wedge};E)$. As both $A(\kappa_{\varrho}(\omega L_{\varrho}\tilde u))$
and $\tilde P_m(\kappa_{\varrho}({\omega}L_{\varrho}\tilde u))$ belong to
$x^{-m/2}L_b^2(Y^{\wedge};E)$, we get from the equations \eqref{PhiAPhiL} and
\eqref{SumSplitting} that necessarily
$$
x^{-m}\!\!\!\sum_{\substack{k+j=\vartheta \\ 0 \leq k,j \leq m-1}}\!\!
\big(P_kx^k\big)\big({\omega}\e_{\sigma_0,j}(\varrho)(\theta \tilde u)\big)
\in x^{-m/2}L_b^2(Y^{\wedge};E)
$$
for all $\varrho \in \R_+$ and all $\tilde{u} \in \tilde{\Sing}_{\sigma_0}$, and, moreover, that
$$
\tilde{\omega}x^{-m}\!\!\!\sum_{\substack{k+j=\vartheta \\ 0 \leq k,j \leq m-1}}\!\!
\big(P_kx^k\big)\e_{\sigma_0,j}(\varrho)(\theta \tilde u) = 0
$$
for $\vartheta \le N(\sigma_0)$ because these functions are of the form
$$
\tilde{\omega}\Bigl(\sum\limits_{\nu}c_{\sigma_0-i(\vartheta-m),\nu}(y)\log^{\nu} x\Bigr)x^{i(\sigma_0-i(\vartheta-m))}.
$$
For $\vartheta > N(\sigma_0)$ every single summand
$\tilde{\omega}x^{-m}\big(P_kx^k\big)\e_{\sigma_0,j}(\varrho)(\theta \tilde u)$ belongs to the space
$x^{-m/2}L_b^2(Y^{\wedge};E)$, and by Lemma \ref{Ltheta} is a polynomial in $\log \varrho$ of degree
at most $\mu$ with coefficients in $x^{-m/2}L_b^2(Y^{\wedge};E)$.

Summing up, we have shown that for every $\tilde{u} \in \tilde{\Sing}_{\sigma_0}$ the function
$$
x^{-m}\!\!\!\sum_{\substack{k+j=\vartheta \\ 0 \leq k,j \leq m-1}}\!\!
\big(P_kx^k\big)\big({\omega}\e_{\sigma_0,j}(\varrho)(\theta \tilde u)\big)
$$
is a polynomial in $\log\varrho$ of degree at most $\mu$ with coefficients in 
$x^{-m/2}L_b^2(Y^{\wedge};E)$, and from the Banach-Steinhaus theorem we now obtain the desired norm
estimates for the family of maps \eqref{KreisSumme}.

On the other hand, 
\begin{align*}
\|\tilde P_m\kappa_{\varrho}\bigl({\omega}L_{\varrho}
\tilde u\bigr)\|_{x^{-m/2}L_b^2} &= 
\|\kappa_{\varrho}^{-1}\tilde P_m\kappa_{\varrho}\bigl({\omega}L_{\varrho}
\tilde u\bigr)\|_{x^{-m/2}L_b^2}\\
&\le \|\kappa_{\varrho}^{-1}\tilde P_m\kappa_{\varrho}\|_{\L(\K^{m,-m/2},
x^{-m/2}L^2_b)}\|{\omega}L_{\varrho}\tilde u\|_{\K^{m,-m/2}}.
\end{align*}
Lemma~\ref{Ltheta} implies 
$\|{\omega}L_{\varrho}\tilde u\|_{\K^{m,-m/2}}\le \textup{const} 
\|\omega \tilde u\|_{\Dom_{\max}}$, and so 
\begin{equation*}
\|\tilde P_m\kappa_{\varrho}\bigl({\omega}L_{\varrho}\tilde u\bigr)\|_{
x^{-m/2}L_b^2} \leq \textup{const} \|\omega \tilde u\|_{\Dom_{\max}}
\end{equation*}
since $\|\kappa_{\varrho}^{-1}\tilde P_m\kappa_{\varrho}\|_{\L(\K^{m,-m/2},
x^{-m/2}L^2_b)}=O(1)$ as $\varrho\to\infty$. Thus \eqref{Kestimate2} is proved. 

Finally, an inspection of the proof reveals that for $t \in \R$ we obtain
\begin{align*}
&\|\tilde{K}(\varrho)\|_{\L(\tilde\Sing_{\max},x^{-m/2}H^t_b)} = 
O(\|\kappa_{\varrho}\|_{\L({\mathcal K}^{t,-m/2})})
\;\text { as } \varrho \to \infty, \\
&\|\tilde{K}(\varrho)\|_{\L(\tilde\Sing_{\max},\Dom^t_{\max})} =
O(\varrho^m\|\kappa_{\varrho}\|_{\L({\mathcal K}^{t,-m/2})})
\;\text { as } \varrho \to \infty,
\end{align*}
and consequently \eqref{Kestimate3} follows because the norm 
$\|\kappa_{\varrho}\|_{\L({\mathcal K}^{t,-m/2})}$ behaves polynomially
as $\varrho \to \infty$.
\end{proof}

\bigskip 
\begin{proof}[Proof of Theorem~\ref{RayMinimalGrowth}]
Fix some complement $\Sing_{\max}$ of $\Dom_{\min}$ in $\Dom_{\max}$
and let $\Sing \subset \Sing_{\max}$ be a subspace such that
$\Dom = \Dom_{\min}\oplus \Sing$. With respect to this decomposition the 
operator $A_\Dom-\lambda$ can be written as
\begin{equation*}
(A_\Dom - \lambda) = 
\begin{pmatrix} 
(A -\lambda)|_{\Dom_{\min}} & (A -\lambda)|_{\Sing}
\end{pmatrix} :
\begin{array}{c} \Dom_{\min} \\ \oplus \\ \Sing \end{array} 
\to x^{-m/2}L_b^2(M;E).
\end{equation*}   
Let $d''=\dim \Sing$. Under the ellipticity condition on $A-\lambda$  
and the injectivity of $A_\wedge-\lambda$ on $\Dom_{\wedge,\min}$  
we already proved in Theorem~\ref{AminuslambdaplusK} the existence of  
a parametrix $B(\lambda)$ of $A-\lambda$ on $\Dom_{\min}$ and a 
generalized Green remainder 
$\begin{pmatrix} 0 & K(\lambda)\end{pmatrix}$ of order $m$ such that
\begin{equation*}
\begin{pmatrix} 
(A -\lambda)|_{\Dom_{\min}} & K(\lambda)
\end{pmatrix} :
\begin{array}{c} \Dom_{\min} \\ \oplus \\ \C^{d''} \end{array} 
\to x^{-m/2}L_b^2(M;E)
\end{equation*}   
is invertible for $\lambda$ sufficiently large with inverse
\begin{equation}\label{LeftDParametrix}
\begin{pmatrix} 
(A -\lambda)|_{\Dom_{\min}} & K(\lambda)
\end{pmatrix}^{-1} = 
\begin{pmatrix} 
B(\lambda) \\  T(\lambda)
\end{pmatrix}, 
\end{equation}
where $\begin{pmatrix}0 \\ T(\lambda)\end{pmatrix}$ 
is a generalized Green remainder of order $-m$. Since
\begin{align*}
\id =
\begin{pmatrix} 
B(\lambda) \\  T(\lambda)
\end{pmatrix}
\begin{pmatrix} 
(A -\lambda)|_{\Dom_{\min}} & K(\lambda)
\end{pmatrix} 
= 
\begin{pmatrix} 
B(\lambda)(A-\lambda)|_{\Dom_{\min}} & B(\lambda)K(\lambda) \\  
T(\lambda)(A-\lambda)|_{\Dom_{\min}} & T(\lambda)K(\lambda)
\end{pmatrix}, 
\end{align*}    
we have $B(\lambda)(A-\lambda)|_{\Dom_{\min}} = 1$, 
$T(\lambda)(A-\lambda)|_{\Dom_{\min}}=0$, and $T(\lambda)K(\lambda)=1$. Then
\begin{equation}\label{CompWithDParametrix}
\begin{pmatrix} 
B(\lambda) \\  T(\lambda)
\end{pmatrix}
\begin{pmatrix} 
(A -\lambda)|_{\Dom_{\min}} & (A-\lambda)|_{\Sing}
\end{pmatrix} = 
\begin{pmatrix} 
1 & B(\lambda)(A-\lambda)|_{\Sing} \\  
0 & T(\lambda)(A-\lambda)|_{\Sing}
\end{pmatrix} 
\end{equation}
which implies that $\begin{pmatrix}(A -\lambda)|_{\Dom_{\min}} 
& (A-\lambda)|_{\Sing}\end{pmatrix}$ is invertible if and only if 
\begin{equation}\label{FLambda}
F(\lambda) = T(\lambda)(A-\lambda): \Sing\to \C^{d''}
\end{equation}    
is invertible. Moreover, we get the explicit representation
\begin{equation}\label{TheResolvent}
(A_\Dom - \lambda)^{-1} =
B(\lambda) + \bigl(1- B(\lambda)(A - \lambda)\bigr)
F(\lambda)^{-1}T(\lambda),
\end{equation}
and \eqref{ResolventeDarstellung} follows from Corollary 
\ref{InverseDarstellung}.

As $F(\lambda)$ and $1 - B(\lambda)(A-\lambda)$ vanish
on $\Dom_{\min}$ for large $\lambda$, they descend to operators 
$F(\lambda) : \tilde\Sing_{\max} \to \C^{d''}$ and 
$1 - B(\lambda)(A-\lambda) : \tilde\Sing_{\max} \to \Dom_{\max}$. 
If $\tilde\Sing=\Dom/\Dom_{\min}$, then the invertibility of \eqref{FLambda} 
is equivalent to the invertibility of
\begin{equation*}
F(\lambda) : \tilde\Sing \to \C^{d''},
\end{equation*}
and in this case, \eqref{TheResolvent} still makes sense in this context.

Let $q : \Dom_{\max} \to \tilde \Sing_{\max}$ be the canonical projection.
The resolvent $(A_{\Dom} - \lambda)^{-1}$ and $F(\lambda)^{-1} : \C^{d''} \to
\tilde\Sing_{\max}$ are related by the formulas
\begin{align*}
F(\lambda)^{-1} = q(A_{\Dom}-\lambda)^{-1}K(\lambda) &: 
\C^{d''} \to \tilde\Sing_{\max}, \\
q(A_{\Dom}-\lambda)^{-1} = F(\lambda)^{-1}T(\lambda) &: 
x^{-m/2}L_b^2 \to \tilde\Sing_{\max}
\end{align*}
in view of $T(\lambda)K(\lambda) = 1$.

Under the assumptions of Theorem \ref{RayMinimalGrowth} we will prove that
$F(\lambda) : \tilde\Sing \to \C^{d''}$ is invertible for large $\lambda$, 
and that the inverse satisfies the norm estimate 
\begin{equation}\label{Festimate}
\|\tilde{\kappa}_{[\lambda]^{1/m}}^{-1}
F(\lambda)^{-1}\|_{\L(\C^{d''},\tilde\Sing_{\max})} = O(1)
\;\text{ as } |\lambda| \to \infty.
\end{equation}
Observe that the parametrix construction from Theorem~\ref{AminuslambdaplusK}
gives the relation
\begin{equation*}
\begin{pmatrix} 
(A_\wedge -\lambda)|_{\Dom_{\wedge,\min}} & K_\wedge(\lambda)
\end{pmatrix}^{-1} = 
\begin{pmatrix} 
B_\wedge(\lambda) \\  T_\wedge(\lambda)
\end{pmatrix}
\end{equation*}
for the $\kappa$-homogeneous principal parts of \eqref{LeftDParametrix}. 
Thus with the same reasoning as above we conclude that
$$
A_{\wedge}-\lambda:\Dom_{\wedge}\to x^{-m/2}L_b^2(Y^{\wedge};E)
$$
is invertible if and only if the restriction of the induced operator
$$
F_{\wedge}(\lambda) = T_{\wedge}(\lambda)(A_{\wedge} - \lambda) : 
\tilde\Sing_{\wedge,\max}\to \C^{d''}
$$
to $\tilde\Sing_{\wedge}=\Dom_{\wedge}/\Dom_{\wedge,\min}$ is invertible. 
Let $q_{\wedge}: \Dom_{\wedge,\max} \to \tilde\Sing_{\wedge,\max}$ be the 
canonical projection. From the relations
\begin{align*}
F_{\wedge}(\lambda)^{-1} = q_{\wedge}(A_{\wedge,\Dom_{\wedge}}-\lambda)^{-1}
K_{\wedge}(\lambda) &: \C^{d''} \to \tilde\Sing_{\wedge,\max}, \\
q_{\wedge}(A_{\wedge,\Dom_{\wedge}}-\lambda)^{-1} = F_{\wedge}(\lambda)^{-1}
T_{\wedge}(\lambda) &: x^{-m/2}L_b^2 \to \tilde\Sing_{\wedge,\max},
\end{align*}
and Proposition \ref{PropSmaxCond}, we deduce that our assumption on
$A_\wedge$ is equivalent to
\begin{equation}\label{Fdachestimate}
\|\kappa^{-1}_{|\lambda|^{1/m}}F_{\wedge}(\lambda)^{-1}\|_{\L(C^{d''},
\tilde\Sing_{\wedge,\max})} = O(1)
\;\text{ as } |\lambda| \to \infty.
\end{equation}
Note that $\|K_{\wedge}(\lambda)\| = O(|\lambda|)$ and
$\|T_{\wedge}(\lambda)\| = O(|\lambda|^{-1})$ as $|\lambda| \to \infty$ 
when considered as operators $\C^{d''} \to x^{-m/2}L^2_b$ and $x^{-m/2}L^2_b 
\to\C^{d''}$, respectively.
 
Write the operator $F(\lambda)\theta^{-1}
F_{\wedge}(\lambda)^{-1} : \C^{d''} \to \C^{d''}$ as
$$
F(\lambda)\theta^{-1}F_{\wedge}(\lambda)^{-1} = 
1 +\bigl(F(\lambda) - F_{\wedge}(\lambda)\theta\bigr)
\tilde{\kappa}_{|\lambda|^{1/m}} \theta^{-1}\kappa^{-1}_{|\lambda|^{1/m}}
F_{\wedge}(\lambda)^{-1},
$$
and let 
\begin{equation*}
R(\lambda)=\bigl(F(\lambda) - F_{\wedge}(\lambda)\theta\bigr)
\tilde{\kappa}_{|\lambda|^{1/m}} \theta^{-1}\kappa^{-1}_{|\lambda|^{1/m}}
F_{\wedge}(\lambda)^{-1}.
\end{equation*}
We will prove in Lemma \ref{FFdach} that
$$
\|(F(\lambda) - F_{\wedge}(\lambda)\theta)
\tilde{\kappa}_{|\lambda|^{1/m}}\|_{\L(\tilde\Sing_{\max},\C^{d''})} \to 0
\;\text{ as } |\lambda| \to \infty.
$$
Thus together with \eqref{Fdachestimate} we obtain that
$\|R(\lambda)\| \to 0$ as $|\lambda| \to \infty$.
Hence $1 + R(\lambda)$ is invertible for large $|\lambda| > 0$, and the
inverse is of the form $1 + \tilde{R}(\lambda)$ with $\|\tilde{R}(\lambda)\|
\to 0$ as $|\lambda| \to \infty$.
This shows that $F(\lambda) : \tilde\Sing \to \C^{d''}$ is invertible from 
the right for large $\lambda$, and by \eqref{Fdachestimate} the right-inverse 
$\theta^{-1}F_{\wedge}(\lambda)^{-1}(1 + \tilde{R}(\lambda))$ satisfies the 
norm estimate \eqref{Festimate}. Since
$$
\dim \tilde\Sing = \dim \tilde\Sing_{\wedge} = d'',
$$
we conclude that $F(\lambda)$ is also injective, and
so the invertibility of $F(\lambda)$ is proved.
In particular, the operator
$$
A_{\Dom} - \lambda : \Dom \to x^{-m/2}L_b^2(M;E)
$$
is invertible for large $\lambda$. It remains to show the norm estimates
\eqref{ResolventEstimates}.

In order to prove \eqref{ResolventEstimates} we make use of the family 
$\tilde{K}(\varrho)$ from Lemma~\ref{kappaSchlangelift} and the representation
\eqref{TheResolvent} of the resolvent. Thus we may write
\begin{align*}
(A_{\Dom} - \lambda)^{-1} &= B(\lambda) + (1 - B(\lambda)(A-\lambda)) 
\tilde{K}(|\lambda|^{1/m})\tilde{\kappa}^{-1}_{|\lambda|^{1/m}}
F(\lambda)^{-1}T(\lambda) \\
&=B(\lambda) + \tilde{K}(|\lambda|^{1/m})\tilde{\kappa}^{-1}_{|\lambda|^{1/m}}
F(\lambda)^{-1}T(\lambda) \\
&\hspace*{6em} - B(\lambda)(A-\lambda)\tilde{K}(|\lambda|^{1/m})
\tilde{\kappa}^{-1}_{|\lambda|^{1/m}}F(\lambda)^{-1}T(\lambda).
\end{align*}
By Remark \ref{PZweiEigenschaften} we have 
$\|B(\lambda)\|_{\L(x^{-m/2}L^2_b,\Dom_{\max})} = O(1)$ as 
$|\lambda| \to \infty$. In view of $\|T(\lambda)\|_{\L(x^{-m/2}L_b^2,\C^{d''})} 
= O(|\lambda|^{-1})$ and \eqref{Festimate} we further obtain
$$
\|\tilde{\kappa}^{-1}_{|\lambda|^{1/m}}F(\lambda)^{-1}
T(\lambda)\|_{\L(x^{-m/2}L_b^2,\tilde\Sing_{\max})}= O(|\lambda|^{-1})
\;\text{ as } |\lambda| \to \infty,
$$
and consequently, using \eqref{Kestimate2} we get
$$
\|\tilde{K}(|\lambda|^{1/m})\tilde{\kappa}^{-1}_{|\lambda|^{1/m}}
F(\lambda)^{-1}T(\lambda)\|_{\L(x^{-m/2}L_b^2,\Dom_{\max})}= O(1)
\;\text{ as } |\lambda| \to \infty.
$$
On the other hand, by \eqref{Festimate} and the estimates \eqref{Kestimate1} 
and \eqref{Kestimate2} we have
$$
\|(A-\lambda)\tilde{K}(|\lambda|^{1/m})\tilde{\kappa}^{-1}_{|\lambda|^{1/m}}
F(\lambda)^{-1}\|_{\L(\C^{d''},x^{-m/2}L_b^2)}= O(|\lambda|)
\;\text{ as } |\lambda| \to \infty.
$$
In view of $\|B(\lambda)\|_{\L(x^{-m/2}L^2_b,\Dom_{\max})} = O(1)$ and 
$\|T(\lambda)\|_{\L(x^{-m/2}L_b^2,\C^{d''})} = O(|\lambda|^{-1})$, 
we conclude that, as $|\lambda| \to \infty$,
$$
\|B(\lambda)(A-\lambda)\tilde{K}(|\lambda|^{1/m})
\tilde{\kappa}^{-1}_{|\lambda|^{1/m}}F(\lambda)^{-1}T(\lambda)\|_{
\L(x^{-m/2}L_b^2,\Dom_{\max})} = O(1).
$$

Summing up, we have proved
$$
\|(A_{\Dom} - \lambda)^{-1}\|_{\L(x^{-m/2}L_b^2,\Dom_{\max})} = O(1)
\;\text{ as } |\lambda| \to \infty,
$$
and the estimates \eqref{ResolventEstimates} follow.
\end{proof}

The following lemma completes the proof of Theorem \ref{RayMinimalGrowth}.

\begin{lemma}\label{FFdach}
With the notation of the proof of Theorem \ref{RayMinimalGrowth}, let
\begin{align*}
F(\lambda) = T(\lambda)(A-\lambda) &: \tilde\Sing_{\max} \to \C^{d''}, \\
F_{\wedge}(\lambda) = T_{\wedge}(\lambda)(A_{\wedge}-\lambda) &: 
\tilde\Sing_{\wedge,\max} \to \C^{d''}.
\end{align*}
Then
\begin{equation}\label{FFdachNull}
\|(F(\lambda) - F_{\wedge}(\lambda)\theta)
\tilde{\kappa}_{|\lambda|^{1/m}}\|_{\L(\tilde\Sing_{\max},\C^{d''})} \to 0
\;\text{ as } |\lambda| \to \infty.
\end{equation}
\end{lemma}
\begin{proof}
For proving \eqref{FFdachNull} it is sufficient to consider the restrictions
$$
(F(\lambda) - F_{\wedge}(\lambda)\theta)\tilde{\kappa}_{|\lambda|^{1/m}} :
\tilde\Sing_{\sigma_0} \to \C^{d''}
$$
for all $\sigma_0 \in \Sigma$. First of all, observe that	
\begin{gather*}
F(\lambda)\tilde{\kappa}_{|\lambda|^{1/m}} 
= T(\lambda)(A-\lambda)\tilde{K}(|\lambda|^{1/m}), \;\text{ and }\\
F_{\wedge}(\lambda)\theta \tilde{\kappa}_{|\lambda|^{1/m}} 
= F_{\wedge}(\lambda)\kappa_{|\lambda|^{1/m}}\theta
= T_\wedge(\lambda)(A_\wedge-\lambda)\kappa_{|\lambda|^{1/m}}\omega\theta
\end{gather*}
with the operator family $\tilde{K}(\varrho)=\omega(\varrho x)
\tilde\kappa_{\varrho}$ from Lemma \ref{kappaSchlangelift}. 
If $\omega_0 \in C_0^{\infty}([0,1))$ is a cut-off function near 
zero with $\omega \prec \omega_0$, then
\begin{align*}
(F(\lambda) - F_{\wedge}(\lambda)\theta)\tilde{\kappa}_{|\lambda|^{1/m}}
&= T(\lambda)(A-\lambda)\tilde{K}(|\lambda|^{1/m})
-  T_\wedge(\lambda)(A_\wedge-\lambda)\kappa_{|\lambda|^{1/m}}\omega\theta\\
&= T(\lambda)\omega_0(A-\lambda)\tilde{K}(|\lambda|^{1/m})
-  T_\wedge(\lambda)\omega_0(A_\wedge-\lambda)\kappa_{|\lambda|^{1/m}}
   \omega\theta\\
&= T(\lambda)\omega_0\left((A-\lambda)\tilde K(|\lambda|^{1/m})
-(A_\wedge-\lambda)\kappa_{|\lambda|^{1/m}}\omega\theta\right)\\
&\hspace*{4em} +(T(\lambda)-T_\wedge(\lambda))\omega_0(A_\wedge-\lambda)
\kappa_{|\lambda|^{1/m}}\omega\theta\\
&= T(\lambda)\omega_0\left(A\tilde K(|\lambda|^{1/m})
 -A_\wedge\kappa_{|\lambda|^{1/m}}\omega\theta\right)\\
&\hspace*{4em} - T(\lambda)\omega_0\lambda \left(\tilde K(|\lambda|^{1/m})
 -\kappa_{|\lambda|^{1/m}}\omega\theta\right)\\
&\hspace*{8em} +(T(\lambda)-T_\wedge(\lambda))\omega_0(A_\wedge-\lambda)
\kappa_{|\lambda|^{1/m}}\omega\theta
\end{align*}

By \eqref{PhiAPhiL}, \eqref{Theta0Estimate}, and Lemma~\ref{Ltheta} 
it follows that the norms of
\begin{gather*}
A\tilde{K}(|\lambda|^{1/m}) - A_{\wedge}\kappa_{|\lambda|^{1/m}}
\omega \theta = A\tilde{K}(|\lambda|^{1/m}) - |\lambda| 
\kappa_{|\lambda|^{1/m}}A_{\wedge}\omega \theta, \\
\lambda\bigl(\tilde{K}(|\lambda|^{1/m}) - \kappa_{|\lambda|^{1/m}} \omega 
\theta\bigr) = 
\lambda \kappa_{|\lambda|^{1/m}} \omega (L_{|\lambda|^{1/m}} - \theta)
\end{gather*}
in $\L(\tilde\Sing_{\sigma_0},x^{-m/2}L^2_b)$ are 
$O(|\lambda|^{1-1/m}\log^\mu |\lambda|)$ as $|\lambda| \to \infty$. 
Finally, because of the norm estimates $\|T(\lambda)\omega_0\|=O(|\lambda|^{-1})$, $\|(A_{\wedge}-\lambda)\kappa_{|\lambda|^{1/m}} \omega \theta\| = 
O(|\lambda|)$, and also $\|\bigl(T(\lambda)- T_{\wedge}(\lambda)\bigr)\omega_0 \| = O(|\lambda|^{-1-1/m})$ as $|\lambda|\to\infty$, the lemma follows.
\end{proof}

Finally, we want to point out that under the assumptions of Theorem
\ref{RayMinimalGrowth} we get the existence of the resolvent with polynomial
bounds for the norm also for closed extensions in Sobolev spaces of arbitrary smoothness.

\begin{theorem}\label{RMG-HigherRegularity} 
Let $A\in x^{-m}\Diff_b^m(M;E)$ be $c$-elliptic with parameter in 
$\Lambda\subset \C$, and let $\Dom^s \subset x^{-m/2}H^{s}_b(M;E)$
be a domain such that $A_{\Dom^s}$ is closed.
Assume that $\Lambda$ is a sector of minimal growth for the closed extension
$A_{\wedge,\Dom^0_\wedge}$ of $A_\wedge$ in $x^{-m/2}L^2_b$,
where $\Dom^0_\wedge \subset x^{-m/2}L^2_b$ is the domain
associated to $\Dom^0$ according to \eqref{DwedgeIdentitaet}.
Then for $\lambda\in\Lambda$ sufficiently large,
\begin{equation*}
 A_{\Dom^s}-\lambda:\Dom^s\to x^{-m/2}H^{s}_b(M;E)
\end{equation*}
is invertible and the resolvent satisfies the equation
\begin{equation*}
 (A_{\Dom^s}-\lambda)^{-1} = B(\lambda) + 
 (A_{\Dom^s}-\lambda)^{-1}\Pi(\lambda)
\end{equation*}
with the parametrix $B(\lambda)$ and the projection $\Pi(\lambda)$
from Theorem~\ref{PDrei}. Moreover, for every $s\in\R$ there exists
$M(s)\in\R$ such that, as $|\lambda|\to\infty$, 
\begin{equation}\label{HsResolventEstimates}
\begin{split} 
 &\bigl\|(A_{\Dom^s}-\lambda)^{-1}\bigr\|_{\L(x^{-m/2}H^s_b)}
 = O(|\lambda|^{M(s)-1}), \\
 &\bigl\|(A_{\Dom^s}-\lambda)^{-1}\bigr\|_{\L(x^{-m/2}H^s_b,
 \Dom^s_{\max})} = O(|\lambda|^{M(s)}). 
\end{split}
\end{equation}
\end{theorem}
\begin{proof}
We know from Proposition~\ref{SpectrumIndependent} that the spectrum does not depend on the regularity $s \in \R$. Consequently, from Theorem 
\ref{RayMinimalGrowth} we obtain the existence of the resolvent 
$(A_{\Dom^s}-\lambda)^{-1}$ for large $\lambda$.

Moreover, as in the proof of Theorem \ref{RayMinimalGrowth} we may write
$$
(A_{\Dom^s} - \lambda)^{-1} = B(\lambda) + (1 - B(\lambda)(A-\lambda))
\tilde{K}(|\lambda|^{1/m})\tilde{\kappa}^{-1}_{|\lambda|^{1/m}}
F(\lambda)^{-1}T(\lambda).
$$
According to what we have proved in this and the previous section we obtain 
that the norms of all operators
\begin{align*}
B(\lambda) &: x^{-m/2}H^s_b \to \Dom^s_{\max}, \\
T(\lambda) &: x^{-m/2}H^s_b \to \C^{d''}, \\
\tilde{\kappa}^{-1}_{|\lambda|^{1/m}}F(\lambda)^{-1} &: \C^{d''} \to 
\tilde\Sing_{\max}, \\
\tilde{K}(|\lambda|^{1/m}) &: \tilde\Sing_{\max} \to \Dom_{\max}^s, \\
(1 - B(\lambda)(A-\lambda)) &: \Dom^s_{\max} \to \Dom^s_{\max}
\end{align*}
behave polynomially as $|\lambda| \to \infty$. This proves the theorem.
\end{proof}


\begin{appendix}
\renewcommand{\theequation}{\Alph{section}.\arabic{equation}}

\section{Invertibility of Fredholm families}
\label{sec-ExtraConditions}

The theorem of this section is essential for the existence of
extra conditions in order to make the family $A_{\wedge} - \lambda$
invertible on the model cone $Y^\wedge$.  The main application of
Theorem~\ref{FredFam} concerns the Fredholm family
$$
a(\lambda) = A_{\wedge} - \lambda : \Dom_{\min}(A_{\wedge}) \to
x^{-m/2}L^2_b(Y^{\wedge};E),
$$
where $\lambda \in \Omega = \{z\in\Lambda : |z|=1\}$ (see also
Corollary~\ref{ExtraCondSector}). 

Theorem~\ref{FredFam} is rather standard and widely
used throughout the literature. However, since several of our key
arguments in the parametrix construction given in
Theorem~\ref{AminuslambdaplusK} rely on this result, we decide to 
give here an independent proof.

\begin{theorem}\label{FredFam}
Let $\Omega$ be a compact connected space
$(C^{\infty}$-manifold$)$, and let $a : \Omega \to
\L(H_1,H_2)$ be a continuous $($smooth$)$ Fredholm family
in the Hilbert bundles $H_1$ and $H_2$. Then there exist
$($smooth$)$ vector bundles $J_-,\, J_+ \in \mathrm{Vect}(\Omega)$
and continuous $($smooth$)$ sections $t,k,q$ such that
$$
\begin{pmatrix}
a & k \\ t & q
\end{pmatrix}
: \Omega \to \L\left(
\begin{array}{c}H_1 \\ \oplus \\ J_-\end{array} ,
\begin{array}{c}H_2 \\ \oplus \\ J_+\end{array}
\right)
$$
is a family of isomorphisms. The difference $[J_+] - [J_-] \in
K(\Omega)$ equals the index $\Ind_K(a)$ of $a$. If $a$ is onto or
one-to-one, we can choose $J_- = 0$ or $J_+ = 0$, respectively. If
$\Omega$ is contractible, then we have $J_{\pm} = {\C}^{N_{\pm}}$
with $N_{\pm} \in {\N}_0$.
\end{theorem}
\begin{proof}
Let $x \in \Omega$ be arbitrary. Choose (smooth) sections
$s_1,\ldots,s_{N(x)}$ of $H_2$ such that
$\{s_1(x),\ldots,s_{N(x)}(x)\}$ forms a basis of a complement of
$\rg(a(x))$ in $H_2$. Define
$$
k_x : \Omega \to \L({\C}^{N(x)}, H_2), \quad
\begin{pmatrix}
c_1 \\ \vdots \\ c_{N(x)}
\end{pmatrix} \mapsto \sum\limits_{j=1}^{N(x)}c_js_j.
$$
It follows that
$$
\begin{pmatrix}
a(x) & k_x(x)
\end{pmatrix}
: \begin{array}{c} H_1 \\ \oplus \\ {\C}^{N(x)} \end{array} 
\to H_2
$$
is surjective, and so $\begin{pmatrix} a & k_x \end{pmatrix}$ is
surjective in an open neighborhood $U(x)\subset \Omega$.
Let $\Omega = \bigcup_{k=1}^MU(x_k)$ be a covering of
$\Omega$ by such neighborhoods, and set
$$
k= \begin{pmatrix} k_{x_1} \ldots k_{x_M} \end{pmatrix} :
\Omega \to \L\left(
\bigoplus\limits_{k=1}^M {\C}^{N(x_k)} , H_2\right).
$$
Then
$$
\begin{pmatrix} a(x) & k(x) \end{pmatrix} :
\begin{array}{c} H_1 \\ \oplus \\ {\C}^{N_-} \end{array}
\to H_2
$$
is surjective for all $x \in \Omega$, where $N_-=
\sum_{k=1}^MN(x_k)$. 

So suppose without loss of generality that $a(x)$ is a surjective 
Fredholm family. Then 
$\dim\ker a(x)$ is independent of $x$, and the disjoint union
\begin{equation*}
J_+=\bigsqcup_{x\in \Omega} \ker a(x)
\end{equation*}
is a locally trivial finite rank continuous (smooth) vector bundle. Let $\pi_x:H_1\to J_+$ be the orthogonal projection. Then
\begin{equation*}
\begin{pmatrix}
a\\\pi
\end{pmatrix}:
H_1\to
\begin{matrix}
H_2\\\oplus\\J_+
\end{matrix}
\end{equation*}
is invertible.

If $a$ is pointwise injective, we obtain from the above argument,
applied to $a^*$, that we may choose $J_+ = 0$. This finishes the
proof of the theorem.
\end{proof}

\begin{remark}\label{DenseSubspace}
Let $H_1, H_2$ be Hilbert spaces and
$$
\begin{pmatrix} a & k \\ t & q \end{pmatrix} : 
\Omega \to \L\left(
\begin{array}{c} H_1 \\ \oplus \\ {\C}^{N_-} \end{array},
\begin{array}{c} H_2 \\ \oplus \\ {\C}^{N_+} \end{array}
\right)
$$
be a smooth family of isomorphisms as in Theorem~\ref{FredFam}.
Moreover, let $D_1' \subset H_1'$ and $D_2 \subset H_2$ be
dense subspaces. Then we can modify $t$ and $k$ such that 
\begin{equation*}
k \in C^{\infty}(\Omega) \otimes ({\C}^{N_-})^* \otimes D_2 
\;\text{ and }\; 
t \in C^{\infty}(\Omega) \otimes D_1' \otimes {\C}^{N_+}.
\end{equation*}
\end{remark}

\begin{corollary}\label{ExtraCondSector}
Let $\Lambda$ be a closed sector in $\C$ as defined in 
Section~\ref{sec-Parametrix}. 
Let $H_1$ and $H_2$ be Hilbert spaces with strongly continuous 
groups $\{\kappa_{\varrho}\}_{\varrho\in\R_+}$ and 
$\{\tilde{\kappa}_{\varrho}\}_{\varrho\in\R_+}$,
and let $a \in C^{\infty}(\Lambda \minus \{0\},\L(H_1,H_2))$ be a 
Fredholm family that satisfies
$$
 a(\varrho^d\lambda) = \varrho^{\mu}\tilde{\kappa}_{\varrho}
 a(\lambda)\kappa_{\varrho}^{-1}
$$
for every $\varrho > 0$, where $d \in \N_0$ and $\mu \in {\R}$ are
given numbers. Then there exist $t$, $k$, and $q$ such that
\begin{equation*}
\begin{pmatrix} a & k \\ t & q \end{pmatrix} \in
C^{\infty}\left(\Lambda \minus \{0\}, \L\left(
\begin{array}{c} H_1 \\ \oplus \\ {\C}^{N_-} \end{array},
\begin{array}{c} H_2 \\ \oplus \\ {\C}^{N_+} \end{array}
\right)\right)
\end{equation*}
is pointwise an isomorphism, and it satisfies
$$
\begin{pmatrix} a(\varrho^d\lambda) & k(\varrho^d\lambda) \\
t(\varrho^d\lambda) & q(\varrho^d\lambda) \end{pmatrix} =
\varrho^{\mu}
\begin{pmatrix} \tilde{\kappa}_{\varrho} & 0 \\ 0 & 1 \end{pmatrix}
\begin{pmatrix} a(\lambda) & k(\lambda) \\ t(\lambda) & q(\lambda)
\end{pmatrix}
\begin{pmatrix} \kappa_{\varrho}^{-1} & 0 \\ 0 & 1 \end{pmatrix}
$$
for every $\varrho > 0$. If $a$ is onto or one-to-one, then we may
choose $N_- = 0$ or $N_+ = 0$, respectively.
\end{corollary}
\begin{proof}
Let $\Omega=\{z\in\Lambda : |z|=1\}$ and let $\hat{a}= a|_\Omega$.
According to Theorem~\ref{FredFam} there exist $\hat{t}$,
$\hat{k}$, and $\hat{q}$ such that the operator function
$$
\begin{pmatrix} \hat{a} & \hat{k} \\ \hat{t} & \hat{q} \end{pmatrix}
\in C^{\infty}\left(\Omega,\L\left(
\begin{array}{c} H_1 \\ \oplus \\ {\C}^{N_-} \end{array},
\begin{array}{c} H_2 \\ \oplus \\ {\C}^{N_+} \end{array}
\right)\right)
$$
is pointwise bijective, and we may choose $N_- = 0$ or $N_+ = 0$
provided that $a$ is everywhere surjective or injective,
respectively. We will be done if we can show that the extension
by $\kappa$-homogeneity
\begin{equation}\label{ExtensionHomogen}
\begin{pmatrix} a(\lambda) & k(\lambda) \\ t(\lambda) &
q(\lambda) \end{pmatrix} =
|\lambda|^{\frac{\mu}{d}}\begin{pmatrix}
\tilde{\kappa}_{|\lambda|^{1/{d}}} & 0 \\ 0 & 1
\end{pmatrix}
\begin{pmatrix} \hat{a}\big(\frac{\lambda}{|\lambda|}\big) &
\hat{k}\big(\frac{\lambda}{|\lambda|}\big) \\
\hat{t}\big(\frac{\lambda}{|\lambda|}\big) &
\hat{q}\big(\frac{\lambda}{|\lambda|}\big) \end{pmatrix}
\begin{pmatrix} \kappa_{|\lambda|^{1/{d}}}^{-1} & 0 \\ 0 & 1
\end{pmatrix}
\end{equation}
for $\lambda \in \Lambda \minus \{0\}$ depends smoothly on
$\lambda$; note that the group actions are assumed to be only
strongly continuous.

In fact, $q$ is clearly $C^{\infty}$ and $a$ was assumed to be
smooth. Thus we only have to check the smoothness of $t$ and $k$.
According to Remark~\ref{DenseSubspace} we may take $\hat{k} \in
C^{\infty}(\Omega) \otimes ({\C}^{N_-})^* \otimes D_2$ and
$\hat{t} \in C^{\infty}(\Omega) \otimes D_1' \otimes {\C}^{N_+}$,
where $D_1' \subset H_1'$ is the space of $C^{\infty}$-elements of
the dual group action $\{\kappa_{\varrho}'\}$ on $H_1'$, and $D_2$
is the space of $C^{\infty}$-elements of the group action
$\{\tilde{\kappa}_{\varrho}\}$ on $H_2$.  With these choices the
operator function defined in \eqref{ExtensionHomogen} is smooth,
as desired.
\end{proof}

\begin{remark}
In our applications the group action involved is always the dilation
group defined in \eqref{DilationGroup}. The space of compactly
supported smooth functions is then an admissible choice for the spaces
$D_1'$ and $D_2$ in the proof of Corollary~\ref{ExtraCondSector} 
(see also Remark~\ref{DenseSubspace}).
\end{remark}

\end{appendix}


\end{document}